\documentclass[12pt]{amsart}
\usepackage{amsbsy}
\usepackage{graphicx,epsfig,subfigure}
\begin{document}

%

%
\newtheorem{theorem}{Theorem}
\newtheorem{proposition}[theorem]{Proposition}
\newtheorem{lemma}[theorem]{Lemma}
\newtheorem{corollary}[theorem]{Corollary}
\newtheorem{definition}[theorem]{Definition}
\newtheorem{remark}[theorem]{Remark}
\numberwithin{equation}{section} \numberwithin{theorem}{section}
\newcommand{\be}{\begin{equation}}
\newcommand{\ee}{\end{equation}}
\newcommand{\re}{{\mathbb R}}
\newcommand{\n}{\nabla}
\newcommand{\ren}{{\mathbb R}^N}
\newcommand{\iy}{\infty}
\newcommand{\pa}{\partial}
\newcommand{\ms}{\medskip\vskip-.1cm}
\newcommand{\mpb}{\medskip}
\newcommand{\ssk}{\smallskip}
\newcommand{\BB}{{\bf B}}
\newcommand{\Am}{{\bf A}_{2m}}
\newcommand{\bL}{\BB^*}
\newcommand{\bLs}{\BB}
\renewcommand{\a}{\alpha}
\renewcommand{\b}{\beta}
\newcommand{\g}{\gamma}
\newcommand{\ka}{\kappa}
\newcommand{\G}{\Gamma}
\renewcommand{\d}{\delta}
\newcommand{\D}{\Delta}
\newcommand{\e}{\varepsilon}
\newcommand{\vp}{\varphi}
\renewcommand{\l}{\lambda}
\renewcommand{\o}{\omega}
\renewcommand{\O}{\Omega}
\newcommand{\s}{\sigma}
\renewcommand{\t}{\tau}
\renewcommand{\th}{\theta}
\newcommand{\z}{\zeta}
\newcommand{\wx}{\widetilde x}
\newcommand{\wt}{\widetilde t}
\newcommand{\noi}{\noindent}
\newcommand{\lb}{\left (}
\newcommand{\rb}{\right )}
\newcommand{\lsb}{\left [}
\newcommand{\rsb}{\right ]}
\newcommand{\lab}{\left \langle}
\newcommand{\rab}{\right \rangle }
\newcommand{\gap}{\vskip .5cm}
\newcommand{\bz}{\bar{z}}
\newcommand{\bg}{\bar{g}}
\newcommand{\Ba}{\bar{a}}
\newcommand{\bt}{\bar{\th}}
\def\com#1{\fbox{\parbox{6in}{\texttt{#1}}}}


\title{\bf Three types of  self-similar blow-up for  the fourth-order
${\mathbf p}$-Laplacian
equation with source: variational and branching approaches}

\author
{V.A. Galaktionov} 

\address{ 
Department of Math. Sci., University of Bath,
 Bath, BA2 7AY, UK}
\email{vag@maths.bath.ac.uk}


\thanks{Research supported by  RTN network
HPRN-CT-2002-00274}
 \keywords{Fourth-order quasilinear
parabolic equation finite propagation, similarity solutions,
blow-up}
 \subjclass{35K55, 35K40, 35K65}
\date{\today}



\begin{abstract}

Self-similar blow-up behaviour  for the fourth-order quasilinear
$p$-Laplacian equation with source,
 $$
 u_t = - (|u_{xx}|^n u_{xx})_{xx} + |u|^{p-1}u  \quad \mbox{in} \,\,\, \re \times \re_+,
  \quad \mbox{where} \,\,\, n>0, \,\,\, p>1,
 $$
   is studied. Using  variational setting for $p=n+1$ and
 branching techniques for $p \not = n+1$,
 finite and countable families  of blow-up patterns of the self-similar form
 $$
  \mbox{$
 u_S(x,t)=(T-t)^{- \frac 1{p-1}}f(y), \quad  \mbox{where} \quad  y=x/(T-t)^{\b},
  \,\,\, \b =- \frac {p-(n+1)}{2(n+2)(p-1)},
   $}
 $$
 are described by an analytic-numerical approach.
Three parameter ranges: $p=n+1$ (regional), $p>n+1$ (single
point), and $1<p<n+1$ (global blow-up) are studied. This blow-up
model
   is motivated by
 the second-order reaction diffusion counterpart
 $$
 u_t = (|u_x|^n u_x)_x+ u^p  \quad (u \ge 0)
 $$
 that was studied in the middle of the 1980s, while first results on blow-up of solutions
 were established by Tsutsumi in 1972.

 This paper is an earlier extended preprint of \cite{GalpCAM}.

\end{abstract}


\maketitle

\section{Introduction: classic and recent blow-up  reaction-diffusion models}

\subsection{Classic second-order model}

The nonlinear $p$-Laplacian operator in $\ren$,
  \be
  \label{p1}
   \D_p u \equiv \n \cdot (|\n u|^{p-2} \n u),
   \quad \mbox{with exponents} \quad p>1 \quad \big(\n={\rm grad}_x \big),
    \ee
    which serves as a natural
  extension of the
Laplacian
 $$
 \D= \D_2, \quad \mbox{i.e.,
  for $p=2$},
   $$
  enters many classic PDEs of
mathematical physics.  One of the key mathematical advantages of
the $p$-Laplacian (\ref{p1}) is that it is {\em nonlinear} and at
the same time remains a {\em monotone} operator in the
$L^2$-metric precisely as the linear Laplacian $\D$ does.
Operators such as (\ref{p1}) appear in many works on nonlinear
parabolic or elliptic PDEs since the 1950s; see various examples,
references, and applications in Lions' classic book \cite{LIO}.
Gradient-dependent nonlinear operators are typical for filtration,
combustion (solid fuels), and non-Newtonian (dilatable,
pseudo-plastic fluids) liquids theory; see
\cite[p.~428]{GSVR}.

\ssk

Concerning parabolic PDEs admitting {\em blow-up solutions} that
compose  the main subject of the present paper, the $p$-Laplacian
also appeared before other well-known nowadays {\em porous medium}
type nonlinearities (see equation (\ref{pm21}) below).
 Namely, it is
remarkable that the first results on blow-up in {\em quasilinear}
parabolic equations were obtained
  by
Tsutsumi in 1972 \cite{Tsut74}
 for the {\em second-order $p$-Laplacian equation} ($p$LE--2) {\em with
 source}
  posed in a bounded smooth domain $\O \subset \ren$ with the zero Dirichlet boundary condition:
 \be
 \label{sss1}
 u_t = \n \cdot\big(|\n u|^n \n u  \big)+ u^p \quad \mbox{in}
 \quad \Omega \times \re_+  \quad (u \ge 0),
  \ee
   where, in comparison with (\ref{p1}),
    we have renamed the exponents by setting $n=p-2>0$ and
    write the source term as $u^p$.
        Concerning the structure of blow-up
   singularities,
  various countable and finite
families of self-similar blow-up
 patterns for the
 one-dimensional equation  of (\ref{sss1}) (and also for the radially symmetric version
 of (\ref{sss1})),
  \be
  \label{s1}
 u_t = \big(|u_x|^n u_x \big)_x+ u^p ,
  \ee
 have been known since the middle of the 1980s; see
  \cite{GPos0, GPos, BuGa}, and other related  references therein.
Surprisingly for the author who initiated the study in
\cite{GPos0, GPos}, it turned out that (\ref{s1}) generates much
wider countable and even uncountable families of self-similar
blow-up patterns than the {\em porous medium equation with source}
({PME with source})
 \be
 \label{pm21}
 u_t = (u^{n+1})_{xx} + u^p \quad (u \ge 0),
 \ee
which was studied by Kurdyumov's Russian School on blow-up and
localization since the beginning of the 1970s; see history,
references, and basic results in \cite[Ch.~4]{SGKM}.

 It is worth
mentioning that
 the set of blow-up similarity
solutions of (\ref{s1}), to say nothing about non-radial patterns
for (\ref{sss1}), is rather complicated (e.g., contains infinite
countable and even uncountable  families of positive solutions for
$p>n+1$), so there are still some difficult open mathematical
problems concerning the structure of blow-up singularities for
(\ref{s1}).

Blow-up results for the $p$-Laplacian equations with source
(\ref{s1}) and (\ref{sss1})
 together with Fujita's pioneering  study of the {\em semilinear heat equation}
  (1966) \cite{Fuj66},
 \be
 \label{Fu1}
 u_t = \D u + u^p \quad (p>1),
  \ee
are crucial for modern singularity and blow-up theory of nonlinear
evolution PDEs.

\smallskip

Nowadays, blow-up and other singularity formation phenomena for
various classes of nonlinear evolution PDEs are rather popular in
mathematical literature and applications in mechanics and physics.
 It is well established that   blow-up
phenomena in nonlinear PDEs not only present principal evolution
patterns of interest in application, but also can give insight
into the deep mathematical nature of nonlinear equations under
consideration and describe general aspects of various fundamental
problems of existence-nonexistence, uniqueness-nonuniqueness,
optimal regularity classes, and admissible asymptotics of proper
solutions.

To emphasize that this is not an exaggeration, let us  mention
that, according to the typical tools of the possible and already
available analysis and proofs,
that the
  two
key open PDE/geometry problems of the twentieth
  and twenty-first
      century  are directly attributed to the area of PDE  blow-up research:

\ssk

 {\bf Problem (I):} {\em Poincar\'e Conjecture}
(a closed connected 3D manifold is homeomorphic to ${\bf S}^3$)
 and the general geometrization problem
  with Perel'man's
recent proof   by introducing  two  new monotonicity formulae and
others to pass through blow-up singularities of Ricci flows with
surgery (see \cite{CZ06} for a full account of
 history, references, and recent development); and

 \ssk

 {\bf Problem (II):} {\em Uniqueness} or {\em nonuniqueness} (and
hence {\em nonexistence} or {\em existence}  of local small-scale
blow-up  singularities) in the 3D Navier--Stokes
equations\footnote{See \cite{GalNSE} as a most recent survey,
where connections with reaction-diffusion theory are discussed.}
(one of the Millennium Prize Problems for the Clay Institute; see
Fefferman
 \cite{Feff00}).



\ssk

 There are several monographs
   \cite{BebEb, SGKM, Pao, MitPoh, AMGV, GalGeom, QSupl}, which
   are
   devoted mainly to space-time stricture of blow-up singularities in second-order
reaction-diffusion PDEs and explain  the role of blow-up phenomena
in general PDE theory. See also \cite{GSVR} presenting
  various exact solutions and some examples of
partial singularity analysis of other classes of thin film,
nonlinear dispersion, and hyperbolic PDEs. In the monograph
\cite{MitPoh}, a nonlinear capacity approach was shown to be
efficient to detect conditions of global nonexistence for a
variety of nonlinear PDEs and systems of different orders and
types.

The questions of the  space-time structure and multiplicity of
possible blow-up asymptotics represent problems of higher
complexity that need another more involved  mathematical
treatment, which often and still cannot be fully justified
rigorously, so a true combination of various approaches including
enhanced numerics is
 in great demand.

\subsection{Fourth-order reaction-diffusion equation}

In this  paper, we study self-similar blow-up for the following
quasilinear parabolic {\em fourth-order  $p$-Laplacian equation
with source} ($p$LE--4 with source):
 \be
 \label{1.5}
  \fbox{$
   \ssk\ssk\ssk
u_t = {\bf A}(u) \equiv - \big(|u_{xx}|^n u_{xx}\big)_{xx} +
|u|^{p-1}u \quad \mbox{in} \,\,\, \re \times \re_+,
 \ssk\ssk\ssk
 $}
 \ee
 where, as above, $n >0$ and $p>1$.
 Here, similar to (\ref{p1}), the fourth-order $p$-Laplacian
 operator, where we set $p=n+2>1$ ($N=1$ in (\ref{1.5})),
  $$
  \D_{p,2} \, u=- \D \big(|\D u|^{p-2} \D u \big)
   $$
   is monotone in the metric of $L^2(\ren)$.
 For $n=0$, (\ref{1.5}) reduces to the semilinear
 equation
 \be
 \label{1.3}
 u_t = -u_{xxxx} +|u|^{p-1}u,
  \ee
 which  describes {\em single point
 blow-up} only for all $p>1$ and is already known to admit various similarity and other blow-up
  solutions, \cite{BGW1}.
  Moreover, it is curious that we have found quite
 fruitful to use the analogy with the linear {\em bi-harmonic
 equation}
  \be
   \label{bh1}
   u_t= - u_{xxxx} + u \quad \mbox{in} \quad \re\times \re_+,
    \ee
    which is obtained from (\ref{1.5}) by both limits
 $n \to 0$ and $p\to 1$. A simple countable subset of exponential
 patterns for (\ref{bh1}) is easy to describe on the basis of
 spectral theory presented in Section \ref{S3.Sp}. Eventually, we will
  detect  certain traces of such countable sets
 (the so-called $p$-branches) of similarity solutions in the nonlinear problem
 (\ref{1.5}).

Being involved in the mathematical study of blow-up for the PME
with source (\ref{pm21}) from the middle of 1970s and for the
$p$LE-2 with source (\ref{sss1}) from the 1980s, the author must
admit that
 the study of blow-up patterns for the proposed
$p$LE-4 with source (\ref{1.5}) was quite a challenge and the
author did not expect that the necessary mathematics should be so
dramatically changed to cover approximately the same concepts
developed twenty or even thirty years earlier. Recall that in
(\ref{1.5}) we just increase by two the order of the diffusion
operator in comparison with the standard model (\ref{s1}).
However, this makes almost all  mathematical tools applied before
very successfully to (\ref{s1}) almost  nonexistent.

\ssk

 Thus,
we consider for (\ref{1.5}) the Cauchy problem with given bounded
compactly supported data
 \be
 \label{u91}
 u(x,0)=u_0(x) \in C_0(\re).
  \ee
 Since the operator ${\bf A}$ in (\ref{1.5}) is potential in the metric of $L^2$
and the $p$-Laplacian is also a monotone operator there, local
existence and uniqueness of a unique weak (continuous) solution,
which is defined in the standard manner, are not principal issues
and follow from classic theory of monotone operators; see Lions
\cite[Ch.~2]{LIO}.
 Finite
propagation  phenomena for the PDE (\ref{1.5}) 
are proved by energy estimates via Saint--Venant's principle; see
\cite{Shi2}, references therein,  and a survey in \cite{GS1S-V}.
Therefore, there exists the unique local solution of the Cauchy
problem (\ref{1.5}), (\ref{u91}), which is a compactly supported
function $u(x,t)$ that can blow up in finite time in the sense
that
 \be
 \label{ub1}
  \mbox{$
 \sup_{x \in \re} \, |u(x,t)| \to + \infty
  \quad \mbox{as} \quad t \to T^-< \infty.
   $}
   \ee
Existence of blow-up in such higher-order quasilinear parabolic
equations is a reasonably well-understood phenomenon; see
references and approaches in \cite{Eg4, Gal2m, GPohTFE} and
Mitidieri--Pohozaev \cite{MitPoh}. For instance (see \cite{Eg4}
and references therein), it is known that, for the equation
similar to (\ref{1.5}) with the absolute value in the source-term,
 \be
 \label{sss1N}
 u_t =- \big(|u_{xx}|^n u_{xx}\big)_{xx} + |u|^p \quad \mbox{in} \quad
 \re \times \re_+,
  \ee
all nontrivial solutions with  data having positive first Fourier
coefficient,
 $$
  \mbox{$
 \int\limits_{\re} u_0(x) \, {\mathrm d}x>0,
 $}
 $$
 blow-up in finite time in the subcritical Fujita range
  $$
   \mbox{$
  n+1 < p < p_0=n +1 + \frac{2(n+2)}N \big|_{N=1}=3n+5 \quad (n \ge
  0),
   $}
   $$
   as well as, most probably, in the critical case $p=p_0$ , which
    needs
     additional study.


\subsection{Layout of the paper: three types of blow-up}

  In Section \ref{SectLocR}, we
describe some local and rather delicate oscillatory properties of
travelling wave solutions near finite interfaces. This is the
first time, where we face difficult and still non fully justified
mathematics concerning higher-order degenerate $p$-Laplacians.
Section \ref{Sect2} is devoted to the setting of blow-up
self-similar solutions and some preliminaries concerning the
linear operator with $n=0$ (even this issue is not that
straightforward and demands essentially non self-adjoint theory).

 Further principal difficulties and important  mathematical problems in the
study of such blow-up solutions concern
 the description and classification of possible types (the
 structure, stability, and multiplicity) of
 blow-up patterns occurring
in finite time. Later on,
 we study three classes of similarity blow-up
solutions of (\ref{1.5}) in the ranges:

\ssk

 {\bf (i)} Section \ref{SectS}: $p=n+1$, {\em regional blow-up}, so the infinite limit (\ref{ub1})
 occurs on a bounded $x$-interval;

\ssk

 {\bf (ii)} Section \ref{SectLS1}: $p>n+1$, {\em single point blow-up}, so (\ref{ub1})
 happens at a single point, say, at $x=0$, and then
 $u(x,T^-)$ is
 bounded for any $x \not = 0$; and

\ssk

{\bf (iii)} Section \ref{SectHS}: $p \in (1,n+1)$, {\em global
blow-up}, and (\ref{ub1}) happens for any $x \in \re$ (and
possibly uniformly on any bounded $x$-interval).

\ssk


A similar classification and various single point blow-up patterns
of the so-called  {\em P-, Q-,
 R-}, and {\em S-type}
for the second-order counterpart (\ref{s1})
  have been known since 1980s; see  \cite{GPos}
  and more references and results in \cite{BuGa}.
  Actually, we show that some concepts of the methodology developed in \cite{GPos, BuGa}
for (\ref{s1}) also apply to the fourth-order reaction-diffusion
equation (\ref{1.5}), but indeed demand a  different and more
difficult mathematics. Several problems remain open still.
 It turns out that, in general, the PDE (\ref{1.5}) admits more
 complicated sets of similarity patterns than
 the {\em fourth-order porous medium equation} (PME$-4$)
 {\em with source} \cite{GPME},
 \be
 \label{pp1}
 u_t= - (|u|^n u)_{xxxx} + |u|^{p-1}u.
  \ee
 The general scheme of blow-up study via variational and branching
approaches applies to  higher-order $p$-Laplacian PDEs such as the
$p$LE--6 with source (or any $2m$th-order one)
 $$
 u_t=\big(|u_{xxx}|^n u_{xxx}\big)_{xxx} + |u|^{p-1}u
 \quad \big(\mbox{or} \quad u_t=(-1)^{m+1} D_x^m(|D_x^m u|^n D_x^m u)+|u|^{p-1}u\big).
  $$

 \subsection{On some other higher-order PDEs with blow-up, extinction, and finite interfaces}

Blow-up  in parabolic PDEs with higher-order diffusion becomes
much more difficult than for second-order reaction-diffusion
equations.
 Even simpler PDEs such as
 the {\em extended Frank-Kamenetskii equation} in one dimension
 \be
 \label{1.2}
 u_t = -u_{xxxx} +{\mathrm e}^u \quad \mbox{or}
 \quad  u_t = u_{xxxxxx} +{\mathrm e}^u,
  \ee
  and their  counterparts with power nonlinearities
\be
 \label{1.3NN}
 u_t = -u_{xxxx} +|u|^{p-1}u \quad \mbox{or}
 \quad  u_t = u_{xxxxxx} +|u|^{p-1}u,
  \ee
revealed several principally new asymptotic blow-up properties
demanding novel mathematical approaches; see details in
  \cite{BGW1, Gal2m}.
  Similar difficulties occur for
    the {\em Semenov--Rayleigh--Benard
  problem} with the leading operator of the form
 \be
 \label{1.4}
 u_t= - u_{xxxx} + \b[(u_x)^3]_x + {\mathrm e}^u \quad (\b \ge 0);
 \ee
see \cite{GW1}. The mathematical difficulties in understanding the
ODE and PDE blow-up patterns  increase dramatically with the order
of differential diffusion-like operators in the equations.
Interesting {\em regional blow-up} and  oscillatory properties
\cite{CG2m} are exhibited by a semilinear diffusion equation with
``almost linear" logarithmic source term
 \be
 \label{1.41}
 u_t= -u_{xxxx} + u \, \ln^4 |u|.
  \ee
 The above models are semilinear and do not admit blow-up patterns
with finite interfaces.

 Concerning quasilinear higher-order PDEs,
 the {\em interface and blow-up phenomena} are natural and most well-known for the
degenerate  unstable
 {\em thin film equations} (TFEs)
  with
 lower-order terms
 such as
 \be
\label{GPP4}
    u_t = -(|u|^n u_{xxx})_x - (|u|^{p-1}u)_{xx,} \quad
 u_t = -\nabla \cdot (|u|^n \nabla \Delta u) \pm \D |u|^{p-1}u,  
  \ee
  where $n>0$ and $p>1$.
Equations of this form are  known to admit non-negative solutions
constructed by special sufficiently ``singular" parabolic
approximations of nonlinear coefficients that lead to
free-boundary problems. This direction was initiated by the
pioneering paper \cite{BF1} and was continued  by many
researchers; we refer to \cite{LPugh, WitBerBer} and the
references therein.
 Blow-up similarity solutions of the fourth-order TFE (\ref{GPP4})
 with the unstable sign ``$-$"
  \be
 \label{TF4.1}
 u_t= - (u^nu_{xxx})_x - (u^p)_{xx} \quad (u \ge 0),
 \ee
 have been also well studied and understood;
 see \cite{BerPugh1, BerPugh2, Bl4, SPugh, WitBerBer}, where
 further references on the mathematical properties of the models
 can be found. Countable sets of blow-up
 patterns for this TFE were described in \cite{Bl4}.

Interface and finite-time {\em extinction  behaviour}, which is
described by  various similarity patterns,  occur for other
reaction-absorption PDEs such as
 \be
 \label{GPP}
\mbox{$
 u_t=  
 -u_{xxxx} - |u|^{p-1}u
 $}
  \ee
in the {\em singular} parameter range
 \be
 \label{prange}
 \mbox{$
 p \in (- \frac 13,1),
  \,\, \Longrightarrow \,\,\, \mbox{$|u|^{p-1}u$ is not Lipschitz continuous
 at $u=0$},
  $}
  \ee
  so that $|u|^{p-1}u$ is not Lipschitz continuous
 at $u=0$;
see \cite{Galp1} and references therein.


\ssk

We have used a simply looking quasilinear model such as
(\ref{1.5}) to demonstrate various  new aspects of higher-order
 reaction-diffusion  blow-up
 phenomena. The  mathematics then
becomes  more difficult than for the second-order PDEs in
(\ref{s1}), where the Maximum Principle reveals its full capacity.
We do not expect straightforward rigorous justifications of
several our conclusions and results, and state key open problems
when necessary.


\section{Local asymptotic properties of solutions near interfaces}
  \label{SectLocR}

  Here, we  describe  generic oscillatory behaviour
  of solutions of (\ref{1.5}) close to finite interfaces.

\subsection{Local properties of  
travelling waves: oscillatory profiles for $\l < 0$}

We use  simple  TW
     solutions,
 \be
 \label{t1}
 u(x,t)=f(y), \quad y=x- \l t,
 \ee
to check generic propagation properties for reaction-diffusion
equations involved. In a wide class of 1D second-order
reaction-diffusion parabolic PDEs, the TWs rigorously describe the
behaviour of finite interfaces for general classes of solutions;
 see \cite[Ch.~7]{GalGeom} and references therein.

We use this approach for the fourth-order PDE (\ref{1.5}).
 The ODE
for $f$ takes the form
 \be
 \label{t2}
 -\l f'= - (|f''|^n f'')'' +|f|^{p-1}f.
  \ee
  By a local analysis near the singular point $\{f=0,f'=0\}$, it
  is not difficult to show that the higher-order term $|f|^{p-1}f$
  on the right-hand side is negligible.
Therefore, near interfaces, assuming that these are propagating,
we can consider the simpler equation
 \be
 \label{le2}
   (|f''|^n f'')'= -f \quad \mbox{for} \quad y>0, \quad f(0)=0,
  \ee
 which is obtained on integration once.
 Here we set $\l=-1$ for propagating waves, by scaling.
 We need to describe its
oscillatory solution of changing sign, with zeros concentrating at
the given interface point $y=0^+$. Oscillatory properties of
solutions are a common feature of related higher-order degenerate
ODEs; see pioneering paper by  Bernis--McLeod \cite{BMc91} for
similar fourth-order ODEs.

It follows from the scaling invariance of (\ref{le2}) that there
exist solutions of the form
 \be
 \label{le3}
  \mbox{$
 f(y) = y^\mu \varphi(s), \quad s= \ln y, \quad
 \mbox{where} \,\,\,\, \mu = \frac {2n+3}{n}> 2 \,\,\,\mbox{for} \,\,\, n>0,
  $}
  \ee
  where $\varphi(s)$ is called the
  {\em oscillatory component} of the given
   solution.
 Substituting (\ref{le3}) into (\ref{le2}) yields
 the following second-order equation for $\varphi(s)$:
 \be
 \label{le4}
 (n+1)|P_2(\vp)|^n P_3(\vp)= - \vp,
  \ee
   where $P_k$ denote linear differential operators
      (see
   \cite[p.~140]{GSVR}) given by the recursion
 $$
 \begin{matrix}
    P_{k+1}(\varphi)= P'_k(\varphi) + (\mu-k) P_k(\varphi), \quad
    k \ge 0; \quad P_0(\varphi)=\varphi, \quad \mbox{so
    that}\smallskip\smallskip\ssk
    \\
P_1(\varphi)=\vp'+ \mu \vp, \quad P_2(\vp)= \vp''+(2\mu-1) \vp'+
\mu(\mu-1)\vp, 
\qquad\quad\,\,
 \smallskip\smallskip \\
  P_3(\varphi)=\vp'''+ 3(\mu-1) \vp'' + (3 \mu^2- 6 \mu +2) \vp'
 + \mu(\mu-1)(\mu-2) \vp,
 \smallskip\smallskip \\
  P_4(\varphi)=\vp^{(4)}+ 2(2\mu-3) \vp''' + (6 \mu^2- 18 \mu +11)
  \vp'' \smallskip\smallskip \\
   +\, 2(2 \mu^3 -9 \mu^2 + 11 \mu -3) \vp'
+ \mu(\mu-1)(\mu-2)(\mu-3) \vp, \,\,\,\mbox{etc.}
 \end{matrix}
 $$

According to (\ref{le3}), we are interested in uniformly bounded
global solutions $\varphi(s)$ that are well defined as $s= \ln y
\to -\infty$, i.e., as $y \to 0^+$. The best candidates for such
global orbits of (\ref{le4}) are periodic solutions $\varphi_*(s)$
that are defined for all $s \in \re$. These  describe suitable
(and, possibly, generic) connections with the interface at
$s=-\infty$. The following result is proved by shooting as in
\cite[\S~7.1]{Gl4} and follows  the arguments in
\cite[\S~2]{GPME}.

\begin{proposition}
 \label{Pr.Per}

For all $n>0$,  $(\ref{le4})$
 has a periodic solution of changing sign $\varphi_*(s)$.

 \end{proposition}


There are two  open problems:

(i) uniqueness of the periodic solution $\vp_*(s)$, and

(ii) stability $\vp_*(s)$ as $s \to + \infty$.

\noi Numerical evidence answers positively to both questions. Then
(i) and (ii) mean  a unique (up to translation) periodic
connection with $s= -\infty$, where the interface is situated.



\begin{figure}
\centering \subfigure[$n=0.75$]{
\includegraphics[scale=0.52]{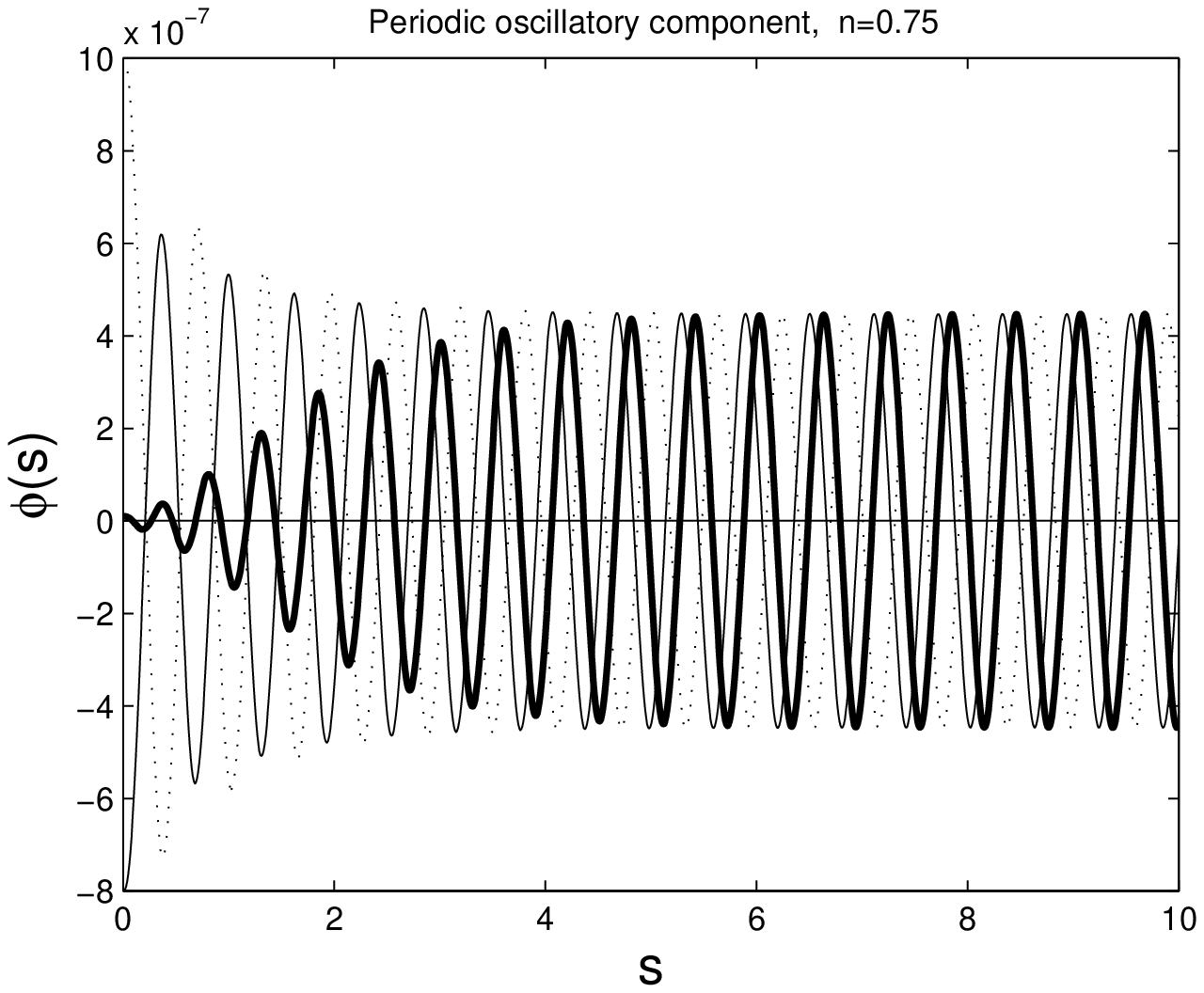}
}
\subfigure[$n=1$]{
\includegraphics[scale=0.52]{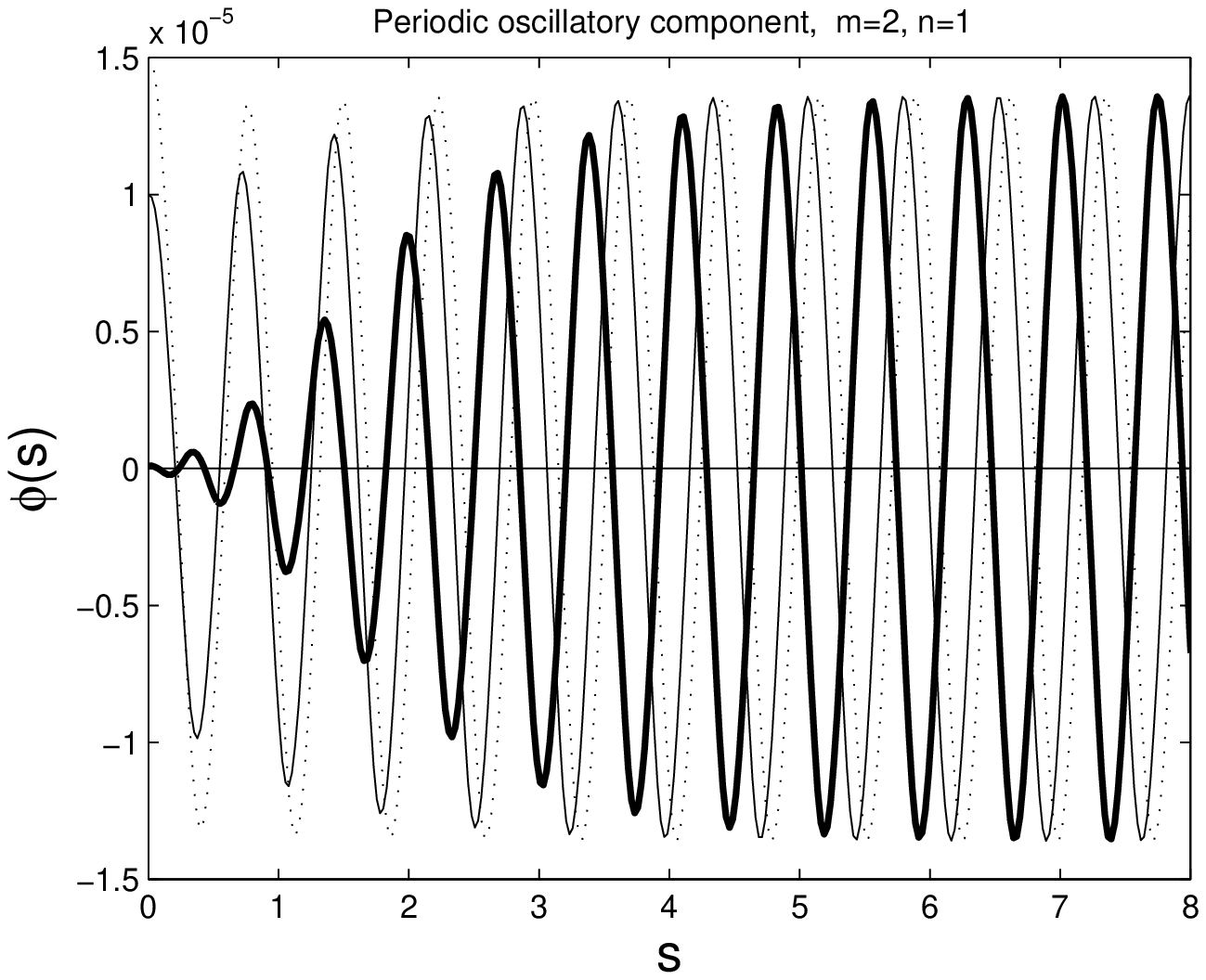}
}
 \subfigure[$n=3$]{
\includegraphics[scale=0.52]{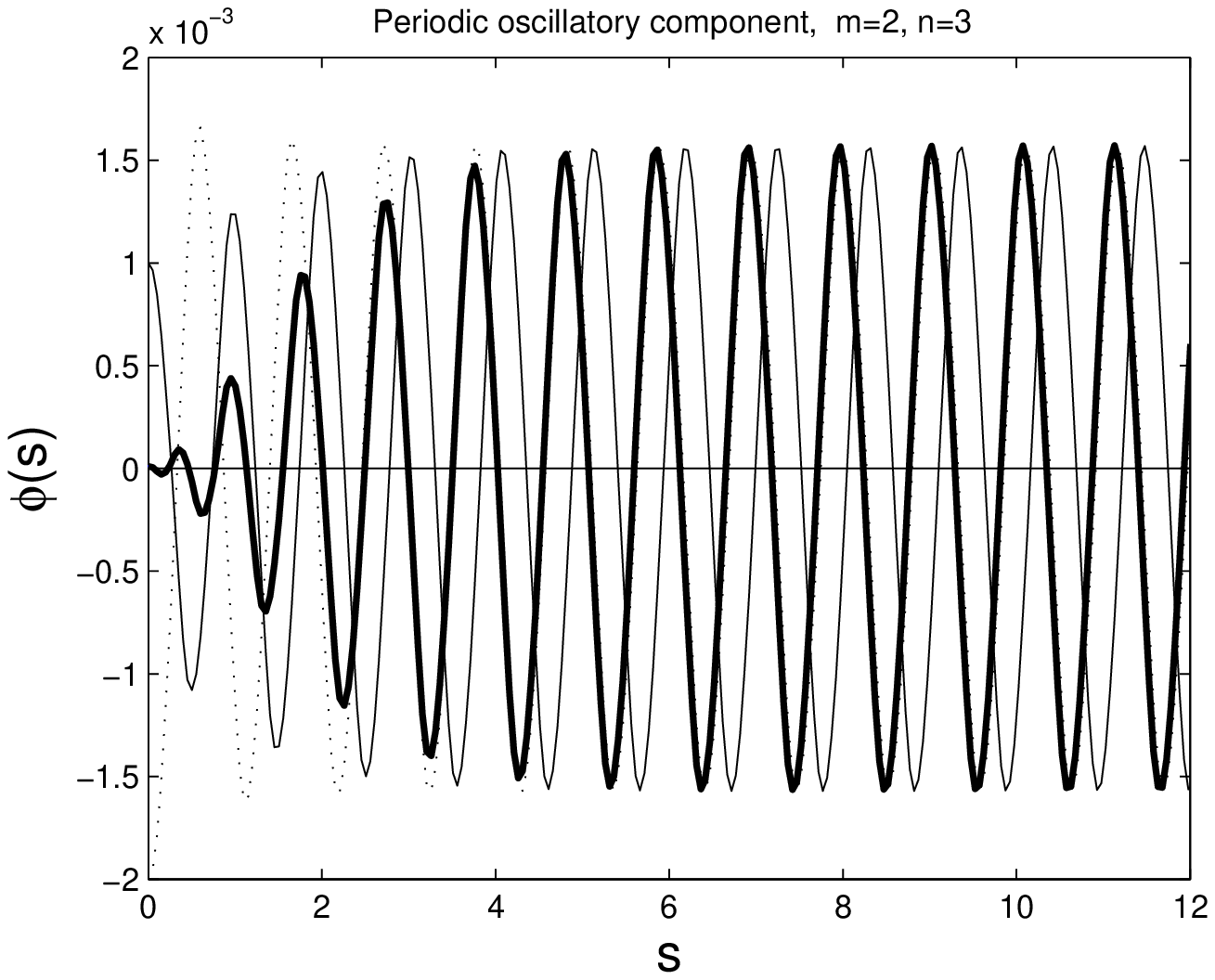}
}
 \subfigure[$n=5$]
{
\includegraphics[scale=0.52]{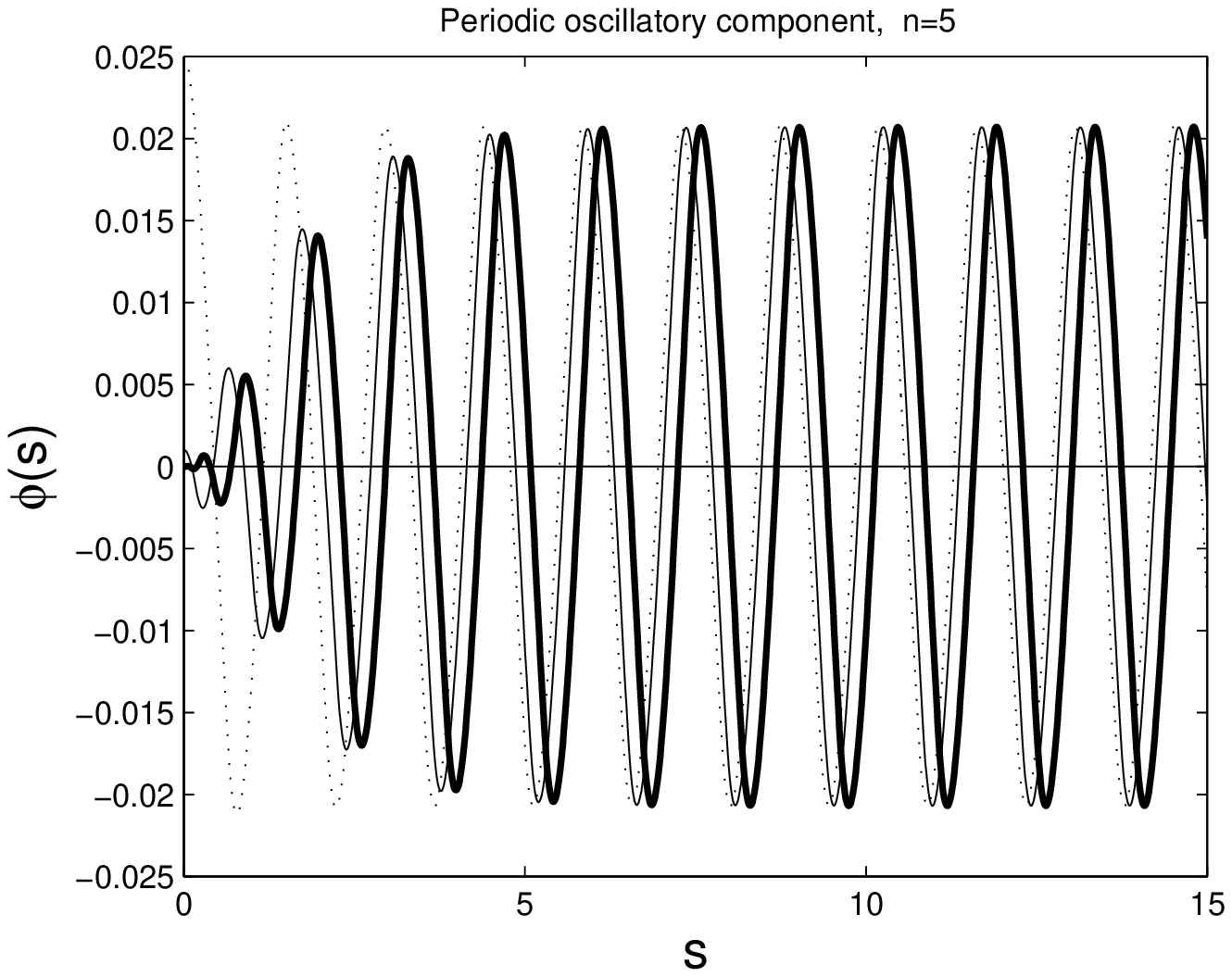}
}
 \vskip -.4cm
\caption{\rm\small Convergence to a stable periodic orbit of the
ODE (\ref{le4}) for  $n= \frac 34$, where $\varphi_* \sim
10^{-7}$, $n=1$, $n=3$, and $n=5$, with  $\varphi_* \sim
10^{-2}$.}
 \label{FOsc2}
\end{figure}

 The convergence to the unique stable periodic
 behaviour of (\ref{le4}) is shown in Figure \ref{FOsc2} for
 various $n = 0.75$ (periodic oscillations are of order $10^{-7}$) and $n=5$
 (order is $10^{-2}$).
 Different  curves therein correspond to different Cauchy data
$\varphi(0)$, $\varphi'(0)$,  $\varphi''(0)$  prescribed at $y=0$.
For $n < \frac 34$, the oscillatory component gets extremely
small, so an extra scaling is necessary as  explained in
\cite[\S~7.3]{Gl4}. A more accurate passage to the limit $n \to 0$
in the degenerate ODEs such as  (\ref{le4}) is presented there in
Section 7.6 and in Appendix B.


Finally, given the 
periodic $\varphi_*(s)$  of (\ref{le4}),
 as a natural
way to approach the interface point $y_0=0$ according to
(\ref{le3}), we have that the ODE (\ref{le2}) and, asymptotically,
(\ref{t2}),  admit  at the singularity set $\{f=0\}$
 \be
 \label{as55}
\mbox{a 2D local asymptotic family with parameters $y_0$ and phase
shift in $s \mapsto s + s_0$}.
 \ee
 We also call (\ref{as55}) an {\em asymptotic bundle} of orbits.


\subsection{Non-oscillatory case $\l >0$}

\begin{figure}
\centering \subfigure[stability of equilibria]{
\includegraphics[scale=0.52]{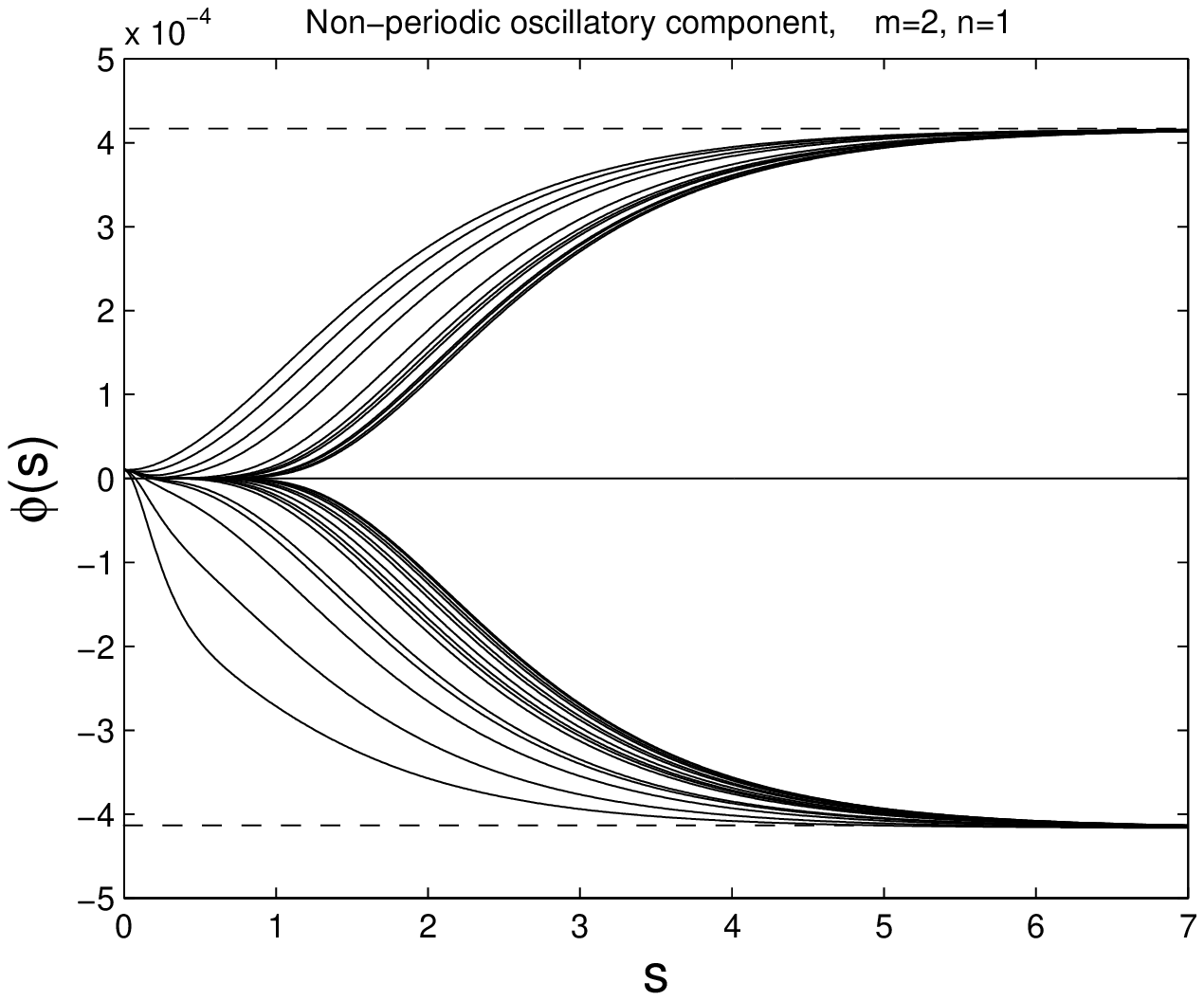}
} \subfigure[enlarged behaviour]{
\includegraphics[scale=0.52]{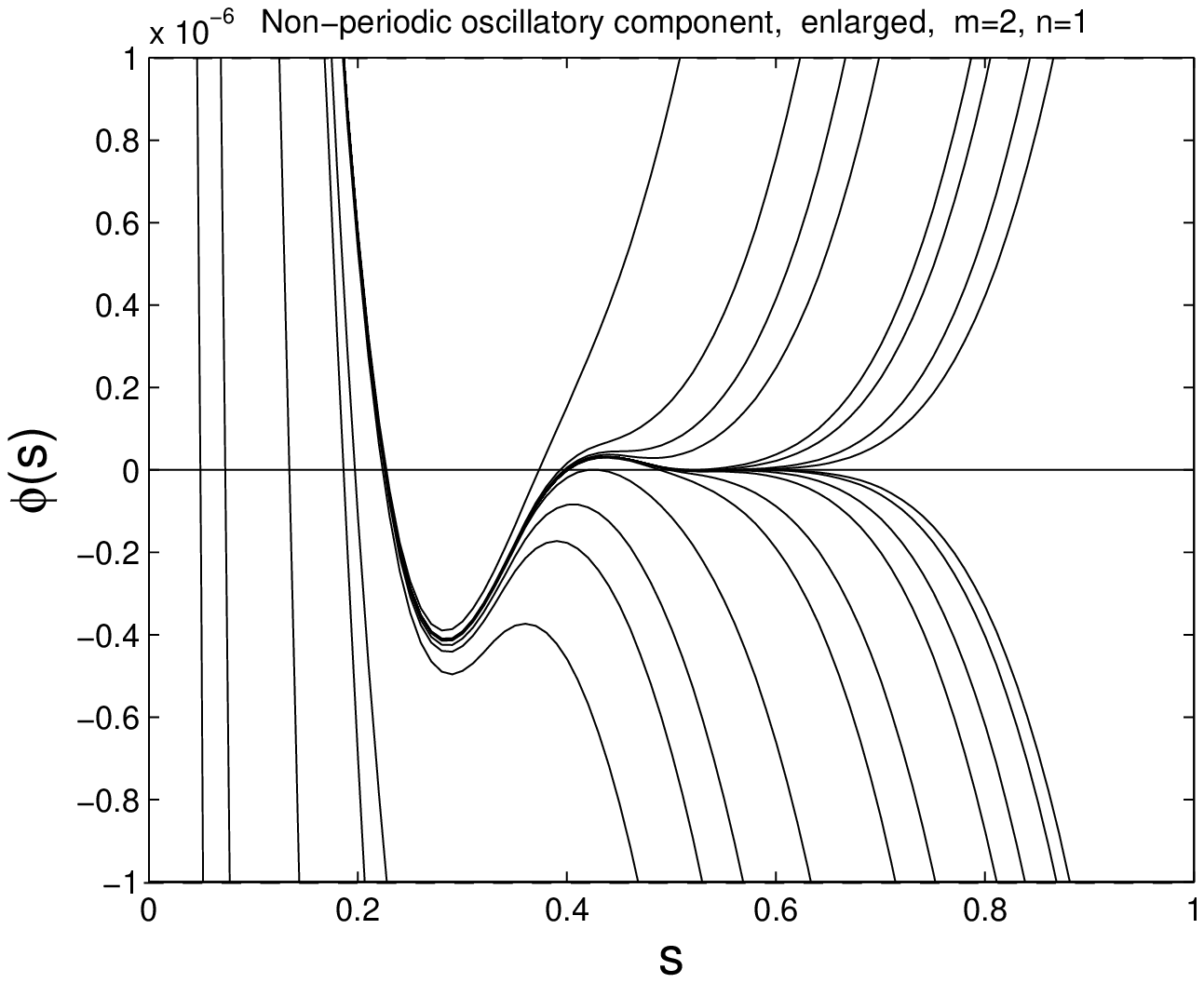}
}
 \vskip -.4cm
\caption{\rm\small Non-oscillatory behaviour for the ODE
(\ref{le4N}) for $n=1$; stability of equilibria (\ref{co1}) (a),
and enlarged non-periodic behaviour in between, (b).}
 \label{FOsc22}
\end{figure}

For $\l=1$, we have the opposite sign in the ODE
 \be
 \label{le4N}
 (n+1)|P_2(\vp)|^n P_3(\vp)=  \vp,
  \ee
  which admits two constant equilibria
 \be
 \label{co1}
 \mbox{$
 \varphi_\pm= \pm [(n+1)(\mu-2)]^{\frac 1 n} [\mu
 (\mu-1)]^{\frac{n+1}n}.
  $}
  \ee
 Figure
\ref{FOsc22}(a) shows that as $s \to + \infty$ the equilibria
(\ref{co1}) are stable (easy to see by linearization). In (b),
which gives the enlarged behaviour from (a) close to $\vp=0$,  we
observe a changing sign orbit, which is not  periodic. This
behaviour cannot be extended as a bounded solution up to the
interface at $s=-\infty$. In other similar  ODEs, which are
induced by other parabolic PDEs, such behaviour between two
equilibria can be periodic; cf. \cite[p.~143]{GSVR}.

These
 results confirm that for $\l>0$, the TWs {\em are not
oscillatory at interfaces},
 and actually such backward propagation via TWs is not possible for {\em almost
 all}
 (a.a.)
 initial data. More precisely, unlike (\ref{as55}), for $\l >0$,
 the asymptotic family (a bundle) as $s \to -\infty$ is 1D, which is not sufficient for
 matching purposes (see typical ideas of construction of similarity
 profiles below).



\section{Blow-up similarity solutions: problem setting and preliminaries}
 \label{Sect2}

\subsection{ODE reduction}
The parabolic PDE (\ref{1.5}) formally possesses the following
similarity solutions describing finite-time blow-up as $t \to
T^-$:
 \be
    \label{RVarsE}
 \mbox{$
 u_S(x,t) = (T-t)^{-\frac 1{p-1}}  f(y),
\quad y = {x}/(T-t)^{\b},  \quad \mbox{with} \quad \b =
\frac{p-(n+1)}{2(n+2)(p-1)}.
 $}
 \ee
The rescaled blow-up profile $f(y)$ satisfies the quasilinear
fourth-order ODE
 \be
  \label{f11E}
 \mbox{$
  {\bf A}(f) \equiv - (|f''|^n f'')''  - \b \, y f' -
\frac 1{p-1} \, f + |f|^{p-1}f=0 \quad \mbox{in} \quad \re.
 $}
 \ee
We impose  at the origin $y=0$ either the symmetry conditions,
   \be
\label{BCs}
 f'(0) = 0 \quad \mbox{and} \quad f'''(0) = 0,
\ee
 or the anti-symmetry ones,
    \be
\label{as9}
 f(0) = 0 \quad \mbox{and} \quad f''(0) = 0.
\ee
 By a standard local analysis of (\ref{f11E}) for small $f
\approx 0$, and
 in view of general results on regularity \cite[Ch.~1,2]{LIO}
  and  finite speed of propagation for such degenerate parabolic equations
  \cite{Shi2},
 a natural setting for the Cauchy problem assumes that, for $p \in
 (1,n+1]$,
   \be
 \label{com1}
 f(y) \,\, \mbox{is sufficiently smooth and compactly supported}.
  \ee
   The actual regularity of $f(y)$ close to
interfaces has been determined in the previous section.

For $p>n+1$, the asymptotic analysis shows that the solutions are
not compactly supported. Note that equation (\ref{f11E}) possesses
the constant equilibria
  \be
 \label{g1N}
 \pm f_*(p)=  \pm(p-1)^{-\frac 1{p-1}} . 
  \ee

\subsection{Blow-up self-similar profiles: preliminaries}
\label{SectEx}

We  next study solvability
   of the
ODE (\ref{f11E}) in $\re$. First of all,
 the local interface analysis from Section \ref{SectLocR}
  applies to 
  (\ref{f11E}). Indeed, close to the interface point
 $y=y_0>0$  of the similarity profile
$f(y)$,
 the ODE (\ref{f11E}) for $p <n+1$ contains the same leading terms as in
 (\ref{le2}) and other linear two are negligible as $y \to y_0^-$.

 For $p=n+1$, where $\b =0$, the leading terms close to the
 interface are
  $$
   \mbox{$
  -(|f''|^n f'')'' - \frac 1n \, f =0.
   $}
   $$
This gives solutions (\ref{le3}) with another exponent
 $$
  \mbox{$
 \mu= \frac{2(n+2)}n,
  $}
 $$
 and a fourth-order ODE for $\varphi(s)$, which admits a periodic
 solution $\varphi_*(s)$; see  examples in
 \cite[Ch.~3-5]{GSVR}.

It is important that, taking into account the local result
(\ref{as55})
 and bearing in mind the two conditions (\ref{BCs}) or (\ref{as9})
 yield two algebraic equations for two parameters $\{y_0,s_0\}$ of
 the bundle. Therefore,
we  expect that
 \be
 \label{co22}
 \mbox{there exists not more than a countable set $\{f_k\}$ of
 solutions}.
  \ee
 Note that this assumes a certain analyticity of the
 dependence on parameters in the degenerate ODE (\ref{f11E}), which
 is not easy to prove.
  In particular,  relative to the parameter
$p>1$, we can expect at most a countable set of $p$-branches of
solutions. This is true for the linear case $n=0$ and $p=1$;
 see below.

\subsection{Fundamental solution and necessary spectral properties}
 \label{S3.Sp}


Here we review some properties of differential operators in the
linear case $n=0$. Consider the linear {\em bi-harmonic equation}
 \be
 \label{bi1}
 u_t=- u_{xxxx} \quad \mbox{in} \quad \re \times \re_+.
  \ee
Its {\em fundamental
  solution} has the form
   \be
   \label{bi2}
    b(x,t)= t^{-\frac 14} F(y), \quad y = x/t^{\frac 14},
     \ee
  where
 the rescaled kernel $F$ is the unique radial solution of the ODE 
  \be
\label{ODEf}
 {\bf B} F \equiv - F^{(4)} +  \mbox{$\frac 1{4}$}\, y
 F' + \mbox{$\frac 1{4}$}\, F =0 \quad
 \mbox{in} \,\,\, \re,  \quad \mbox{with} \,\,\, \mbox{$\int\limits_\re$} \, F\, {\mathrm d}y = 1.
 \ee
 On integration once, we obtain a third-order equation,
 \be
 \label{F3New}
  \mbox{$
 - F''' + \frac 14 \, y F=0 \quad \mbox{in} \,\,\, \re.
  $}
  \ee

The kernel $F=F(|y|)$ is radial, has exponential decay, oscillates
as $|y| \to \infty$,
 and
\be \label{fbar}
 |F(y)| \le  D \, {\mathrm e}^{-d|y|^{4/3}}
\quad \mbox{in} \,\,\, \re,
 \ee
  for a
positive constant $D$ and
  $d = 3 \cdot 2^{-11/3}$; see  \cite[p.~46]{EidSys}.
The necessary spectral properties of the linear non self-adjoint
operator ${\bf B}$ and the corresponding adjoint operator ${\bf
B}^*$ are of importance in the  asymptotic analysis
 and are explained in \cite{Eg4} for general $2m$th-order
 operators (see also \cite[\S~4]{Bl4}).
 In particular, ${\bf B}$ has a discrete (point) spectrum $\s({\bf B})$
in a weighted space $L^2_\rho(\re)$, with $\rho(y)={\mathrm e}^{a
|y|^{4/3}}$, $a \in (0,2d)$ is a constant,
 \begin{equation}
\label{spec1}
  \sigma({\mathbf B}) = \{\lambda_l = -\mbox{$\frac
l4$}, \,\, l = 0,1,2,...\}.
\end{equation}
The corresponding eigenfunctions are given by
 \be
 \label{psi1}
 \mbox{$
 \psi_l(y)= \frac {(-1)^l}{\sqrt{l!}} F^{(l)}(y), \quad l=0,1,2,... \,
 .
  $}
  \ee

The adjoint operator
 \be
 \label{ad1}
  \mbox{$
 {\bf B}^*=- D_y^4 - \frac 14 \, y D_y
  $}
  \ee
  has the same spectrum (\ref{spec1}) and polynomial
  eigenfunctions
  \begin{equation}
 \label{psi**1}
  \mbox{$
 \psi_l^*(y) = \frac 1{\sqrt{l !}}
\sum\limits_{j=0}^{\lfloor-\lambda_l\rfloor} \frac {1}{j !}D^{4j}_y y^l,
\quad l=0,1,2,... \, ,
 $}
 \end{equation}
which form a complete subset  in $L^2_{\rho^*}(\re)$, where
$\rho^*=\frac 1 \rho$. As ${\bf B}$, the adjoint operator ${\bf
B}^*$ has compact resolvent $({\bf B}^*-\lambda I)^{-1}$. It is
not difficult to see  by integration by parts that  the
eigenfunctions
 (\ref{psi1})
 are orthonormal to polynomial eigenfunctions
$\{\psi^*_l\}$ of the adjoint operator ${\bf B}^*$, so
 \be
 \label{ort1}
 \langle \psi_l, \psi^*_k \rangle = \d_{lk},
  \ee
  where $\langle \cdot, \cdot \rangle$ denotes the standard (dual)
  scalar product in $L^2(\re)$.

\subsection{Countable set of similarity solutions for
$n=0$, $p=1$}

Performing in the equation (\ref{bh1})
 the change
  \be
  \label{w11}
  u(x,t)= {\mathrm e}^t \, w(x,t)
   \ee
    reduces it to
the pure bi-harmonic equation  (\ref{bi1})
 for $w(x,t)$. By the scaling as for the fundamental
solution $b(x,t)$ in (\ref{bi2}),
 \be
 \label{v11}
 w(x,t)= t^{-\frac 14} v(y,\t), \quad y=x/t^{\frac 14}, \,\, \t =
 \ln t,
  \ee
  we obtain the rescaled equation with the ${\bf B}$ in
  (\ref{ODEf}) having eigenfunctions (\ref{psi1}), so
   \be
   \label{v12}
   v_\t = {\bf B} v \quad \Longrightarrow \quad \exists \,\,\, v_l(y,\t)=
   {\mathrm e}^{\l_l \t} \psi_l(y).
    \ee
Setting $\l_l=- \frac l4$ as in (\ref{spec1}) and $t={\mathrm
e}^\t$, we obtain a countable set of different asymptotic patterns
for the linear PDE (\ref{bh1}) corresponding to $n=0$ and  $p=1$:
 \be
 \label{u11}
 \mbox{$
 u_l(x,t) = {\mathrm e}^{-t} \,\, t^{-\frac {1+l}
4} \psi_l \big(\frac
 x{t^{1/4}}\big), \quad l=0,1,2,... \, .
  $}
  \ee

 It turns out that the blow-up similarity patterns (\ref{RVarsE}) can be deformed
as $n \to 0$ and  $p \to 1$ to those in (\ref{u11}) (though
entirely rigorous proof is very difficult and not fully completed
for such degenerate equations, as will happen for some other
related homotopy questions). Then (\ref{u11})
suggests that 
there exists a countable number of branches $\{f_l(y;n,p)\}$,
which
 appear from the {\em branching point}
  $\{n=0, \,  p=1\}$ according to classic theory, \cite[\S~56]{KrasZ}.  
We claim that the above two (linear for $n=0$, $p=1$ and nonlinear
for $n>0$, $p>1$) asymptotic problems admit a continuous {\em
homotopic} connection as $n \to 0$, $p \to 1$, so that, after
necessary scaling, (\ref{u11}) is obtained in the limit from
nonlinear eigenfunctions. For such ODEs, this reduces to a matched
asymptotic expansion analysis, which is rather technical and is
not studied here.

What is key for the future study
 is that the oscillatory behaviour
of linear patterns in (\ref{u11}) is then inherited by nonlinear
blow-up patterns at least for small $n>0$ and  $p>1$. This shows
once more that similarity profiles $f(y)$ corresponding to the
Cauchy problem must be oscillatory near interfaces.
Homotopy  approaches can play a role for specifying correct
settings of the Cauchy problem  for variety of nonlinear PDEs with
non-smooth or singular coefficients, if they share the same {\em
homotopy class} with a well-posed
 linear equation; see \cite[Ch.~8]{Gl4}.

 It follows from the ODE (\ref{f11E}) that
 \be
 \label{g1}
 \|f\|_\infty \sim f_*(p)= (p-1)^{-\frac 1{p-1}} \to +\infty \quad \mbox{as \,
 $p \to 1^+$},
  \ee
  so the divergence (in fact, towards the rescaled linear problem)  is exponentially fast.

\section{Regional blow-up  profiles for $p=n+1$: variational approach}
 \label{SectS}

We begin with the special case $p=n+1$, where $\b=0$ in
(\ref{RVarsE}) (so $y=x$) and $f(y)$ in (\ref{f11E}) solves an
autonomous fourth-order ODE of the form
 \be
  \label{f11ES}
 \mbox{$
  {\bf A}(f) \equiv - (|f''|^n f'')''   -
\frac 1{n} \, f + |f|^{n}f=0 \quad \mbox{in} \quad \re.
 $}
 \ee
This is a variational problem that can be studied in greater
detail. Later on, we apply these patterns and classification for
$p=n+1$ in neighbouring parameter ranges $p>n+1$ and $p<n+1$ by
using a natural idea of $p$-branches of solutions.

For convenience, we perform in (\ref{f11ES}) an extra scaling
 \be
 \label{1}
  \mbox{$
 f= \big(\frac 1n\big)^{\frac 1n} F \quad \Longrightarrow
 \quad - (|F''|^n F'')''   -
 F + |F|^{n}F=0 \quad \mbox{in} \quad \re.
 $}
  \ee
  For any $n>0$, this equation admits three constant equilibria
 $$
 F \equiv -1, \,\, 0,\,\, 1.
  $$

\subsection{Variational setting and compactly supported solutions}

 Operators involved in the ODE (\ref{1}) are potential in $L^2$, so the problem
admits a variational setting and solutions can be obtained as
critical points of a $C^1$ functional   of the form
 \be
 \label{V1}
  \mbox{$
 {E}(F)= - \frac 1{n+2} \,  \int\limits |F''|^{n+2}\, {\mathrm d}y - \frac 12 \int
 F^2 \, {\mathrm d}y +\frac{1}{n+2} \, \int\limits |F|^{n+2} \, {\mathrm d}y.
  $}
  \ee
  Then we are looking for critical points in $W^{n+2}_2(\re)
  \cap L^2(\re) \cap L^{n+2}(\re)$.
  For 
 compactly supported solutions (see below), we choose a sufficiently
  large interval $B_R=(-R,R)$ and consider the variational problem for (\ref{V1})
  in $W_{2,0}^{n+2}(B_R)$,  assuming Dirichlet boundary
  conditions at the end points $\partial B_R=\{\pm R\}$.
  By Sobolev embedding theorem, $W_{2,0}^{n+2}(B_R)$ is  compactly embedded into
    $L^2(B_R)$ and $L^{n+2}(B_R)$.
   Continuity of any bounded solution $F(y)$ is guaranteed
 by Sobolev embedding $H^2(\re) \subset C(\re)$.

Thus, we will be looking for compactly supported solutions. This
demand is associated with the well-known fact that the
corresponding parabolic flow with the elliptic operator as in
(\ref{1}),
 \be
 \label{par1}
 w_t= - (|w_{xx}|^n w_{xx})_{xx} - w + |w|^n w,
  \ee
describes processes with finite propagation of interfaces.
 By energy estimates, such
  results have been proved for a number of
  quasilinear higher-order parabolic equations with
 potential $p$-Laplace-type operators;
 see \cite{Shi2}. Therefore, our blow-up patterns are indeed
 nontrivial compactly supported stationary solutions of (\ref{par1}).
 Examples of  ODE proofs via typical energy estimates can be
 found in \cite[\S~7]{BMc91}.

\smallskip


 Thus, in what follows, to revealing  compactly supported
patterns $F(y)$, we will pose the problem in bounded sufficiently
large intervals $(-R,R)$  with Dirichlet data at $\pm R$.

 \subsection{L--S theory and direct
 application of fibering method}

 The functional (\ref{V1}) is $C^1$,
uniformly differentiable, and weakly continuous, so we can apply
classic
 Lusternik--Schnirel'man (L--S) theory of calculus of variations
\cite[\S~57]{KrasZ} in the form of the fibering method
  \cite{Poh0, PohFM}.




 According to L--S theory
   and the fibering
 approach,
  the number of critical points of the
functional (\ref{V1}) depends on the {\em category} (or {\em
genus}) of functional subset on which the fibering is taking
place. The critical points of ${E(F)}$ are convenient to obtain by
the {\em spherical fibering} in the form
 \be
 \label{f1}
 F= r(v) v \quad (r \ge 0).
  \ee
  Here $r(v)$ is a scalar functional, and $v$ belongs to a subset
  in  $W_{2,0}^{n+2}(B_R)$ given by
   \be
   \label{f2}
    \mbox{$
    {\mathcal H}_0=\bigl\{v \in W_{2,0}^{n+2}(B_R): \,\,\,H_0(v)
     \equiv  -  \int |v''|^{n+2}\, {\mathrm d}y +  \int
 |v|^{n+2}\, {\mathrm d}y =1\bigr\}.
    $}
    \ee
Then the new functional
 \be
 \label{f3}
 \mbox{$
H(r,v)= E(r v) \equiv \frac 1{n+2} \, r^{n+2} - \frac 1{2}\, r^{2}
\int\limits v^{2}\, {\mathrm d}y
 $}
  \ee
 has the absolute minimum point, where
 \be
 \label{f31}
  \begin{matrix}
  \mbox{$
 H'_r \equiv r^{n+1}-  r \int v^{2}\, {\mathrm d}y =0
  \,\,\Longrightarrow \,\,
   r_0(v)=\big(\int v^{2}\, {\mathrm d}y\big)^{\frac 1{n}},
   $}
   \ssk \ssk\ssk\ssk \\
    \mbox{at which} \quad
    \mbox{$
   H(r_0(v),v)=- \frac {n}{2(n+2)} \,  r_0^{n+2}(v).
    $}
    \end{matrix}
    \ee
Therefore, introducing
 \be
 \label{f4}
 \mbox{$
 \tilde H(v) = \bigl[ - \frac {2(n+2)}{n}H(r_0(v),v)
 \bigr]^{\frac {n}{n+2}} \equiv \int v^{2}\, {\mathrm d}y,
  $}
  \ee
  we arrive at the  quadratic, even, non-negative, convex,
 and uniformly differentiable functional, to which
    L--S  theory applies, \cite[\S~57]{KrasZ}.
   Searching
  for critical points of $\tilde H$ in the set ${\mathcal H}_0$,
  one needs to estimate the category-genus $\rho$
  of the set ${\mathcal H}_0$.
The details on this notation and basic results for semilinear
equations can be found in Berger \cite[p.~378]{Berger}. The Morse
index $q$ of the quadratic form $Q$ in Theorem 6.7.9 therein is
precisely the dimension of the space where the corresponding form
is negatively definite. This includes all the multiplicities of
eigenfunctions involved in the corresponding subspace. Note that
Berger's analysis and most of others  are dealing  with
perturbation theory of  linear operators,
 which makes it easier to get the genus of necessary functional
 sets involved.
 For the quasilinear operators
that define the set (\ref{f2}) by their potentials, an extra study
of genus is needed
 (to be performed below).




For detecting  geometric shapes of patterns, we recall that by the
minimax analysis of L--S category theory \cite[p.~387]{KrasZ},
  \cite[p.~368]{Berger}, the
critical values $\{c_k\}$ and the corresponding critical points
$\{v_k\}$ are given by
 \be
 \label{ck1}
  \mbox{$
 c_k = \inf_{{\mathcal F} \in {\mathcal M}_k} \,\, \sup_{v \in {\mathcal
 F}} \,\, \tilde H(v),
  $}
  \ee
where  ${\mathcal F} \subset {\mathcal H}_0$ are  closed sets,
 and
 ${\mathcal M}_k$ denotes the set of all subsets of the form
  $$
  B S^{k-1}
\subset {\mathcal H}_0,
 $$
 where $S^{k-1}$ is a suitable sufficiently
smooth $(k-1)$-dimensional manifold (say, sphere) in ${\mathcal
H}_0$ and $B$ is an odd continuous map.
 Then each member of ${\mathcal M}_k$ is of  genus at least $k$
 (available in ${\mathcal H}_0$).
   It is also important to remind that the
definition of genus \cite[p.~385]{KrasZ} assumes  that
$\rho({\mathcal F})=1$, if no {\em component} of ${\mathcal F}
\cup {\mathcal F}^*$, where
 $$
 {\mathcal F}^*=\{v: \,\, -v \in {\mathcal F}\},
 $$
 is the {\em reflection} of ${\mathcal F}$ relative to 0,
 contains a pair of antipodal points $v$ and $v^*=-v$.
 Furthermore, $\rho({\mathcal F})=n$ if each compact subset of
${\mathcal F}$ can be covered by, minimum, $n$ sets of genus one.

According to (\ref{ck1}),
 $$
 c_1 \le c_2 \le ... \le c_{l_0},
 $$
 where $l_0=l_0(R)$ is the category of ${\mathcal H}_0$ satisfying (see
 below)
  \be
  \label{l01}
  l_0(R) \to + \infty \quad \mbox{as} \quad R \to \infty.
  \ee
  Roughly speaking,
since the dimension of the sets ${\mathcal F}$ involved in the
construction of ${\mathcal M}_k$ increases with $k$, this
guarantees that the critical points delivering critical values
(\ref{ck1}) are all different.


\subsection{Category of ${\mathcal H}_0$ gets arbitrarily large as $R \to +
\infty$}

  It follows from \cite[p.~385]{KrasZ}, \cite[p.~376]{Berger} (see also \cite{PohFM})
   that according to (\ref{f2}),  the category
$l_0=\rho({\mathcal H}_0)$  of the set ${\mathcal H}_0$ can be
associated with  the maximal number $K=K(R)$  of nonlinear
eigenvalues $\l_k <1$
 of the corresponding elliptic problem
 \be
 \label{f55}
 - (|\psi''|^n \psi'')'' + \l_k |\psi|^n \psi=0, \quad \psi \in W^2_{2,0}(B_R).
  \ee
This problem is solved by L--S theory and gives at least a
countable set of critical values and different critical points of
the positive homogeneous functional
 \be
 \label{ls1}
 \mbox{$
 \int |v|^{n+2}\, {\mathrm d}y \quad \mbox{on the unit sphere}
 \quad S_1=\big\{ \int |v''|^{n+2}\, {\mathrm d}y=1 \big\}.
 $}
  \ee
Indeed, given an eigenfunction $\psi_k \not = 0$ with $\l_k<1$,
multiplying (\ref{f55}) by $\psi_k$ yields
 $$
   \begin{matrix}
 - \int |\psi_k''|^{n+2}\, {\mathrm d}y + \int |\psi_k|^{n+2}\, {\mathrm d}y = (1- \l_k)
\int |\psi_k|^{n+2}\, {\mathrm d}y>0 \quad \Longrightarrow
\smallskip\smallskip\smallskip
\\
 \bar \psi_k=B_k \psi_k \in {\mathcal H}_0,
 \quad  |B_k|^{n+2}= \big[(1-\l_k) \int |\psi_k|^{n+2} \, {\mathrm
 d}y\big]^{-1},
  \end{matrix}
 $$
 where $B_k>0$
is the necessary normalization factor.
 By L--S theory, all such nonlinear eigenfunctions are different
 (since correspond to different critical values of the
 functional), so that all of them $\{\bar\psi_k, \,\, k=1,...,K\}$ are
 linearly independent.
 In order to estimate the genus of ${\mathcal H}_0$, we take their
 linear combination
  \be
  \label{ll1}
  v= C_1 \bar \psi_1 + ...+C_K \bar \psi_K \in {\mathcal H}_0,
   \ee
   so on substitution into the functional in (\ref{f2}) we get the
   following algebraic equation for the coefficients
   ${\bf C}=\{C_1,...,C_K\} \in \re^K$:
   \be
   \label{al11}
    \mbox{$
    G({\bf C}) \equiv
- \int |C_1 \bar \psi_1''+...+C_K \bar \psi_K''|^{n+2}\,{\mathrm
d}y + \int |C_1 \bar \psi_1+...+C_K \bar \psi_K|^{n+2}\, {\mathrm
d}y =1,
 $}
 \ee
 which is an equation of a surface ${\mathcal L}_K$ in $\re^K$
 being  symmetric under the reflection
  \be
  \label{ss33}
  {\bf C} \mapsto -{\bf C}.
   \ee
One can see that, by construction of the normalized eigenfunctions
$\bar \psi_k$, for any fixed  $k=1,2,...,K,$
 \be
 \label{g11}
G({\bf C})= |C_k|^{n+2}(1+o(1)) \quad \mbox{as} \quad C_k \to \infty.
 \ee
It is not difficult to see (using the variational and extremal
nature of nonlinear eigenfunctions) that ${\mathcal L}_K$ contains
a simple closed connected component, which, in view of
(\ref{ss33}), is homotopic to the unit sphere $S^{K-1}$ in
$\re^K$. By the ``additivity"   properties of the genus, this
implies that
  \be
  \label{hh0}
   \rho({\mathcal H}_0) \ge K(R)-1.
    \ee
    We do not know whether
 this estimate is sharp:
 optimal estimates of the category (genus) of the sets and even multiplicity of nonlinear
  eigenfunctions for such functionals compose a difficult open
  problem, which persists even for classic $p$-Laplacian operators
  as in (\ref{p1}).


  Since the dependence of the spectrum on the length $R$ for (\ref{f55}) is,
  by simple scaling,
   \be
   \label{f56}
   \l_k(R)= R^{-4-2n} \l_k(1) \to 0^+
    \,\,\, \mbox{as} \,\,\, R \to \iy, \quad k=0,1,2,... \, ,
    \ee
we have that the category $\rho({\mathcal H}_0)$ can be
arbitrarily large for $R \gg 1$, and (\ref{l01}) holds:

\begin{proposition}
 \label{Pr.MM}
The ODE problem $(\ref{1})$ has at least a countable set of
different solutions denoted by $\{F_l, \, l \ge 0\}$, and each one
 $F_l(y)$ is
obtained as a critical point of the functional $(\ref{V1})$ in
$W^2_{m,0}(B_R)$ with sufficiently large $R=R(l) > 0$.
 \end{proposition}



 \subsection{First basic pattern and local structure of zeros}

\begin{figure}
 \centering
\includegraphics[scale=0.8]{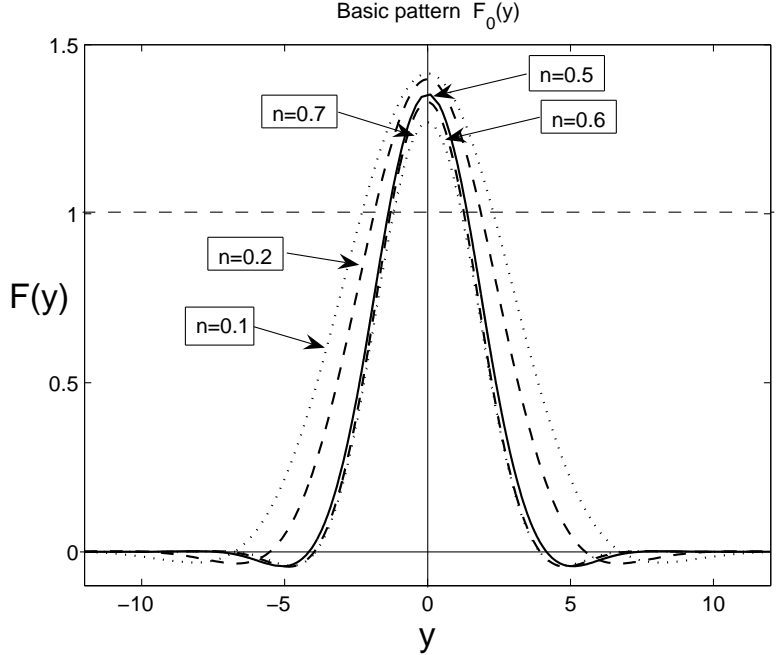}  
 \vskip -.4cm
\caption{\rm\small The first solution $F_0(y)$  of  (\ref{1}) for
various $n$.}
   \vskip -.3cm
 \label{G1}
\end{figure}

 Let us  present numerical results concerning existence
 and multiplicity of solutions for equation (\ref{1}).
 In Figure
\ref{G1}, we show the first basic 
pattern for (\ref{1}) called the  $F_0(y)$ for various $n \in
[0.1,0.7]$.
 These profiles are  constructed by {\tt MatLab}
 by using a natural regularization in the singular term,
 \be
 \label{4.1}
  \mbox{$
 - \big[(\e^2+(F'')^2)^ {\frac n 2} F''\big]''   -
 F + |F|^{n}F=0 \quad \mbox{in} \quad \re \quad (\e>0).
 $}
 \ee
 Here,  the regularization parameter $\e$
 and both tolerances in the {\tt bvp4c}
solver, typically,   take the values
 \be
 \label{eps1}
  \e= 10^{-2} \,\,\, \mbox{or} \,\,\, 10^{-3} \quad  \mbox{and} \quad  {\rm Tols}=10^{-3}
  \,\,\, \mbox{or} \,\,\, 10^{-4}.
  \ee
 For $n>0.5$, convergence gets rather slow.
 For $n \le 0.7$, the global structure of blow-up profiles (excluding their fine zero
 structure, see below) is  stable with respect to reasonable variations
 of $\e$ and Tols. In fact, this  reflects the structural
 stability of first basic blow-up patterns, which the author observed in
 dozens of other nonlinear parabolic  models with blow-up. Note that
 proving
  stability  even in the linearized setting involves
 non self-adjoint operators  with non-constant coefficients that
 leads to several technical difficulties and remains open.
 On the other hand,
   for $n \ge 1$, i.e., for strongly nonlinear diffusion operators in (\ref{4.1}), we
 did not get  reliable enough numerical results with the
 necessary accuracy, so we will avoid using such cases for further
 illustrations.

\begin{figure}
\centering \subfigure[scale $10^{-2}$]{
\includegraphics[scale=0.52]{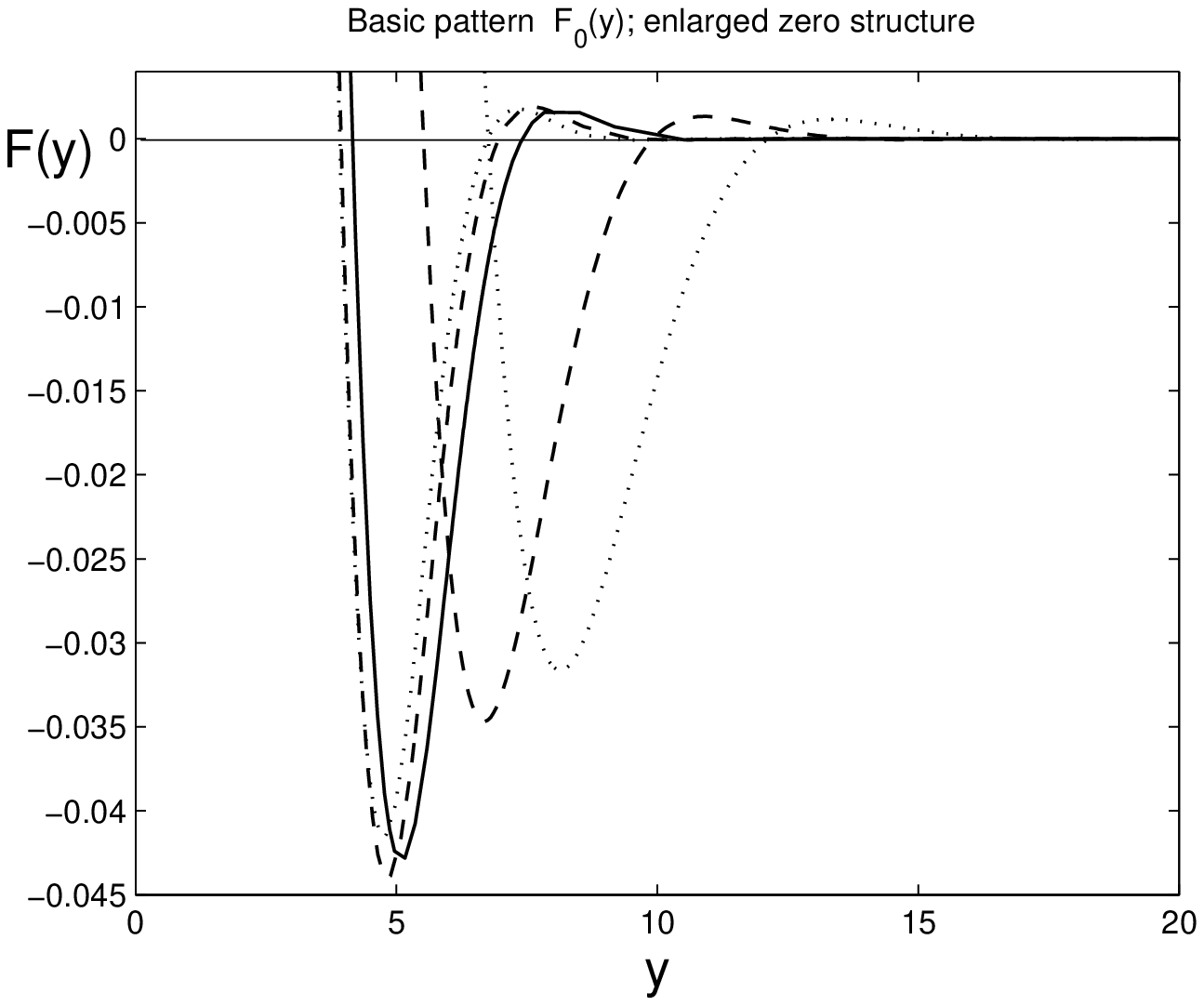}
} \subfigure[scale $10^{-3}$]{
\includegraphics[scale=0.52]{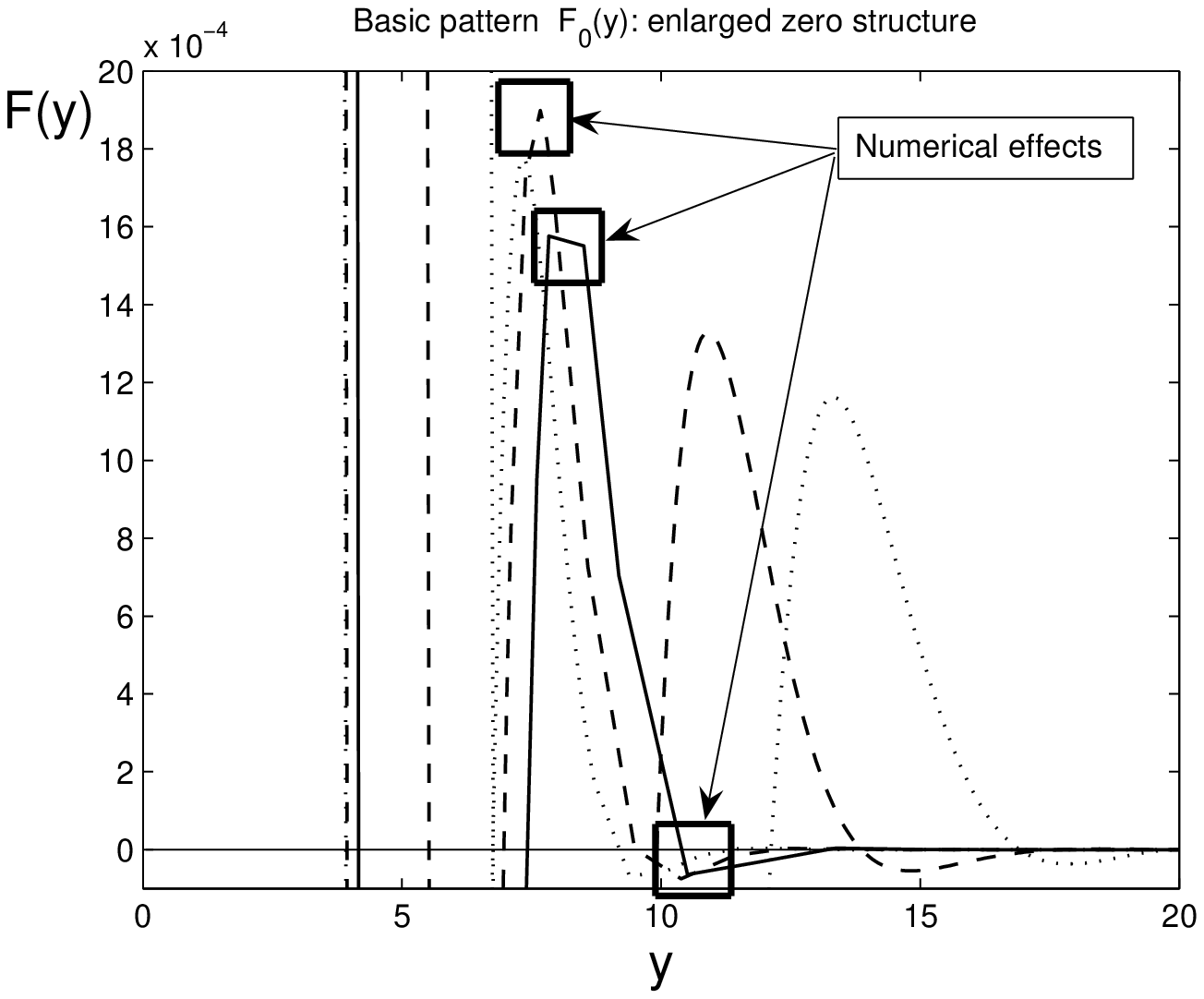}               
}
 \vskip -.2cm
\caption{\rm\small  Enlarged zero structure of the profile
$F_0(y)$  in  the linear scale.}
 \label{G2}
\end{figure}

 Incidentally, this
  makes it possible to reveal some features of the  local structure of
 multiple zeros close to the interface. Figure \ref{G2}
   shows how
the zero structure of profiles $F_0$ from Figure \ref{G1}  repeats
itself in a ``self-similar manner" from one zero to another in the
usual linear scale. In Figure (b),  a ``discrete", piece-wise
continuous structure for $n=0.5$ is already revealed, and this is
the best we have been able to achieve numerically. However, this
makes no problem, since the accuracy $10^{-3}$ achieved in (b) is
already in agreement with parameters in (\ref{eps1}), so further
improvements make no practical sense. In addition, this shows that
the discrete and continuous solutions of this difficult
variational problem remain very similar even for the present rough
meshes, when the discrete features become clearly visually
observable (as usual, it is a key fact for such numerics).


Further revealing zero structure and eventually the behaviour such
as (\ref{le3}) as $s=y_0-y \to 0^+$
cannot be reliable done in the parameter range (\ref{eps1}). In
\cite{Bl4, Gl4}, for similar thin film models,
 this demanded $\e$
and Tols to achieve at least $10^{-12}$, which is not possible for
the current model in view of slow convergence for higher-order
$p$-Laplacians.
 It is also quite a challenge to detect numerically the free-boundary
point.
 The main difficulty is to distinguish the nonlinear oscillations
 via (\ref{le3}) and the linear ones in
 the ``linearized area", where (\ref{4.1}) implies an exponential
 behaviour for $y \gg 1$ governed by the ODE
   \be
   \label{lin1}
   \mbox{$
 F^{(4)}= - 
 \e^{-n} F+...
   \quad \Longrightarrow \quad
   F(y) \sim {\mathrm e}^{- \frac {\sqrt 2}2 \, \e^{-  n/4}} \,
    \cos\big( \frac {\sqrt 2}2 \,  \e^{-\frac n4} y +c \big),
     $}
   \ee
  where $c$ is a constant.
 Actually, we saw not more than first 1--3  nonlinear zeros of the
 type (\ref{le3})
    and the
  rest of zeros corresponded to the linear behaviour (\ref{lin1}).


  \subsection{Basic countable family: approximate Sturm's
property}
 \label{Sect54}

\begin{figure}
\centering \subfigure[ $F_0(y)$ ]{
\includegraphics[scale=0.52]{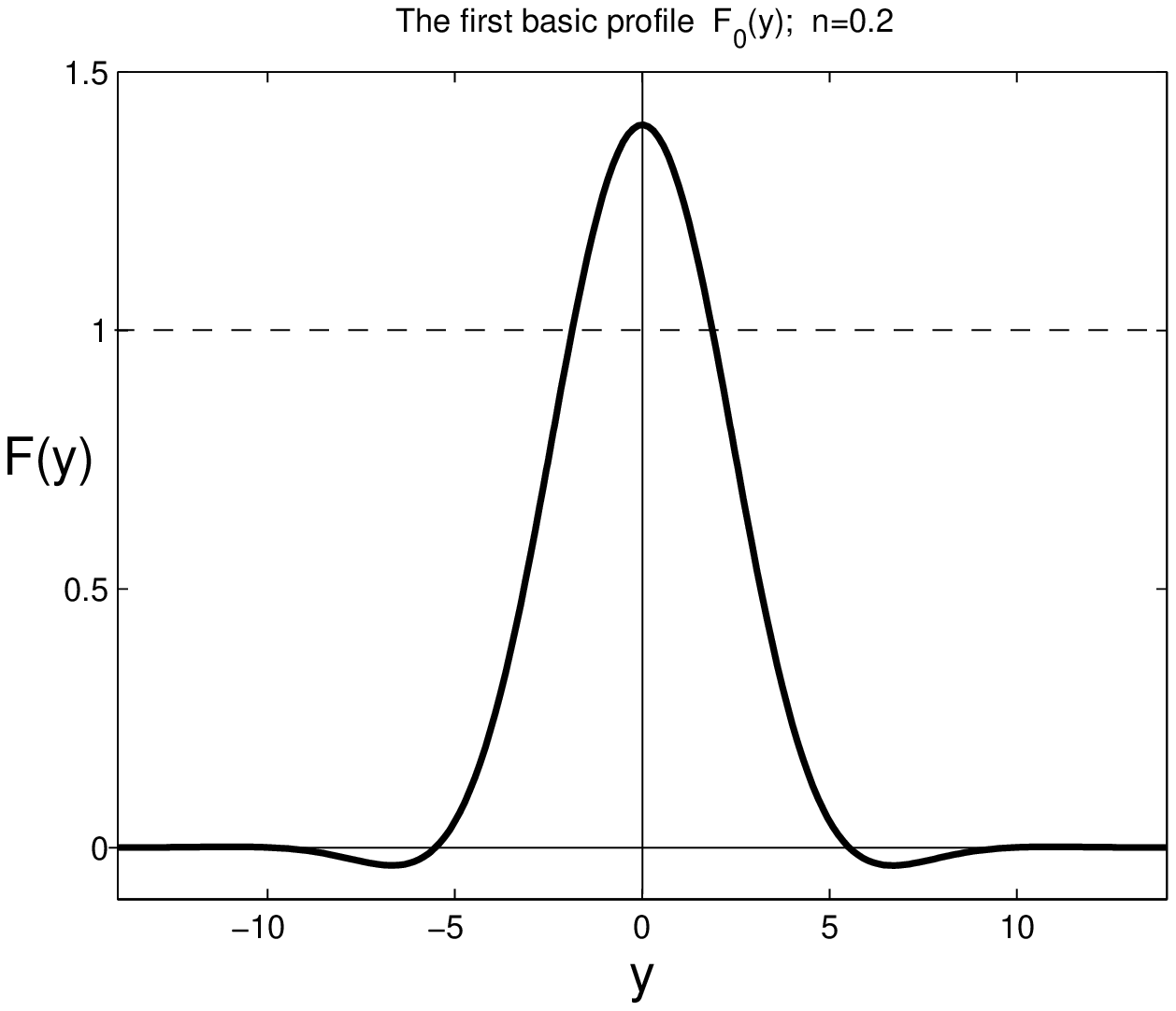}
} \subfigure[$F_1(y)$]{
\includegraphics[scale=0.52]{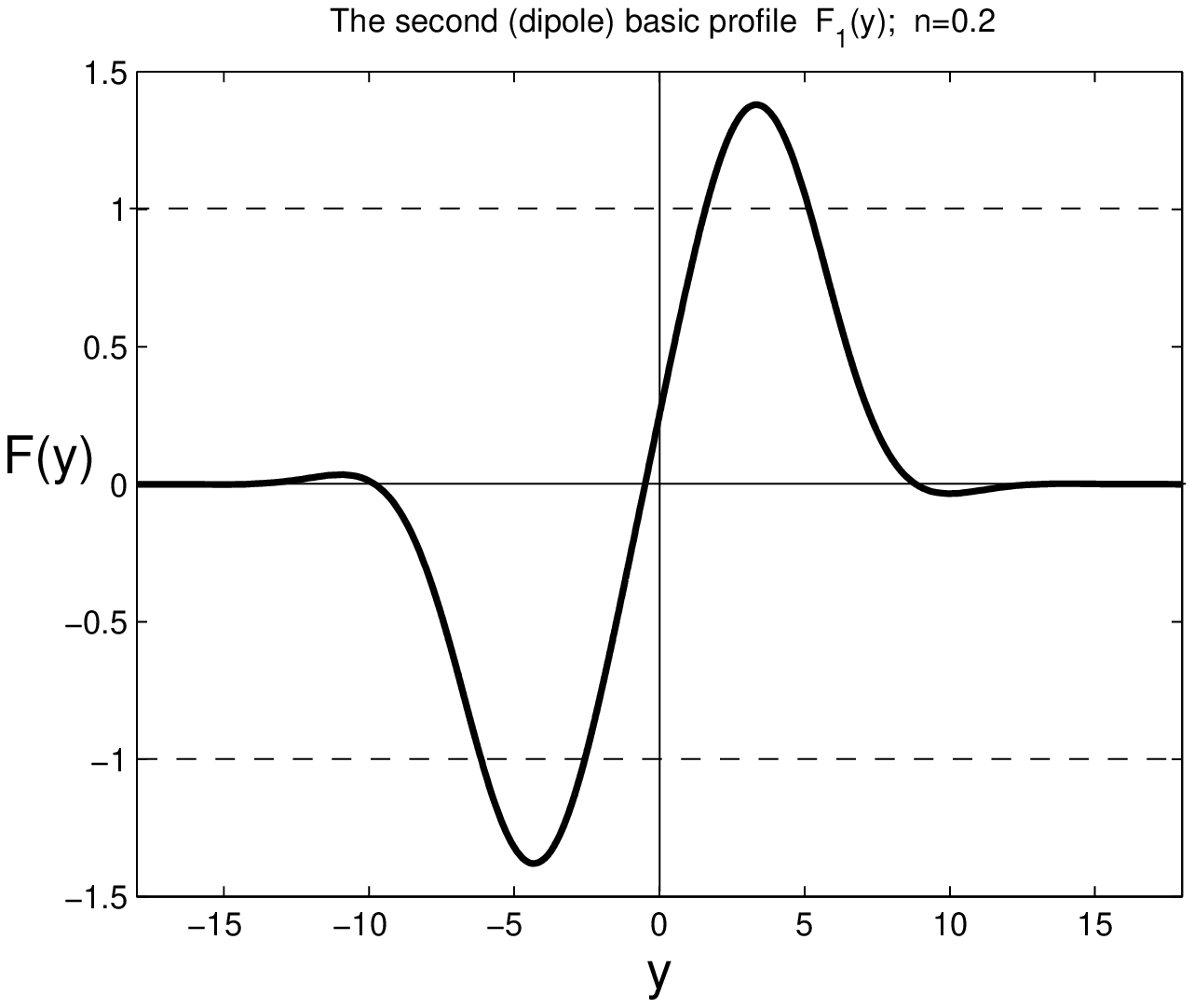}
} \subfigure[$F_2(y)$]{
\includegraphics[scale=0.52]{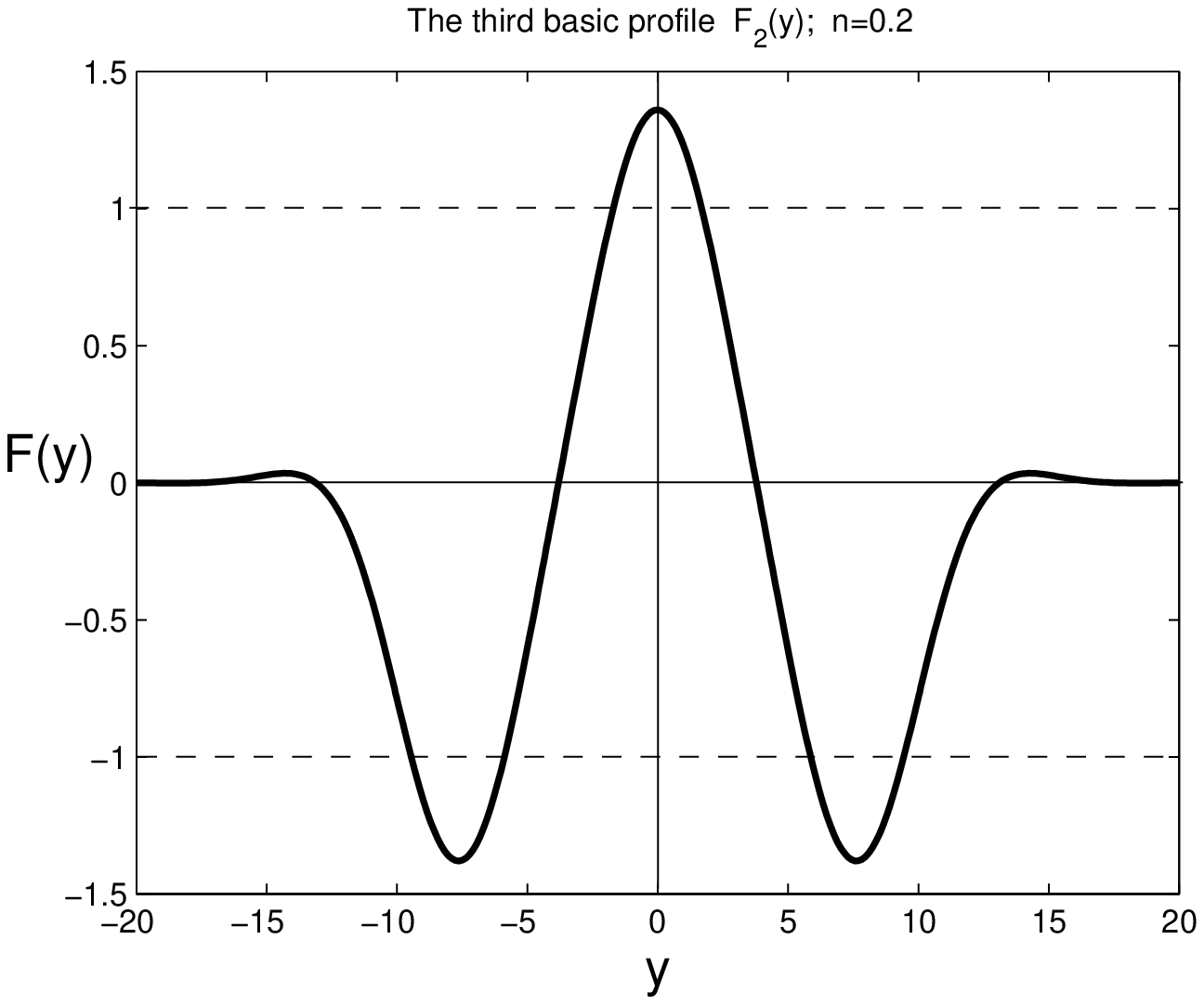}
} \subfigure[$F_3(y)$]{
\includegraphics[scale=0.52]{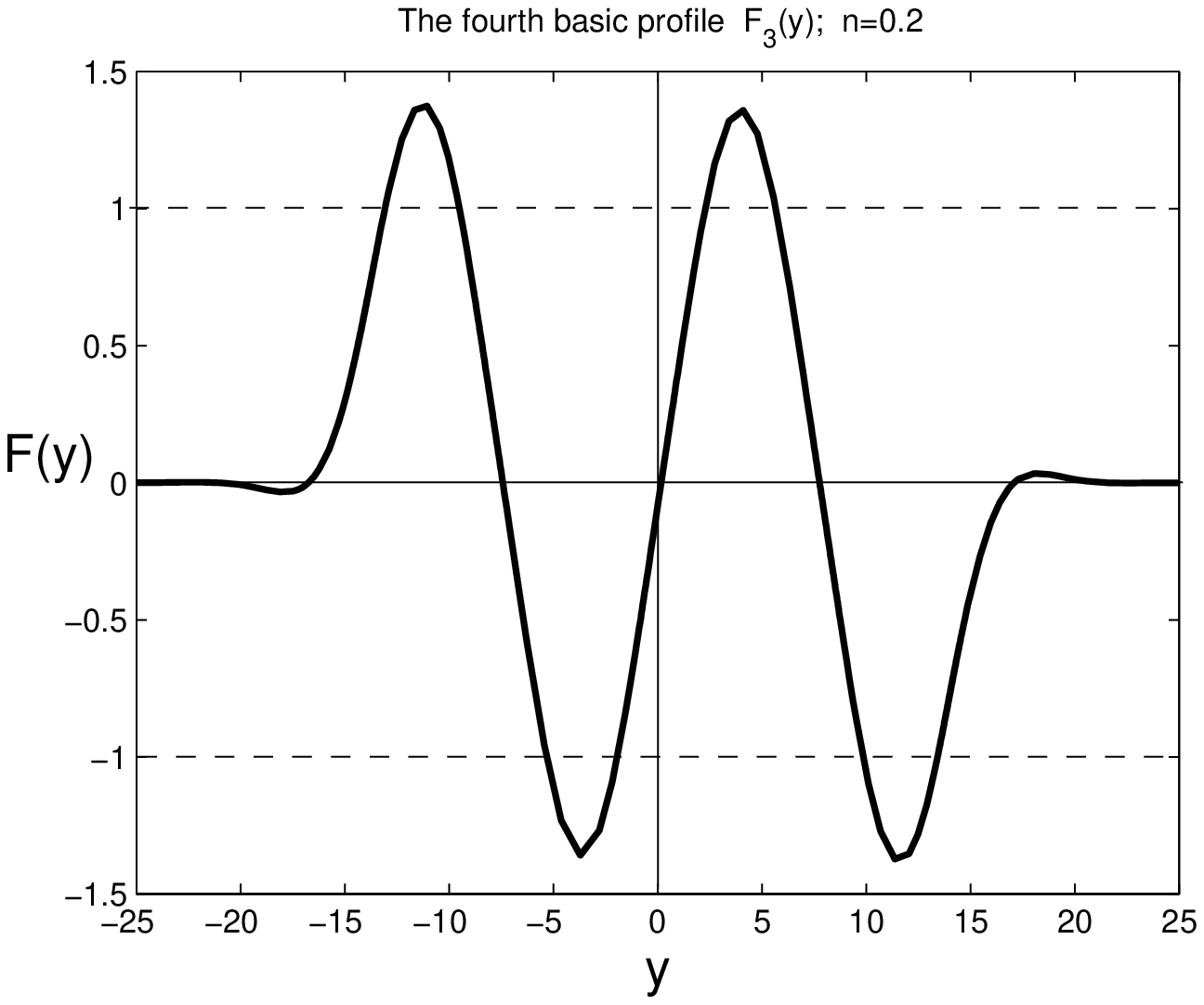}
} \subfigure[$F_4(y)$]{
\includegraphics[scale=0.52]{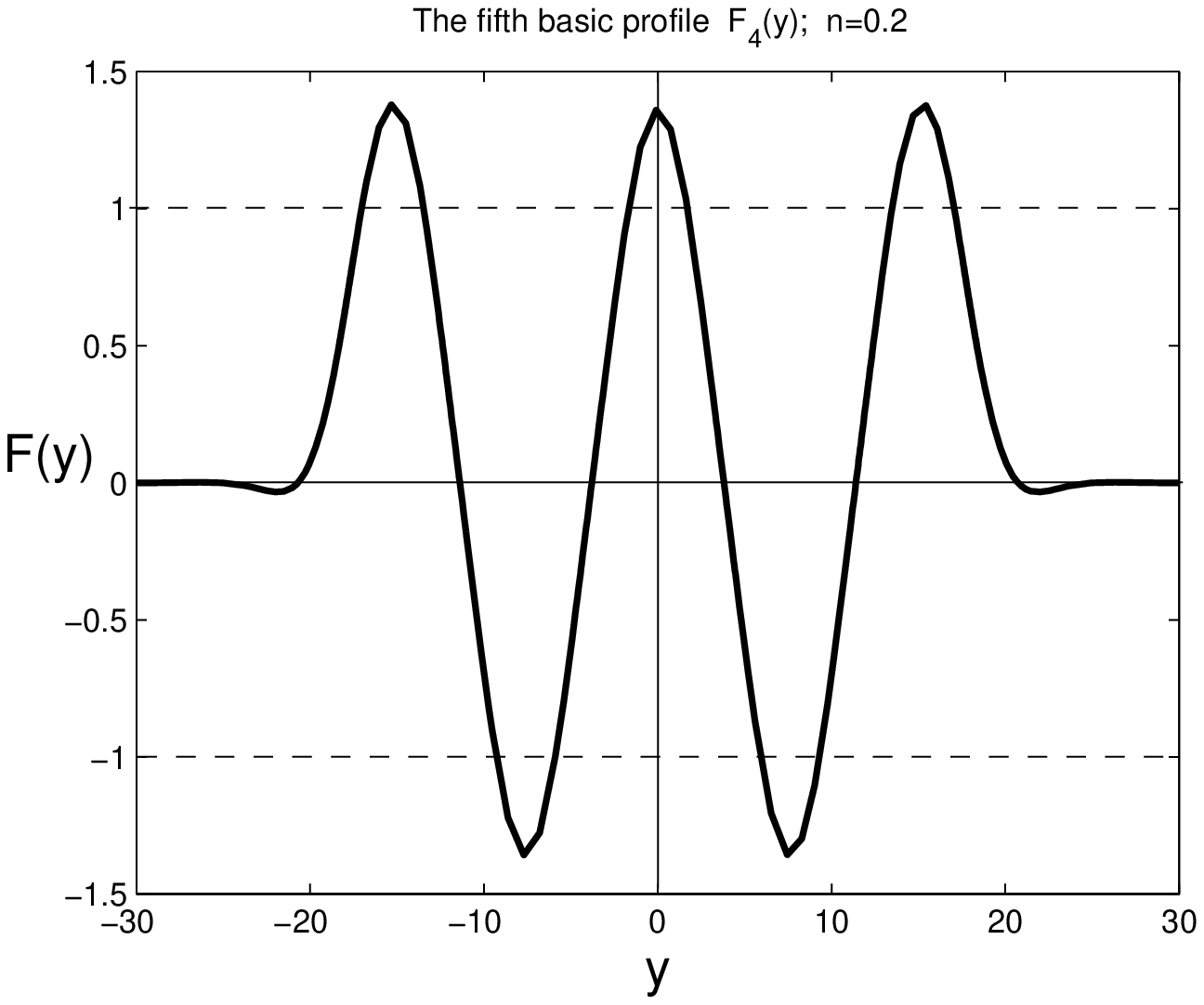}
} \subfigure[$F_5(y)$]{
\includegraphics[scale=0.52]{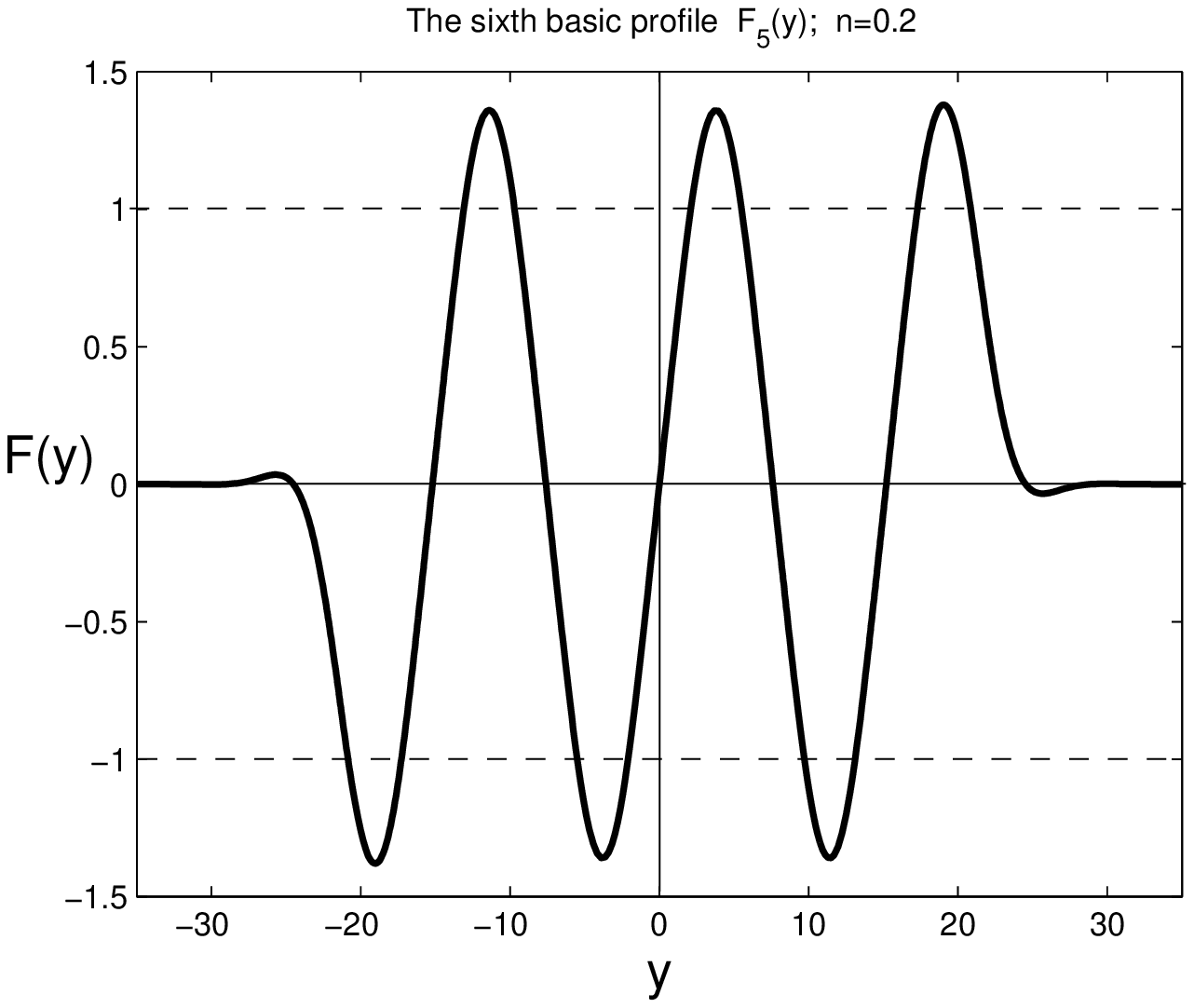}
}
 \vskip -.2cm
\caption{\rm\small The first six patterns of the basic family
$\{F_l\}$ of the ODE (\ref{1}) for $n=0.2$.}
 \label{G4}
\end{figure}

 In Figure \ref{G4}, we show the basic family denoted by
 $$
 \{F_l, \, l=0,1,2,...\}
  $$
   of solutions of (\ref{1}) for $n=0.2$. This family is connected
    with the application of  L--S and fibering theory;
   see
   \cite{GMPSob, GMPSobI}.
Each profile $F_l(y)$ has  $l+1$ ``dominant" extrema and $l$
``transversal" (not from the  tail) zeros; see \cite[\S~5]{GMPSob},
\cite[\S~5]{GMPSobI},
and \cite[\S~4]{GHUni} for further details. It is important that
 $$
\mbox{all the internal zeros of $F_l(y)$ are  {\em transversal},}
 $$
excluding the oscillatory end points of the support.
  In other words, each profile $F_l$ is
approximately obtained by a simple ``interaction" (gluing
together) of $l+1$ copies of the first pattern $\pm F_0$ taking
with necessary signs. Such a gluing of oscillatory tails  is
illustrated in Figures \ref{G6} and \ref{G7}; see also further
comments below. There is some analytic evidence \cite[\S~5]{GMPSobI}
that exactly this basic family $\{F_l\}$ is obtained by the
classic Lusternik--Schnirel'man construction of critical points of
the reduced functional (\ref{f4}). A rigorous definition of gluing
assumes formation of all the internal transversal zeros while the
outer ones at the end point of the support are the only ones that
remain oscillatory according to the behaviour (\ref{le3}) with the
periodic orbit $\varphi=\varphi_*(\ln(y_0-y))$. Some question of
global behaviour of such patterns $F_l(y)$ for large $l$ remain
open, \cite{GMPSob, GMPSobI}.

\ssk

Let us  forget for a moment about the complicated oscillatory
structure of solutions near interfaces, where an infinite number
of extrema and zeros occur. Then the dominant geometry of profiles
in Figure \ref{G4} looks like it approximately obeys Sturm's
classic zero set property, which is true rigorously for the case
$m=1$ only, i.e., for the second-order ODE
 \be
 \label{4.4}
  \mbox{$
 F''=-F + \big|F \big|^{-\frac n{n+1}}F \quad \mbox{in} \quad \re.
 $}
  \ee
For (\ref{4.4}), the basic family $\{F_l\}$ is explicitly
constructed by direct gluing together simple patterns $\pm F_0$
given explicitly; see \cite[p.~168]{GSVR}. Therefore, each $F_l$
consists of precisely $l+1$ patterns  (with signs $\pm F_0$), so
that Sturm's property is clearly true
    by direct  application of L--S
    category theory.



 \subsection{Countable family of $\{F_0,F_0\}$-gluing}

Further patterns to be introduced do not exhibit as clear a
``dominated" Sturm property and are associated with a double
fibering technique where both the {\em Cartesian} and {\em
spherical} representations of critical points are used; see
\cite[\S~3]{GMPSob}, \cite[\S~3]{GMPSobI}. Let us present some explanations.

\begin{figure}
 \centering \subfigure[profiles]{
\includegraphics[scale=0.52]{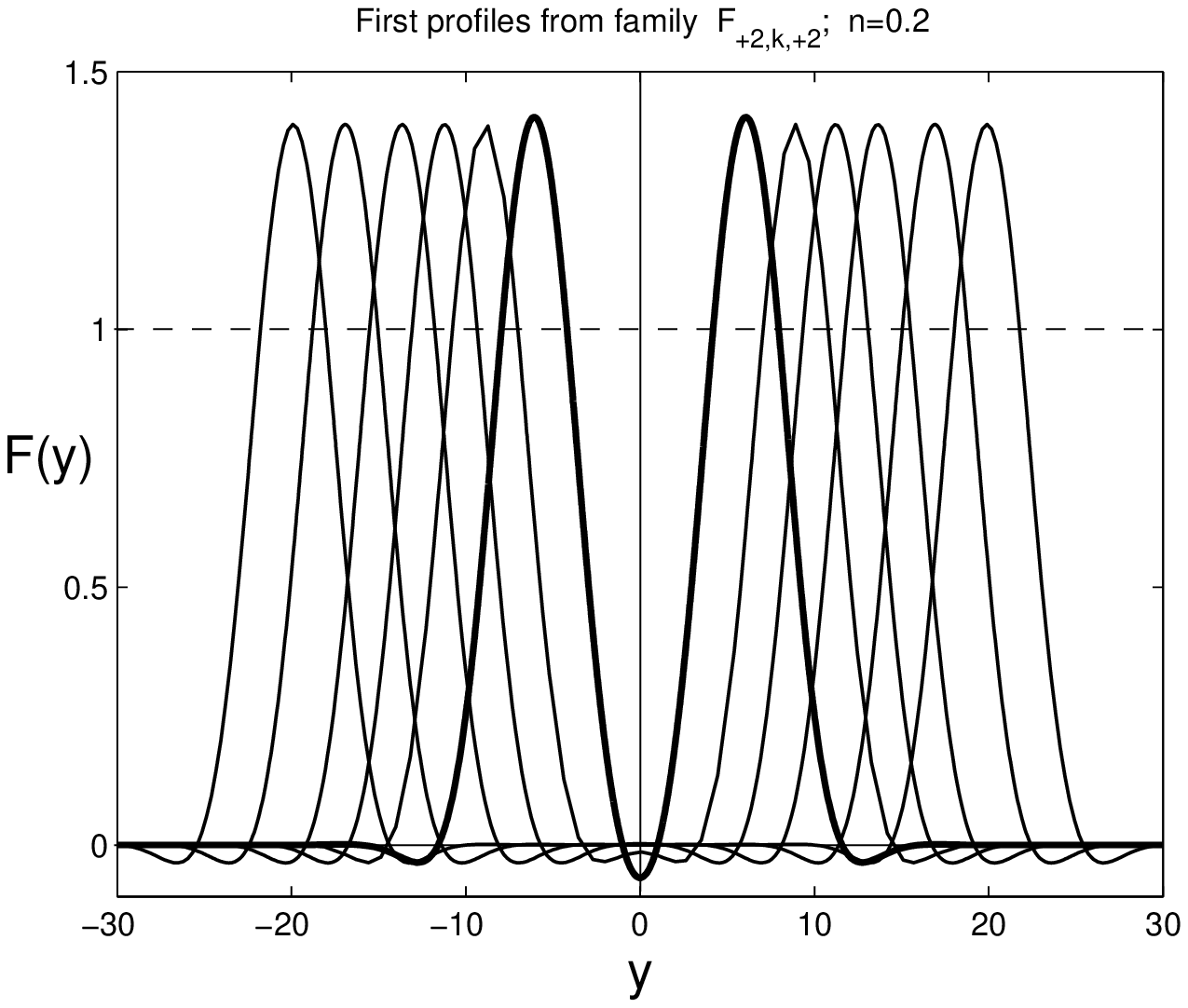}
} \subfigure[zero structure, enlarged]{
\includegraphics[scale=0.52]{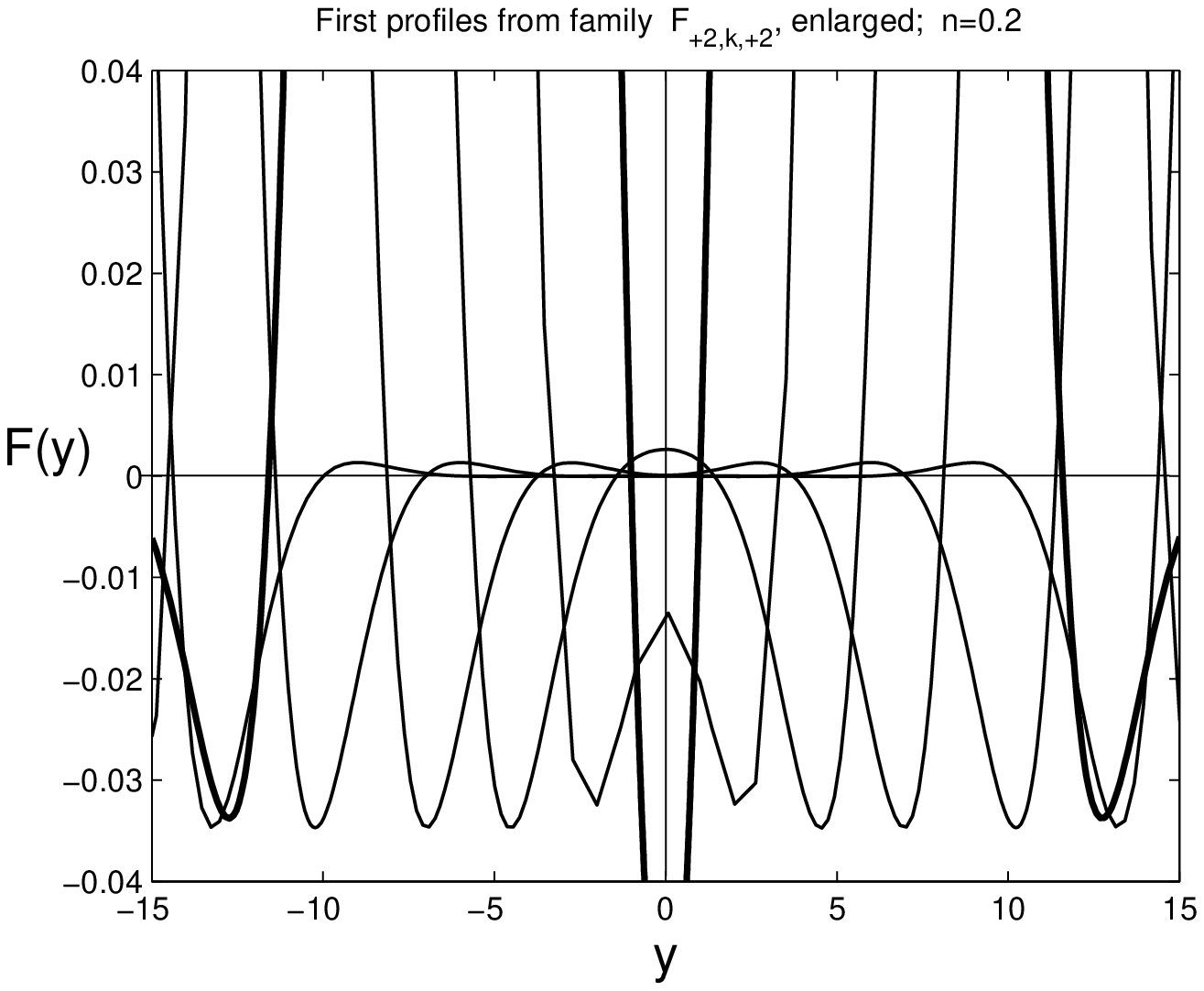}               
}
 \vskip -.2cm
\caption{\rm\small First six patterns from the family
$\{F_{+2,k,+2}\}$ of the  $\{F_0,F_0\}$-gluing;  $n=0.2$.}
 \label{G6}
\end{figure}

The nonlinear interaction of the two first
 patterns $F_0(y)$ leads to a new family of profiles.
In Figure \ref{G6} for $n=0.2$, we show the first six profiles
from this family denoted by $\{F_{+2,k,+2}\}$. In the refined zero
structure of the last profile in (b), we already see some
numerical effects of rather rough meshes, which again do not deny
the sufficient overall quality of numerics.
 In each function
$F_{+2,k,+2}$ the multiindex
 $$
 \s=\{+2,k,+2\},
  $$
   from left to right
denotes: +2 means two intersections with the equilibrium +1, then
next $k$ intersections with zero, and final +2 stands again for
two intersections with +1. Later on, we will use such a multiindex
notation to classify other  patterns obtained.

 As a general rule, we point out again that any finite gluing of a pair of patterns
 $\pm F_0(y)$ actually means that all internal zeros become
 transversal. Note that this and all the Figures involved are not enough to explain the essence
 of this complicated and mathematically not fully understood
 procedure for non-homotopic variational problems
 \cite{GMPSobII}. The resulting patterns  have zeros of infinite
 order only at the end points of its support.



In view of the infinite oscillatory character of $F_0(y)$ at the
interfaces, we expect that the family $\{F_{+2,k,+2}\}$ is
countable, and such functions exist for any even $k=0,2,4,...\, $.
Then $k=+\infty$ corresponds to the non-interacting pair
  \be
 \label{F0+}
 F_0(y+ y_0) + F_0(y-y_0), \quad \mbox{where} \,\,\,\,\,{\rm supp}\,
 F_0(y) = [-y_0,y_0].
  \ee

It is expected that there exist various triple $\{F_0,F_0,F_0\}$
and any multiple interactions $\{F_0,...,F_0\}$ of $k$ single
profiles, with different distributions of zeros between any pair
of neighbours (proof is an open problem).

 \subsection{Countable family of  $\{-F_0,F_0\}$-gluing}

\begin{figure}
 \centering
\includegraphics[scale=0.6]{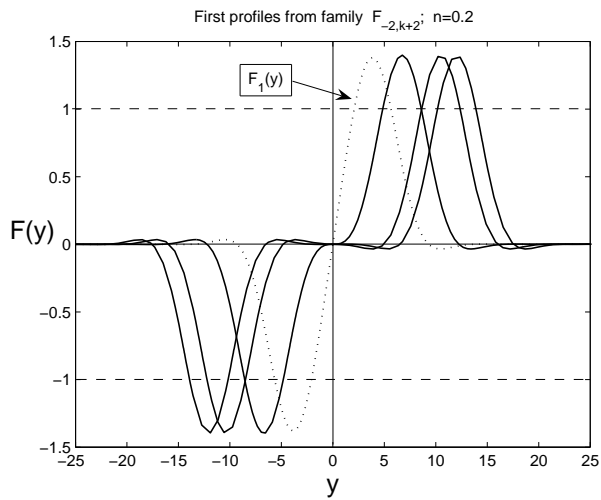}
 \vskip -.3cm
\caption{\rm\small First four patterns from the family
$\{F_{-2,k,+2}\}$ of the  $\{-F_0,F_0\}$-gluing;  $n=0.2$.}
   \vskip -.3cm
 \label{G7}
\end{figure}

 We now describe the interaction of $-F_0(y)$ with
  $F_0(y)$.
In Figure \ref{G7}  for the case $n=0.2$ (which is convenient in
terms of rather fast convergence of the numerical method
employed), we show the first profiles from this
family denoted by $\{F_{-2,k,+2}\}$, where for  
 the
multiindex  $\s=\{-2,k,+2\}$,
  the first number $-2$ means two
intersections with the equilibrium $-1$, etc.  It
can be seen that the first two profiles belong to the same class
 $
 F_{-2,1,2},
 $
i.e., both have a single zero for $y \approx 0$.
  The last solution shown is $F_{-2,5,+2}$.
 Again, we expect that the family $\{F_{-2,k,+2}\}$ is countable,
and such functions exist for any odd $k=1,3,5,...$, and
$k=+\infty$ corresponds to the non-interacting pair
 \be
 \label{F0-}
  \mbox{$
- F_0(y+ y_0) + F_0(y-y_0).
  $}
  \ee
   We expect that
there exist families of an arbitrary number of gluing $\{\pm F_0,
\pm F_0,..., \pm F_0\}$ consisting of any $k \ge 2$ members (again
an open problem).




\subsection{Periodic solutions in $\re$ as new types of oscillations about $\pm 1$}

\begin{figure}
 \centering \subfigure[about $F \equiv 1$: $F(0)=0.4135$]{
\includegraphics[scale=0.52]{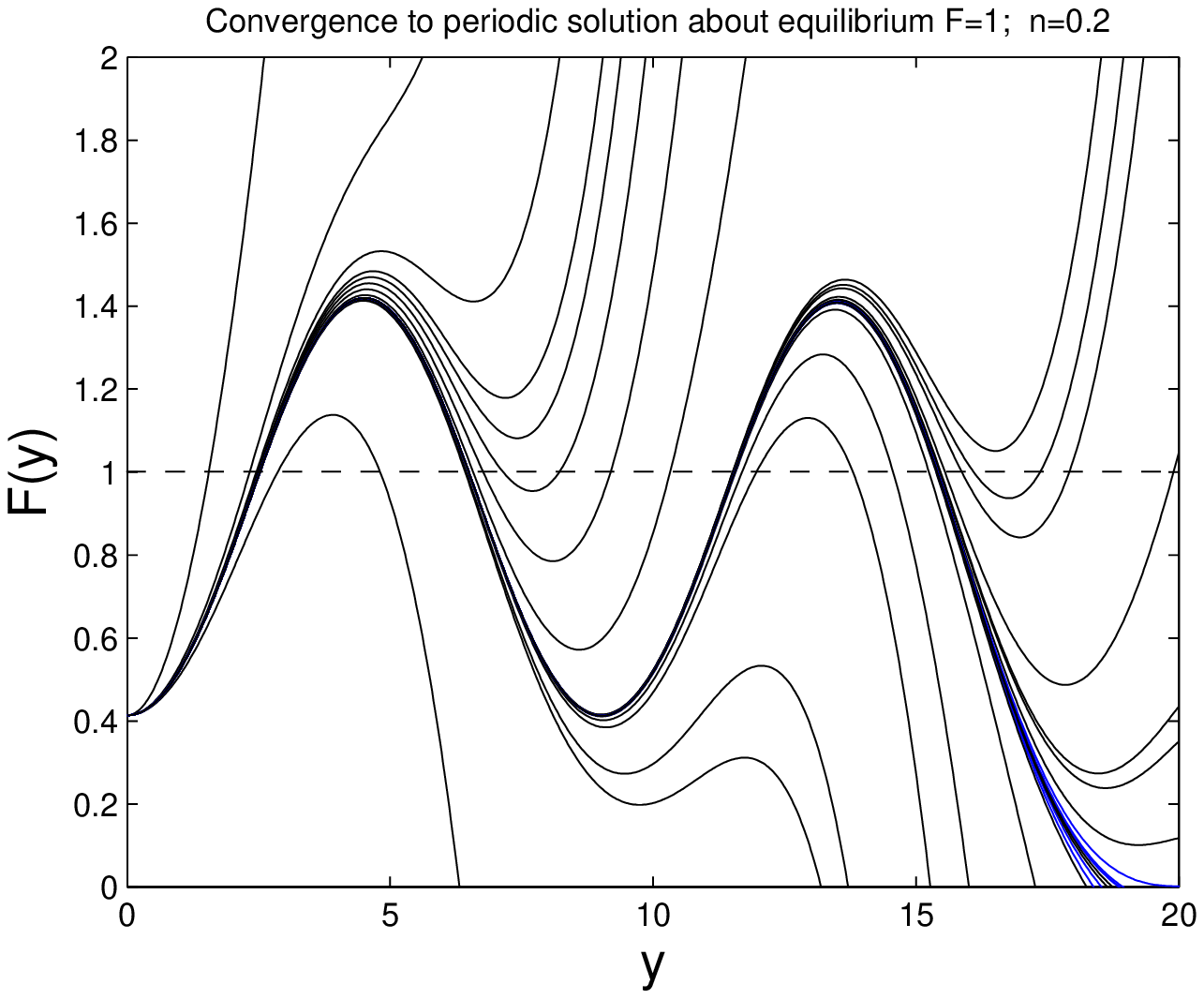}
} \subfigure[about $F \equiv -1$: $F(0)=1.4085$]{
\includegraphics[scale=0.52]{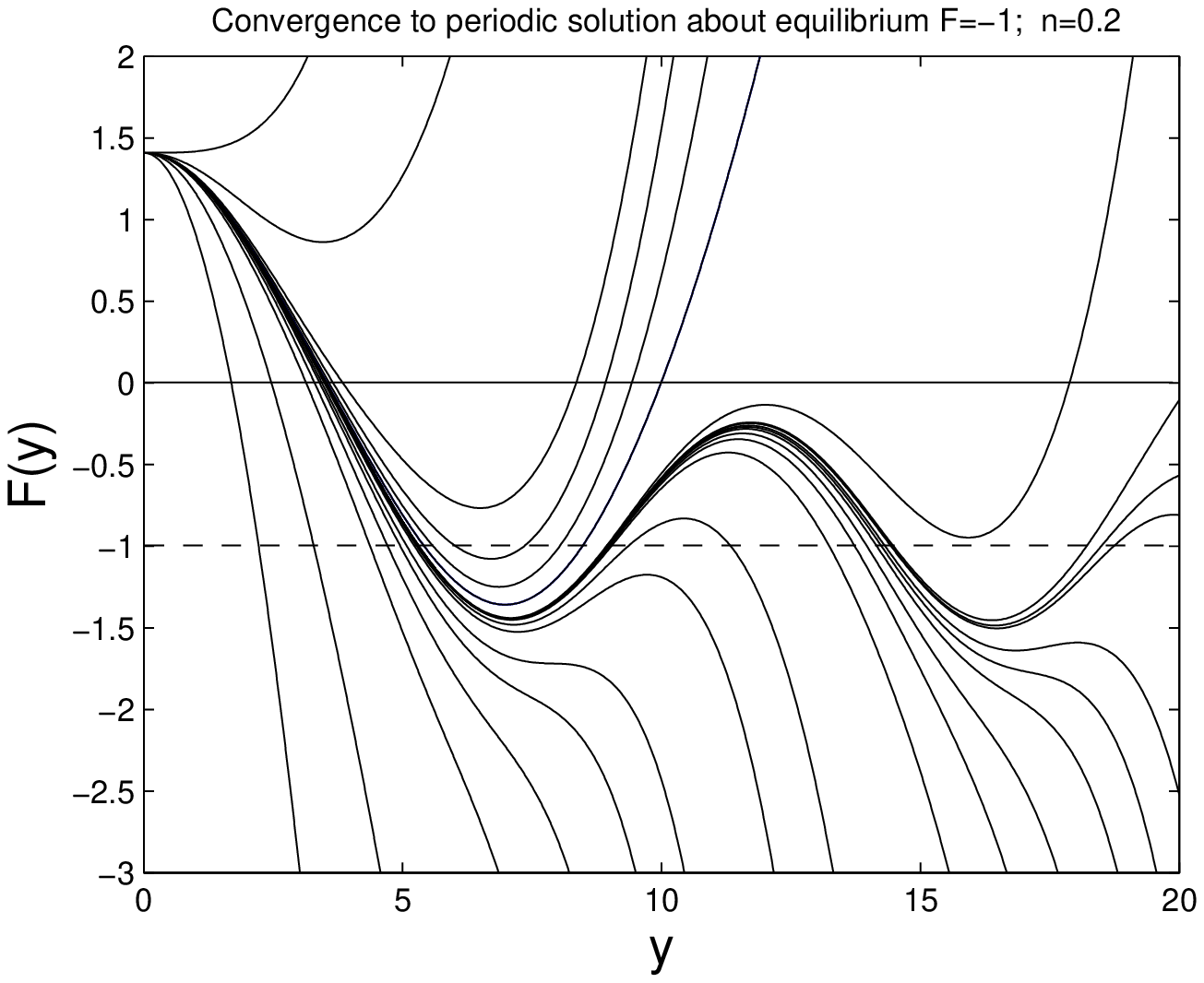}               
}
 \vskip -.2cm
\caption{\rm\small Examples of convergence to periodic solutions
of the ODE (\ref{1}) for $n=0.2$; about $F \equiv 1$ (a) and $F
\equiv -1$ (b).}
 \label{GP1}
\end{figure}

Before introducing new types of patterns, we need to describe
other non-compactly supported solutions in $\re$. As a variational
problem, equation (\ref{1}) admits  infinitely many periodic
solutions; see e.g., \cite[Ch.~8]{MitPoh}.  Figure \ref{GP1} for
$n=0.2$ reveals unstable periodic solutions obtained by shooting
from the origin with various Cauchy data at $y=0$. In (b), the
periodic orbit $F_*(y)$ is oscillating about the equilibrium $F
\equiv -1$.
 It turns out that precisely the periodic orbit $F_*(y)$ in (a) with the
 range
 \be
 \label{**1}
 \min \, F_*(y) =0.4135..., \quad \max \, F_*(y)= 1.4085...
  \quad (n=0.2)
 \ee
  plays an important part in the
 construction of other families of compactly supported
patterns. Namely, all the variety of solutions of (\ref{1}) that
have oscillations about equilibria $\pm 1$ are  close to $\pm
F_*(y)$ there.


\subsection{Family $\{F_{+2k}\}$}

\begin{figure}
 \centering
\includegraphics[scale=0.65]{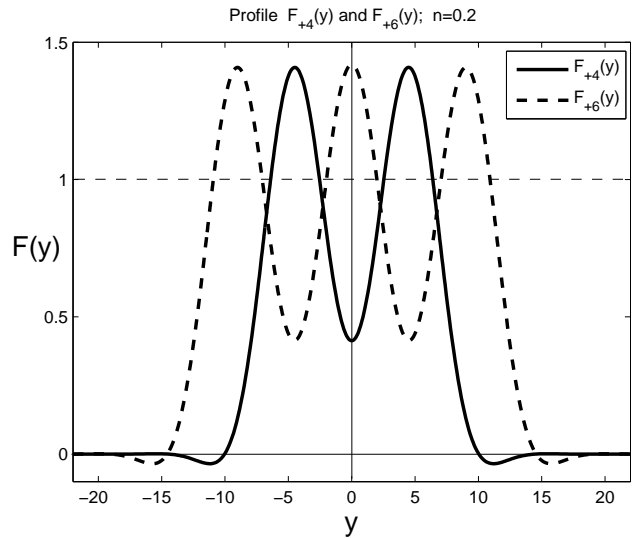}
 \vskip -.4cm
\caption{\rm\small The first profiles $F_{+4}(y)$ and $F_{+6}(y)$
from the family $\{F_{+2k}, \, k \ge 2\}$; $n=0.2$.}
   \vskip -.3cm
 \label{FF4p}
\end{figure}

\begin{figure}
\centering
{
\includegraphics[scale=0.65]{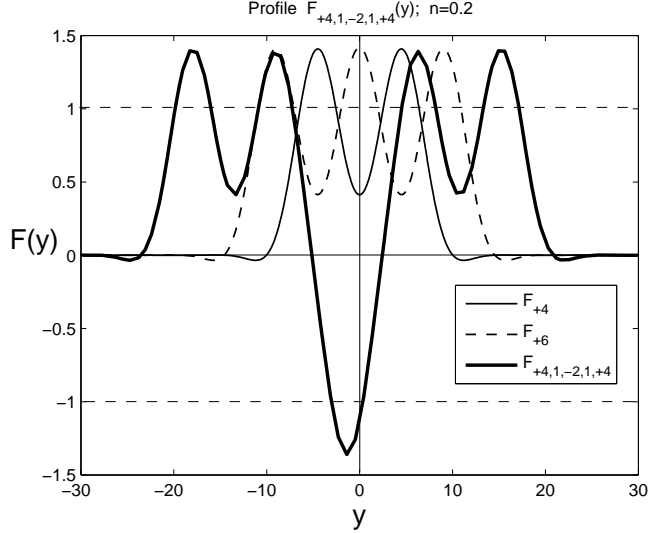}
}
 \vskip -.2cm
\caption{\rm\small The profile $F_{+4,1,-2,1,+4}(y)$ (the boldface
line) from the family $F_{+k,l,-m,l,+k}$ of solutions of
(\ref{1}); $n=0.2$. For comparison, the profiles  $F_{+4}(y)$ and
$F_{+6}(y)$ from Figure \ref{FF4p} are presented.}
 \label{G8}
\end{figure}

Such functions $F_{+2k}$ for $k \ge 1$ have $2k$ intersections
with the single equilibrium +1  only and have a clear ``almost"
periodic structure of oscillations about. The number of
intersections denoted by $+2k$ gives an extra Strum index to such
a pattern.
In this notation, for $k=1$, we have
 $$
 F_{+2}=F_0.
  $$
  Two profiles $F_{+4}$ and $F_{+6}$ are  shown in Figure \ref{FF4p} for $n=0.2$.
The further profile $F_{+4,1,-2,1,+4}(y)$   comprising two
sub-structures $F_{+4}$ from the family $\{F_{+2k}\}$ is shown in
Figure \ref{G8} by the boldface line.

\subsection{More complicated patterns: towards  chaotic structures}

\begin{figure}
 \centering
\includegraphics[scale=0.65]{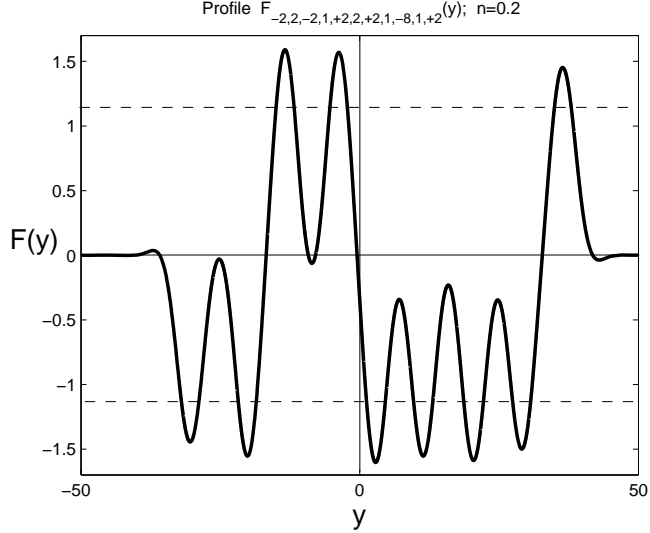}
 \vskip -.4cm
\caption{\rm\small A complicated pattern $F_{\s}(y)$ of the index
(\ref{ch1}) for (\ref{1}) for $n=0.2$.}
   \vskip -.3cm
 \label{FC2}
\end{figure}

By combining the above rather simple families of patterns, we
claim that a pattern (more precisely, a class of patterns) with an
arbitrary admissible   multiindex of any length
 \be
 \label{mm1}
 \s=\{\pm \s_1, \s_2, \pm \s_3, \s_4,..., \pm \s_l\}
  \ee
  can be constructed.
  For example, in Figure \ref{FC2}, a single  complicated pattern
  with
   \be
   \label{ch1}
   \s=\{-2,2,-2,1+2,2,+2,1,-8,1,+2\}
    \ee
    is given.
   The computation (the convergence is rather slow for this type of
   $p$-Laplacian operators) is performed for $n=0.2$ as usual.
  We must admit that, as it is seen from Figure \ref{FC2},  the iterations have not been properly
    converged within the parameter range (\ref{eps1}), but we can guarantee the
     convergence with the accuracy
    at least $\sim 10^{-1}$.
 This is not that bad, since such patterns are not structurally stable
  and have multi-dimensional unstable manifolds, so the
  convergence {\em must be} extremely slow. It is worth mentioning that,
     using {\tt bvp4c} solver, this computation took
     a few hours with the maximal number of 75000
   points on the interval $(-50,50)$.
     Nevertheless, regardless such
    a lack of accuracy, we are sure that such complicated critical point  profiles
    really exist, since we have seen a lot of those in other
    similar (and simpler numerically)
   higher-order variational problems \cite{GMPSob, GMPSobI} that were
    not associated with such awkward and strongly degenerate operator as $p$-Laplacian
    ones.
   Special more conservative and divergent numerical techniques
   are necessary for tackling higher-order $p$-Laplacian operators in the
   ODEs, but here we demonstrate what an average (numerically, non-professional)
    PDE user can extract from standard {\tt MatLab}
   codes.  Theoretically, via the L--S/fibering theory, all those
    patterns are well defined.

We claim that the multiindex (\ref{mm1}) can be rather arbitrary
taking   finite parts of any non-periodic ``fraction". Actually,
this means {\em chaotic features}
 of the whole family of solutions $\{F_\s\}$. In fact, there is no any exiting news in such
 a chaotic proclamation: one can see
 that even the basic simple countable family $\{F_l\}$ is indeed {\em
 chaotic}, since  the choice of the sequence of elementary profiles
 $\pm F_0$ in $F_l$ for $l \gg 1$ can be {\em arbitrary} long with an {\em arbitrary} sequence
 $\{\pm\}$ of sign changes, thus exhibiting no finite periodic order in the index $\s$.
 These chaotic types of behaviour are known for other simpler fourth-
 order ODEs with
coercive operators and definite {\em homotopic} features,
\cite[p.~198]{PelTroy}.


A variety of  complicated patterns of these types for different
variational problems associated with the PME-type nonlinearities
in (\ref{pp1}) can be found in \cite[\S~5.6]{GMPSob},
\cite[\S~5.6]{GMPSobI} (see also \S~6 in \cite{GMPSobI} therein
for sixth- and eighth-order models) and in \cite[\S~5]{GPME}. The
convergence of standard numerical methods for these problems are
much faster, with $\e$, Tols up to $\sim 10^{-10}$.

\section{Single point blow-up for $p>n+1$: P- and Q-type profiles}
 \label{SectLS1}

 We now return to the similarity  ODE (\ref{f11E}) in the case $p>n+1$, which,
in view of the spatial rescaled variable $y$ in (\ref{RVarsE}),
corresponds to single point blow-up. It is key that (\ref{f11E})
for $p \not = n+1$ {\em is not variational}. Formally, solutions
of (\ref{f11E}) can be traced out by shooting and matching
procedures, which are too complicated. Instead,  we will use a
continuation in parameters approach, which allows us to predict
solutions by using those in the variational case $p=n+1$.

\subsection{Asymptotics at infinity and single point blow-up}

We begin with simpler asymptotics of the solutions of (\ref{f11E})
as $y \to + \infty$. Unlike the previous case of regional blow-up
for $p=n+1$, for $p>n+1$, equation  (\ref{f11E}) admits
non-compactly supported solutions with the following behaviour:
 \be
 \label{as1}
 f(y) \sim (C_0 y^\g +...) + (C_1 \, {\mathrm e}^{-b_0
 y^\nu}+...) \quad \mbox{as} \quad y \to + \iy,
 \ee
 where $C_0 \not = 0$ and $C_1 \in \re$ are arbitrary constants
 and
 $$
 \mbox{$
  \g = - \frac {2(n+2)}{p-(n+1)}<0, \quad
 \nu= \frac{2(n+2)(p-1)}{3[p-(n+1)]}>0, \quad
   b_0= \frac 1 \nu \big[ \frac \b{(n+1) \g^n(\g-1)^n} \, C_0^{-n}\big]^{\frac
   13}> 0.
   $}
   $$
The first term in (\ref{as1}) represents an ``analytic" part (can
be truly analytic for some parameters) of the expansion, while the
second one gives the essentially ``non-analytic" part.
 Such a structure in (\ref{as1}) is usual for saddle-node-type
 equilibria \cite[p.~311]{Perko}, but, f
 or the fourth-order
 ODE (\ref{f11E}), this expansion does not admit a simple
 phase-plane interpretation. However,  existence of such asymptotic
 expansions can be justified by fixed point arguments, which
 becomes quite a technical issue and is not done here.

 One can see passing to the limit $t \to T^-$  in (\ref{RVarsE})
that the first term in the asymptotic expansion (\ref{as1}) gives
the following {\em final-time profile} of this single point
blow-up for even patterns $f=f(|y|)$:
 \be
 \label{ft1}
 u_S(x,T^-)= C_0 |x|^{-\frac {2(n+2)}{p-(n+1)}} < \infty \quad \mbox{for
 all \, $x \not = 0$.}
  \ee

Returning to the asymptotic expansion, we conclude that
(\ref{as1}) represents
 \be
 \label{as2}
 \mbox{a 2D asymptotic family (bundle) of solutions.}
  \ee
Hence, the family (\ref{as1}) is well suitable for matching with
also two symmetry conditions at the origin (\ref{BCs}), so we
expect not more than  a countable set of solutions. For first
patterns, we keep the same notation $f_l(y)$  as in Section
\ref{SectS} for $p=n+1$.


\subsection{Oscillatory behaviour about constant equilibrium
$f_*$}

In order to predict the multiplicity of solutions of (\ref{f11E}),
we need to study its oscillatory properties.
 To this end, we perform the linearization about the constant equilibrium $f_*$ in (\ref{g1N})
  of the ODE
(\ref{f11E}),
 \be
  \label{lin11}
  f=f_*+Y
   \ee
   formally
assuming  that $|Y| \ll 1$ on some bounded intervals. This yields
the ``linearized" nonlinear equation
 \be
 \label{BB1}
  \mbox{$
 {\bf B}_n(Y) \equiv -(|Y''|^n Y'')'' - \b Y' y + Y=0 \quad
  \big(\,\b =
\frac{p-(n+1)}{2(n+2)(p-1)}\,\big).
 $}
 \ee

We are going to study oscillatory, sign-changing properties of
solutions of (\ref{BB1}) for various $n>0$.
 Notice that the ``linearized" ODE (\ref{BB1}) remains a difficult
fourth-order equation. Indeed, in view of the invariance with
respect to the group of scalings
 $$
 Y \mapsto \e^{\frac{2(n+2)}n} Y, \quad y \mapsto \e \, y \quad
 (\e >0)
 $$
the transformation
 $$
 Y(y)= y^{\frac{2(n+2)}n} \varphi(s), \quad s= \ln y,
 $$
 reduces (\ref{BB1}) to an autonomous fourth-order ODE.
 Setting $P(\varphi)=\varphi'$ yields a third-order ODE,
 but further reductions are impossible. Thus, (\ref{BB1}) cannot be
 studied on the phase-plane in principle; cf. \cite{BuGa}, where, for (\ref{s1}), oscillatory
 analysis on the phase-plane is a convenient and exhaustive tool.

Therefore, we will need another further investigation of
(\ref{BB1}), and we begin with the following useful comment:

\ssk

\noi{\bf Linear case $n=0$.} Then the quasilinear operator ${\bf
B}_n$ in (\ref{BB1}) becomes linear,
 \be
 \label{BB2}
  \mbox{$
   {\bf B}_0 Y = 0, \quad \mbox{where} \quad
  {\bf B}_0 = - D_y^4  - \frac 14 \, y D_y + I \equiv {\bf B}^* +
  I \quad \big(\b= \frac 14 \big)
  $}
  \ee
  and ${\bf B}^*$ is the adjoint linear operator
  (\ref{ad1}). According to its point spectrum (\ref{spec1}),
  equation (\ref{BB1}) for $n=0$ has a non-oscillatory solution,
  being the eigenfunction $\psi^*_l(y)$ for $l=4$, i.e.,
   \be
   \label{BB3}
    \mbox{$
   Y(y) =\psi^*_4(y)= \frac 1{\sqrt{24}} \,(y^4+24) \quad (n=0),
 $}
    \ee
which is an example of a {\em non-oscillatory} solution.
Nevertheless, one can see from the operator (\ref{BB2}) that the
linear ODE (\ref{BB1}) for $n=0$ has other oscillatory solutions
with an increasing envelope as $y \to +\infty$; see below.

\ssk



\begin{figure}
\centering \subfigure[ $n=0$ ]{
\includegraphics[scale=0.52]{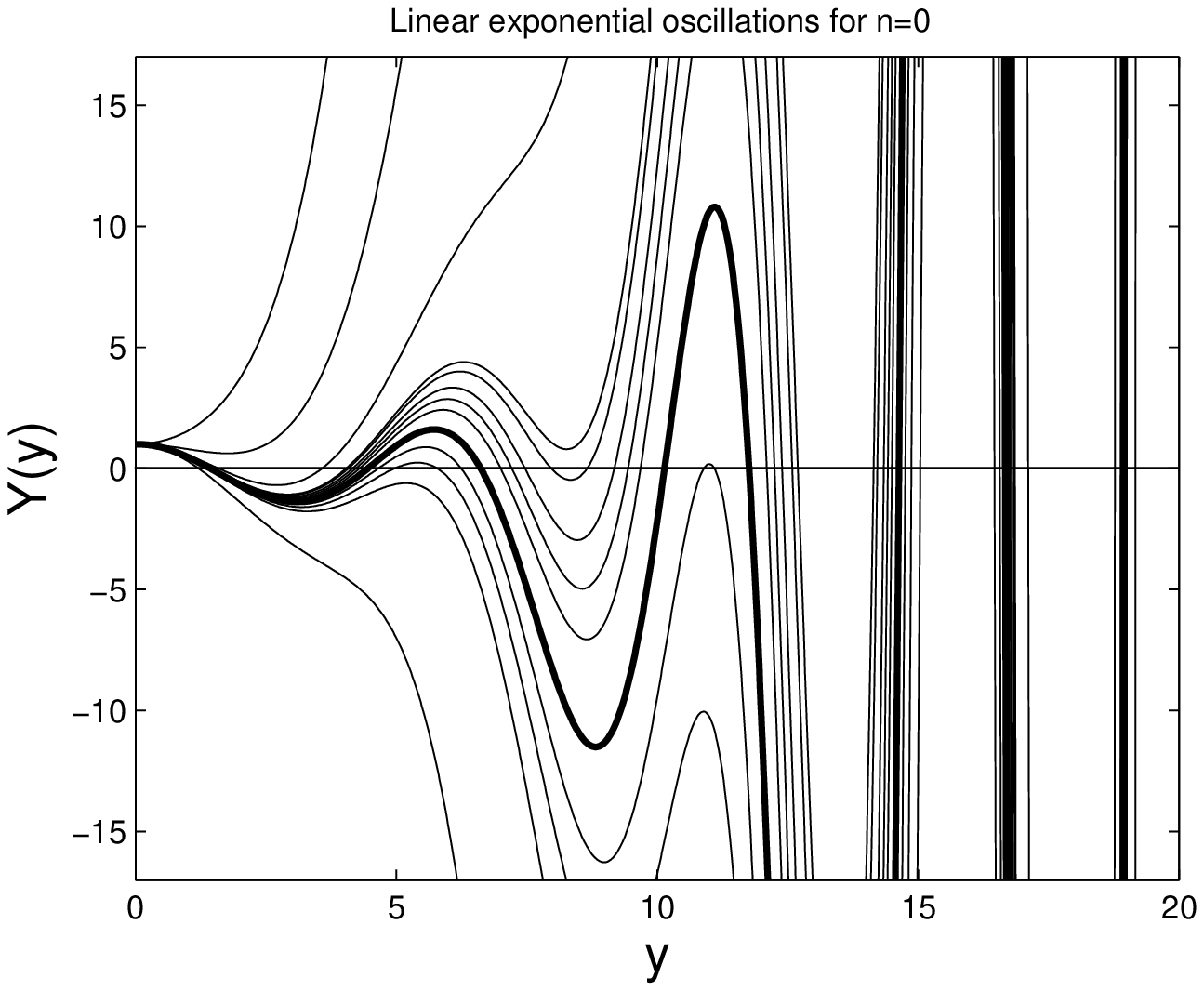}
} \subfigure[$p=n+1: \, n=0.2$]{
\includegraphics[scale=0.52]{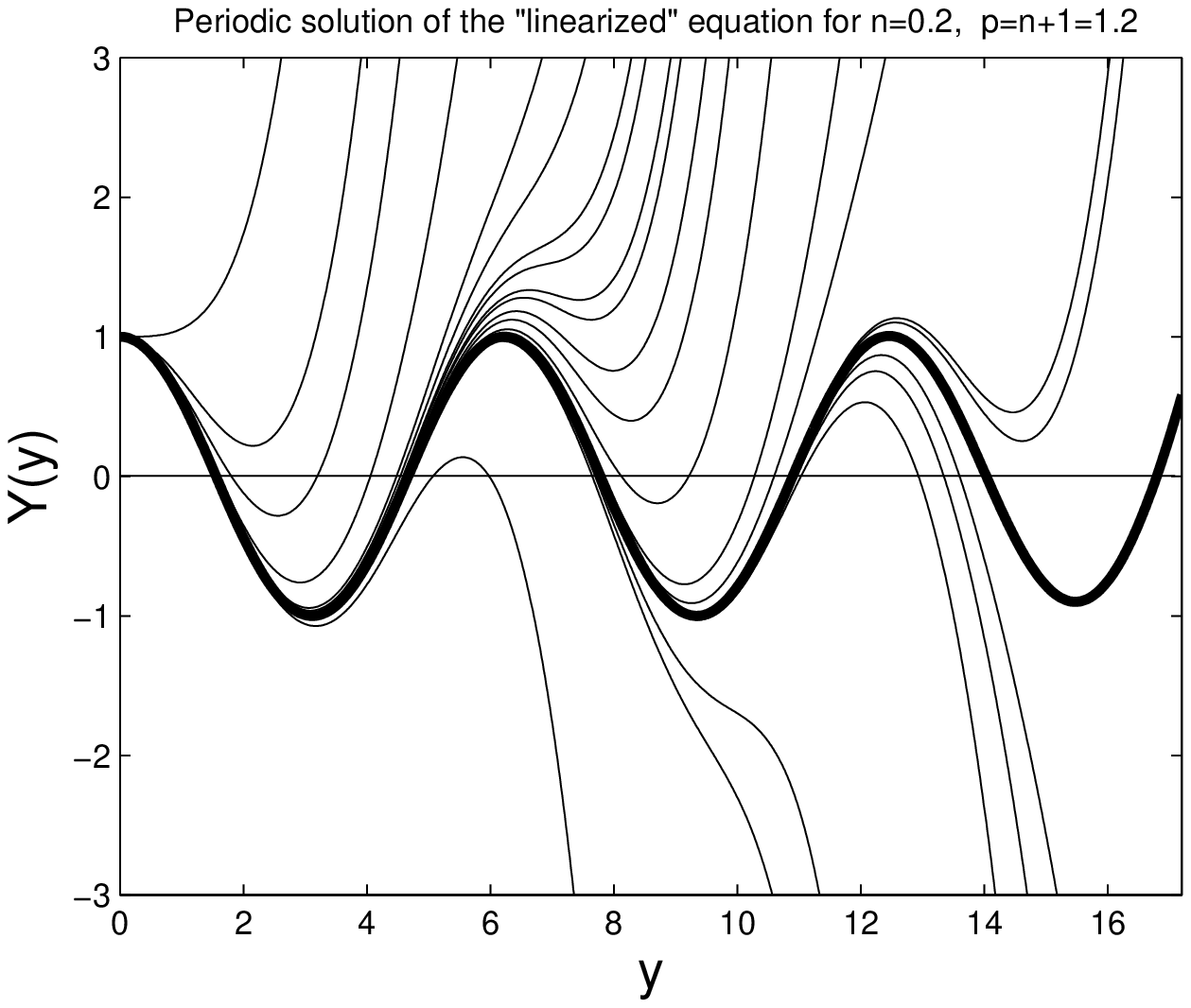}
} \subfigure[$p>n+1: \, p=1.4,\,n=0.2$]{
\includegraphics[scale=0.52]{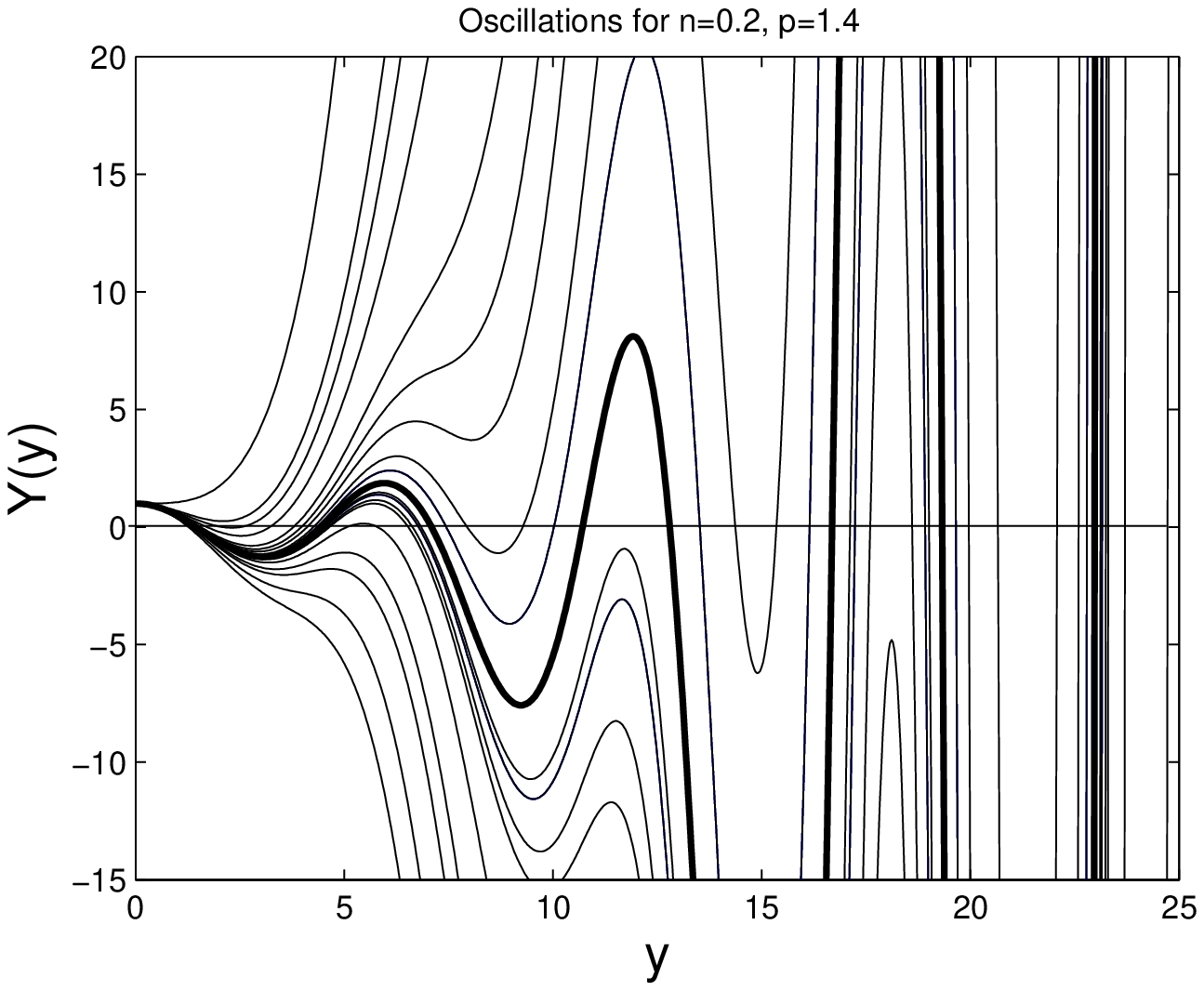}
}
\subfigure[$p<n+1: \, p=1.1, \, n=0.2$]{
\includegraphics[scale=0.52]{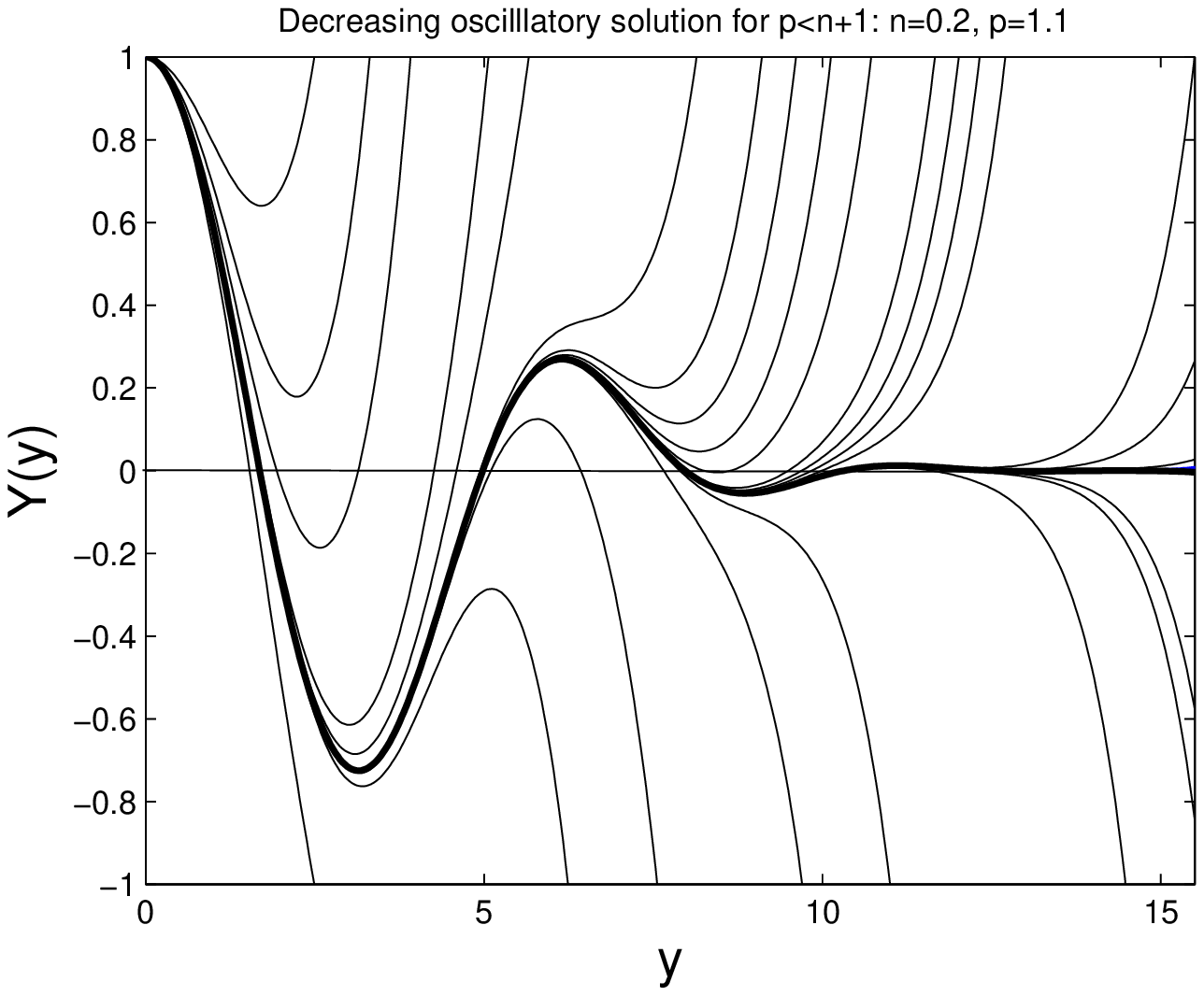}
}
 \vskip -.2cm
\caption{\rm\small Examples of oscillatory patterns of the
``linearized" ODE (\ref{BB1}).}
 \label{FG65}
\end{figure}

\noi{\bf Quasilinear case $n>0$.}
  In Figure \ref{FG65}(a)--(d), we
show that, for any $p \ge n+1$, the ODE (\ref{BB1}) admits
infinitely oscillatory solutions with increasing amplitude of
oscillations.

Here, (a) shows linear increasing amplitude of oscillations for
$n=0$. It is curious that such behaviour persists in the nonlinear
range $n>0$, $p \ge n+1$, so that $n=0$ is a {\em branching point}
for (\ref{BB1}) from solutions of the linear equation (\ref{BB2}).
 In (b), we show the bounded periodic
solution for the variational case $p=n+1$ (cf. Figure \ref{GP1}),
which generates a $p$-branch of non-periodic patterns for $p>n+1$;
see (c). This suggests  that basic blow-up similarity patterns
$\{f_l(y)\}$ are expected to exist for $p >n+1$ sufficiently close
to $n+1$.
Figure \ref{FG65}(d) shows that for $p<n+1$, the amplitude of
oscillations becomes decreasing, so we expect a {\em single}
P-type profile $f_0(y)$ for $p \in (1,n+1)$.

 \subsection{Basic patterns: first numerical conclusions}

 For convenience, similar to the change  in (\ref{1}),
  we perform the following scaling in
 (\ref{f11E}) for $p \not = n+1$:
 \be
 \label{sc1}
 f = A F, \quad y \mapsto a y, \quad \mbox{where}
  \quad A=f_*=(p-1)^{-\frac 1{p-1}}, \,\,\, a=(p-1)^\b.
   \ee
   Then  $F=F(y)$ solves the equation
    \be
    \label{1N}
     \mbox{$
  -(|F''|^n F'')''- \tilde \b F'y - F + |F|^{p-1} F=0, \quad
  \mbox{with}
  \quad \tilde \b= (p-1)\b= \frac{p-(n+1)}{2(n+2)},
 $}
   \ee
which has the scaled equilibria $F_*= \pm 1$ that are convenient
for numerical experiments for small $n>0$ and $p$ close to $1^+$.

\begin{figure}
 \centering
\includegraphics[scale=0.75]{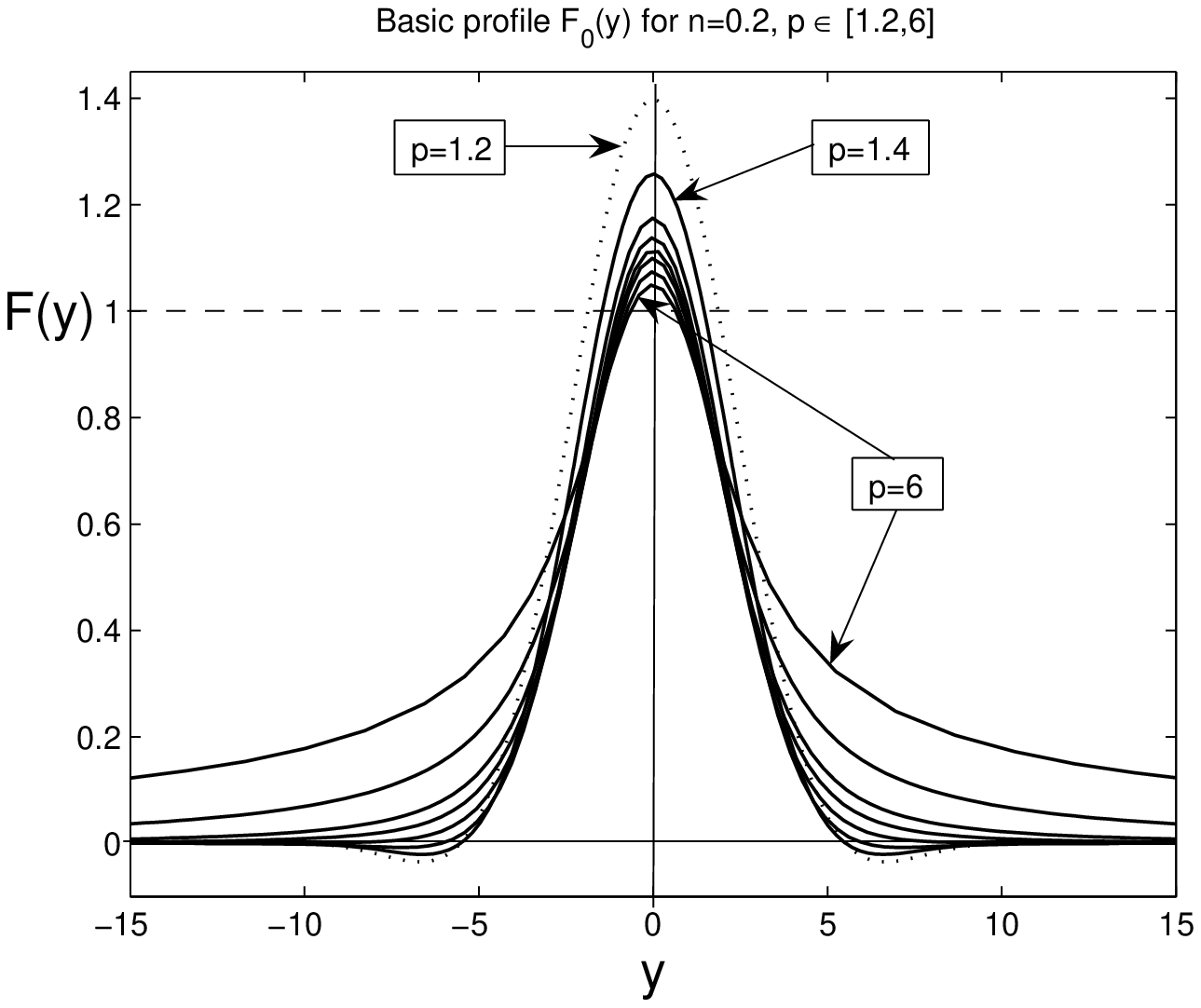}
 \vskip -.4cm
\caption{\rm\small The first basic single point blow-up patterns
$F_0(y)$ of (\ref{1N}) for $n=0.2$ and various $p \in [1.2,6]$.}
   \vskip -.3cm
 \label{FLS1}
\end{figure}

In Figure \ref{FLS1}, we present the first pattern $F_0(y)$ for
$n=0.2$ for  $p =n+1=1.2$ (the dotted line for comparison), 1.4,
  1.8,  2.2, 2.6,  3, 4, and 6.
In particular,  it is clearly seen that, for larger $p$, the
profiles approach the positive asymptotic behaviour (\ref{as1})
for $y \gg 1$ with $C_0>0$, and become strictly positive in $\re$.

\begin{figure}
 \centering
\includegraphics[scale=0.75]{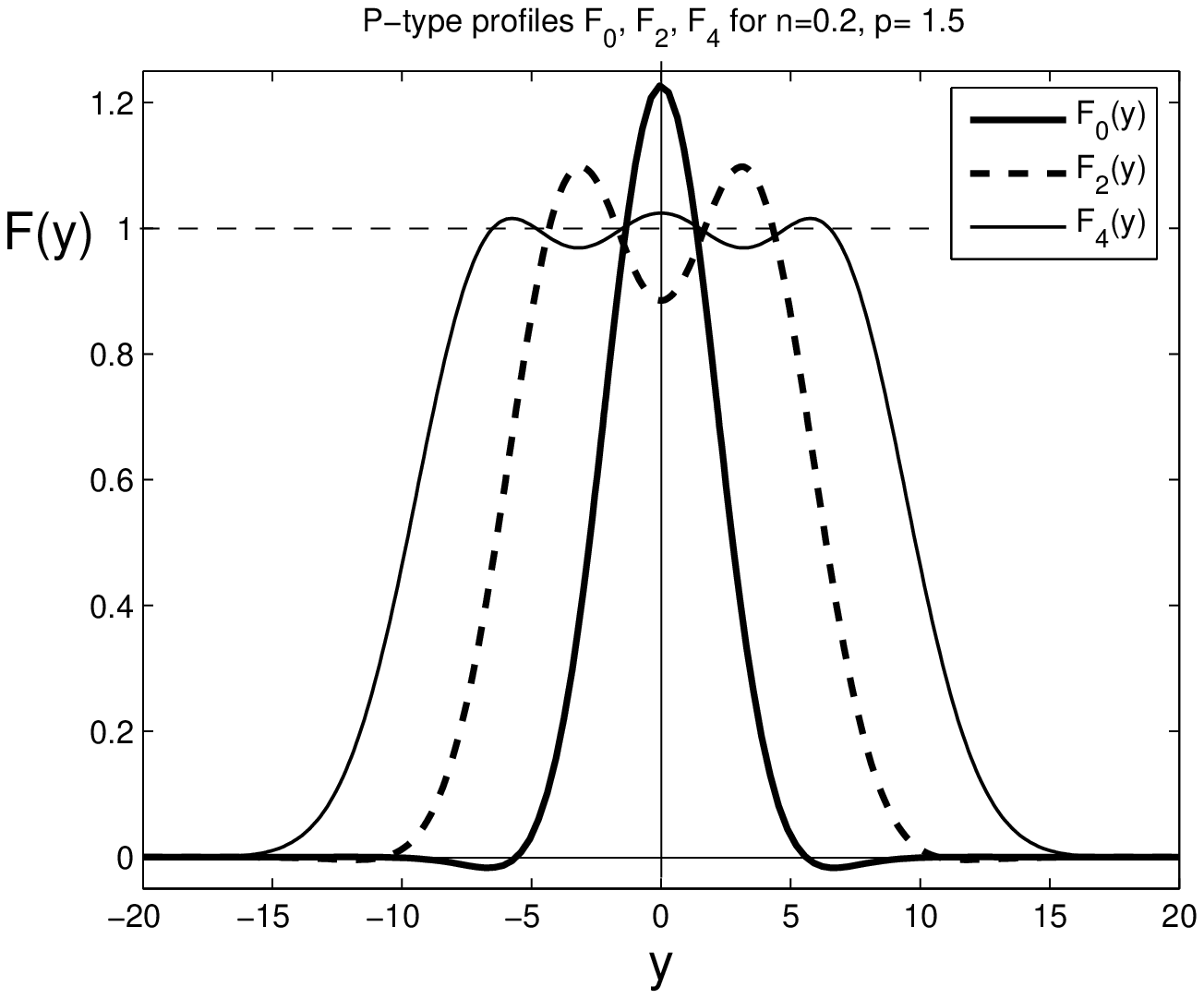}    
 \vskip -.4cm
\caption{\rm\small Three single point blow-up patterns $F_0(y)$,
$F_2(y)$, and $F_4(y)$ of (\ref{1N});
 $n=0.2$,
$p= 1.5$.}
   \vskip -.3cm
 \label{FLS2N}
\end{figure}

\begin{figure}
 \centering
\includegraphics[scale=0.75]{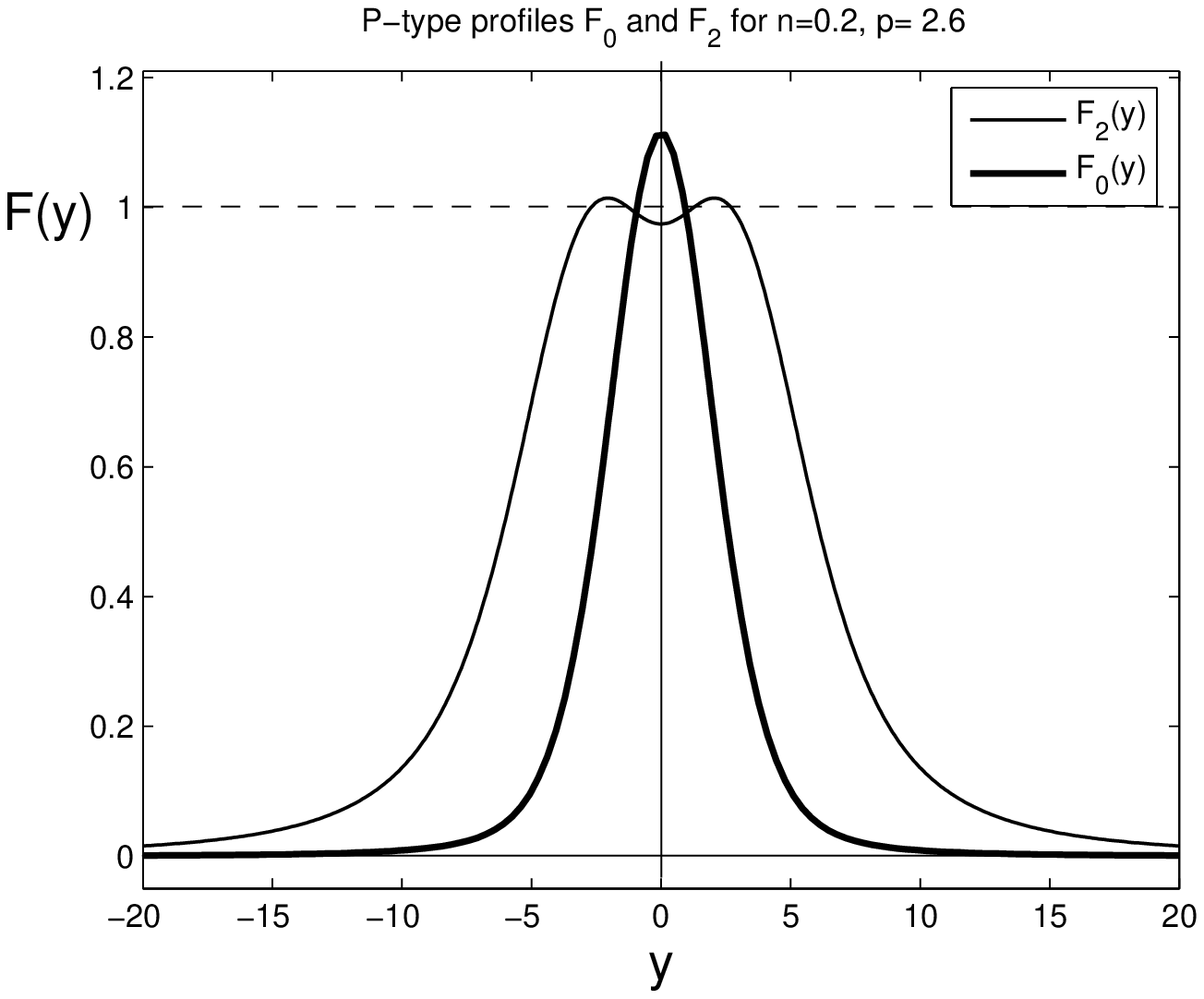}    
 \vskip -.4cm
\caption{\rm\small Two  patterns $F_0(y)$ and $F_2(y)$  of
(\ref{1N});
 $n=0.2$,
$p=2.6$.}
   \vskip -.3cm
 \label{FLS2NN}
\end{figure}

 Figure \ref{FLS2N}
 shows
   the first  $F_0(y)$, the third  $F_2(y)$, and the fifth $F_4(y)$ P-type patterns for
$n=0.2$ and  $p =1.5$. In Figure \ref{FLS2NN}, we demonstrate two
profiles $F_0(y)$ and $F_2(y)$ for the same $n=0.2$ and $p=2.6$.
In the last case, the next even profile $F_4(y)$ was not detected
numerically, and this nonexistence will be confirmed later by the
$p$-branching approach.



 \subsection{Q-type profiles}

 It follows from (\ref{f11E}) that an oscillatory expansion can be
 started at any finite point $y=y_0>0$, at which
  \be
  \label{vv1}
  f(y_0)=f_* \quad \Longrightarrow \quad f(y) \equiv f_* \,\,\,
  \mbox{for}
  \,\,\, y \le y_0.
   \ee
This yields the so-called Q-type solutions; see  classification in
\cite{GPos, BuGa}. Then setting
 $$
 y \mapsto y_0+y, \quad y>0,
  $$
  we again arrive asymptotically  at a linearized equation similar to
  (\ref{BB1}),
   \be
   \label{g11NN}
   \tilde {\bf B}_n(Y) \equiv -(|Y''|^n Y'')''- \l_0 Y'=0, \quad \l_0= \b
   y_0>0,
    \ee
  so on integration we obtain the TW equation (\ref{le2}), where
  the constant $\l_0$ is scaled out. Therefore, we use the change
  (\ref{le3}) to get the oscillatory ODE (\ref{le4}). According to
  (\ref{as55}), this gives a 2D asymptotic family to be matched
  with the bundle (\ref{as1}) at infinity.

\begin{figure}
 \centering
\includegraphics[scale=0.75]{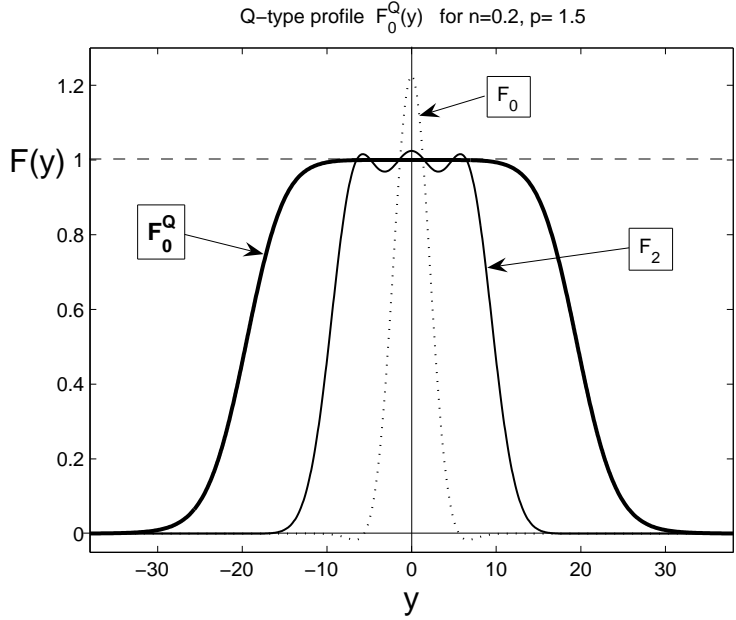}    
 \vskip -.4cm
\caption{\rm\small The first Q-type solution $F_0^Q(y)$ (the
boldface line)  of the ODE (\ref{1N});
 $n=0.2$,
$p= 1.5$.}
   \vskip -.3cm
 \label{FLS2Nq}
\end{figure}

  Analytically, as well as numerically, the problem of existence
  of a countable subset of such Q-type similarity profiles is more
  difficult.
 Figure \ref{FLS2Nq}
 shows
   the first Q-type profile  $F_0^Q(y)$  for
$n=0.2$ for  $p =1.5$.
  The convergence here is  slower and we do
not succeed in getting other Q profiles.
 For the sake of comparison, we also
   present here P solutions $F_0(y)$ and $F_2(y)$.

\subsection{On branching of solutions from variational critical points}

Consider the ODE (\ref{1N}) from the point of view of a
perturbation approach.
 For
$p=n+1$, i.e., for $\b=0$, the ODE has been studied in Section
\ref{SectS}. Setting $\e=p-(n+1)$ and
 assuming that $|\e|>0$ is sufficiently small, we write
  (\ref{f11E}) in the form
  \be
  \label{2.11}
   \mbox{$
 {\bf F}(f) \equiv - (|f''|^n f'')''  -
\frac 1{n} \, f + |f|^{n}f= \frac {\e}{2(n+2)(p-1)} \, y f' -
\frac {\e}{n(n+\e)} \, f + |f|^n f(1-|f|^\e).
 $}
 \ee
 On the right-hand side,
  the key perturbation term satisfies
 \be
 \label{per66}
 g(f,\e)=|f|^n f(1-|f|^\e) \to 0 \quad \mbox{as} \quad \e \to 0
  \ee
   and
  is at least continuous  at the singular point $f=0$, $\e=0$.
 On the left-hand side of (\ref{2.11}), we have the
 variational operator from Section \ref{SectS}.
 Thus, the
 non-perturbed problem (\ref{2.11}) for $\e=0$ admits families of solutions
 described in Section \ref{SectS} (see also extra details on other
 families of solutions in
 \cite[\S~5]{GMPSob} and \cite[\S~5,6]{GMPSobI}). Therefore, classic perturbation and branching
 theory \cite{BufTol, VainbergTr, KrasZ} suggests (and does not prove since $g \not \in C^1$
 at $(0,0)$) that, under natural hypotheses, the variational
 problem for $\e=0$ generates a countable family of $p$-branches,
 which can be extended for some sufficiently small $|\e|>0$.

The analysis of bifurcation, branching, and continuous extensions
should be  performed for the equivalent to (\ref{2.11}) integral
equation with Hammerstein compact operators; see typical examples
in \cite{BGW1, GHUni, GW2}, where similar perturbation problems
for blow-up and global patterns were investigated.
 Namely, denoting by ${\bf D}=(D_y^2)^{-1}$ the inverse of $D^2_y$ in a
 sufficiently large interval $(-R,R)$ (for $p \le n+1$) and in
 $\re$ (for $p>n+1$), we write (\ref{2.11}) as follows:
  \be
  \label{ss1}
   \begin{matrix}
  f= {\bf A}(f,\e) \equiv {\bf D}\big[\big|{\bf D} h \big|^{-\frac
  n{n+1}} {\bf D} h \big], \quad \mbox{where}\qquad\quad
  \ssk\ssk\ssk\\
   h(f,\e)=\frac {\e}{2(n+2)(p-1)} \, y f' -
\frac {\e}{n(n+\e)} \, f + |f|^n f(1-|f|^\e) -
 \big(-\frac 1{n} \, f + |f|^{n}f \big).\qquad\quad
\end{matrix}
 \ee
 For $\e=0$, this gives the integral equation with a potential
 operator, which is equivalent to the differential one studied
 before.

As customary, parameter $p$-branches of solutions of integral
equations with compact operators are fully extensible and can end
up either at a singularity point or at another bifurcation values;
see  \cite{BufTol} for a modern sounding of such  results. An
efficient way to prove  branching is using degree-index theory,
which
 establishes branching from an
isolated solution\footnote{Proving that a given solution is
isolated is also a difficult problem, especially for $p$-Laplacian
operators.}, say, the first one $f_0$ for simplicity, from the
branching point $\e=0$ provided its index  satisfies
\cite[p.~353]{KrasZ}
 \be
 \label{ind1}
\g
\not = 0.
  \ee
  For differentiable operators, where the index ${\rm ind}(f_0,I-{\bf A}'(f_0,0))$ is equal to the
rotation of the vector field $I-{\bf A}'(f_0,0)$, there are
special techniques for its calculations; see
\cite[Ch.~20,\,24]{KrasZ}.

   However, in view of non-potential structure in (\ref{ss1}), and
   of a non-sufficient (in the usual sense) differentiability of (\ref{per66})
   at $(0,0)$, it is convenient to use other
alternatives of bifurcation-branching theory {\em without} direct
differentiability hypotheses that assume sufficient regularity of
the perturbations; cf. \cite[\S~28]{Deim}.
 Namely,
in view of our difficulties with the differentiability, using
 Theorem 28.1
 in \cite[p.~381]{Deim}
  replaces  differentiability  by a slightly weaker control of
  smallness of
  nonlinear terms in a neighbourhood of $(0,0)$.
  As usual, the key principle of branching
 is that it occurs
  if the corresponding
 eigenvalue has odd multiplicity
 (the even multiplicity case needs an additional treatment, which is also a
   routine procedure not to be treated here);
 see further comments below.
 Justification of
branching phenomena for such quasilinear degenerate $p$-Laplacian
operators (including also questions of compactness) in the present
problem needs further deeper analysis and more involved functional
topology/constraints. Therefore, most of further analytical
conclusions remain formal and are open problems.
 Nevertheless, it turns out and will be checked numerically, the
 predicted branching behaviour from $p=n+1$ actually occurs, and
 thus becomes a key tool of the proposed study of non-variational
 problems at hand.

Thus, anyway, for actual
 applications, one needs to know the
spectrum of
the self-adjoint operator ${\bf F}'(f_0)$.
 We have that
(\ref{ss1}) contains nonlinearities that are hardly differentiable
at 0, so these applications
lead to  difficult technical problems, where
 branching analysis from
oscillatory profiles $f_0(y)$ with interior transversal zeros
needs taking into account more complicated ``functional topology".
 Anyway, continuing the application, it is worth mentioning that
 $$
  \l=0 \quad (\mbox{with the eigenfunction} \,\,\, \psi_0 \sim
  f_0')
  $$
  is always an eigenvalue of ${\bf F}'(f_0)$.
   Actually, this corresponds to the invariance
  of the original PDE (\ref{1.5}) relative to  the group of translations
   with the infinitesimal generator
  $D_y$. In particular, this implies that $\l=1$ is an eigenvalue of ${\bf A}'(f_0,0)$,
  so this corresponds to the {\em critical case}, where computing
  of the index is more difficult and is  performed as in
  \cite[\S~24]{KrasZ}.
   It is more important that
 $$
  \l=1 \quad (\mbox{with the eigenfunction} \,\,\, \psi_1 \sim
  f_0)
  $$
 is also an eigenvalue of ${\bf F}'(f_0)$ (this is associated with the
 group of translations
  with the generator $D_t$).
In other words, the index condition (\ref{ind1}) or smallness
assumption via \cite[Th.~ 28.1]{Deim}  need special additional
treatment associated with spectral properties of the linearized
operator ${\bf F}'(f_0)$. In this connection, we conjecture that
 branching at $p=n+1$
 is valid in appropriate functional setting, and
 there exist continuous $\e$-curves, from any profiles from the family of
basic patterns $\{f_l, \, l \ge 0\}$ constructed in Section
\ref{Sect54}.

Finally, let us note
  another convenient (but not that efficient)
way
to use Schauder's Theorem  applied to (\ref{ss1}) to get solutions
of (\ref{2.11}) for small $\e>0$ and to trace out $p$-branches of
the suitable profiles. This approach effectively applies in the
case of porous medium operators as in (\ref{pp1}),
 \cite[\S~6]{GPME}.

\subsection{Numerical construction of $p$-branches}

    We
recall that (\ref{1N}) is not  variational, and we are going to
use  a certain continuity feature concerning the limit $p \to
n+1$.

\begin{figure}
 \centering \subfigure[$p_0$-branch]{
\includegraphics[scale=0.52]{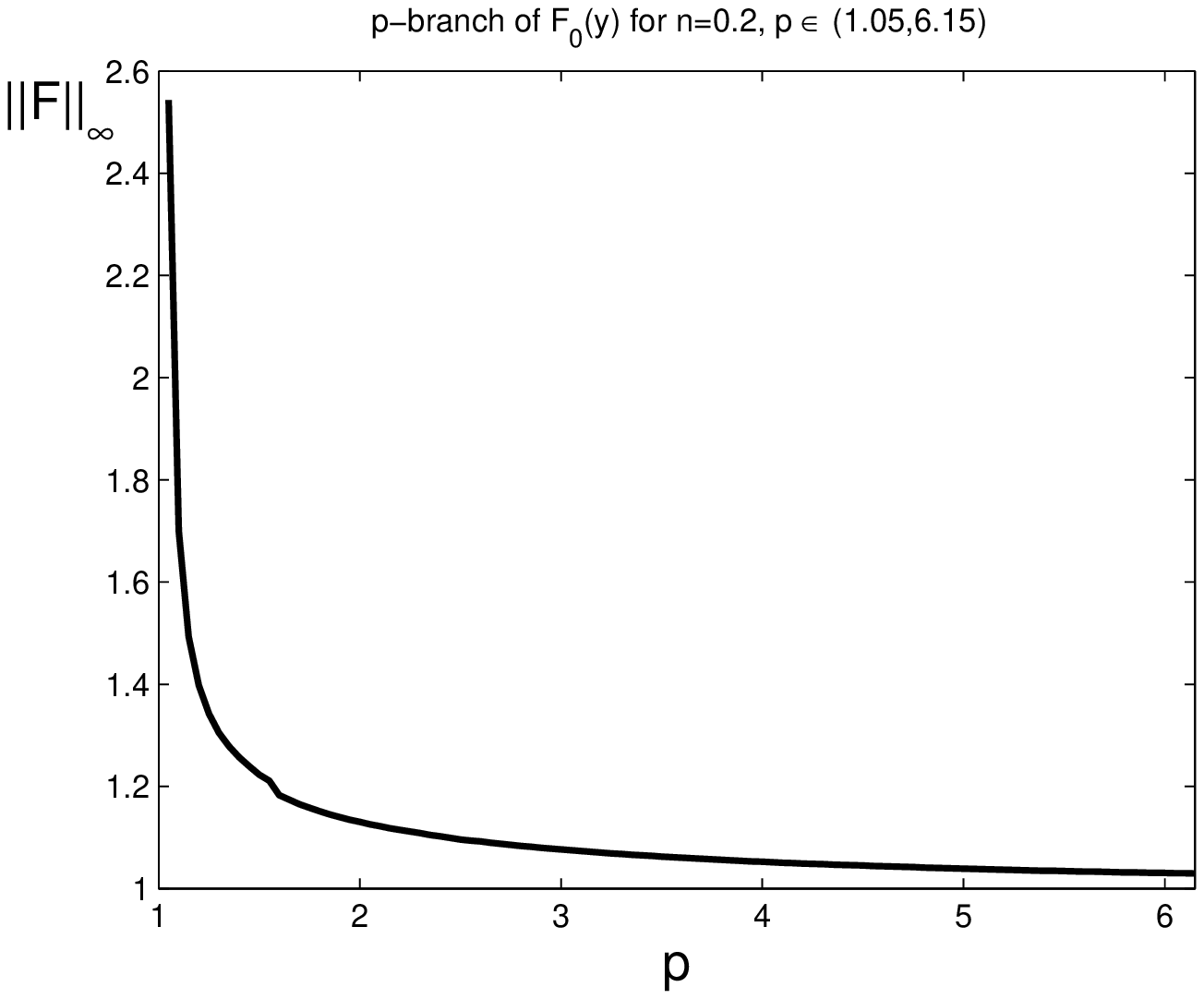}                
} \centering \subfigure[$F_0$ profiles]{
\includegraphics[scale=0.52]{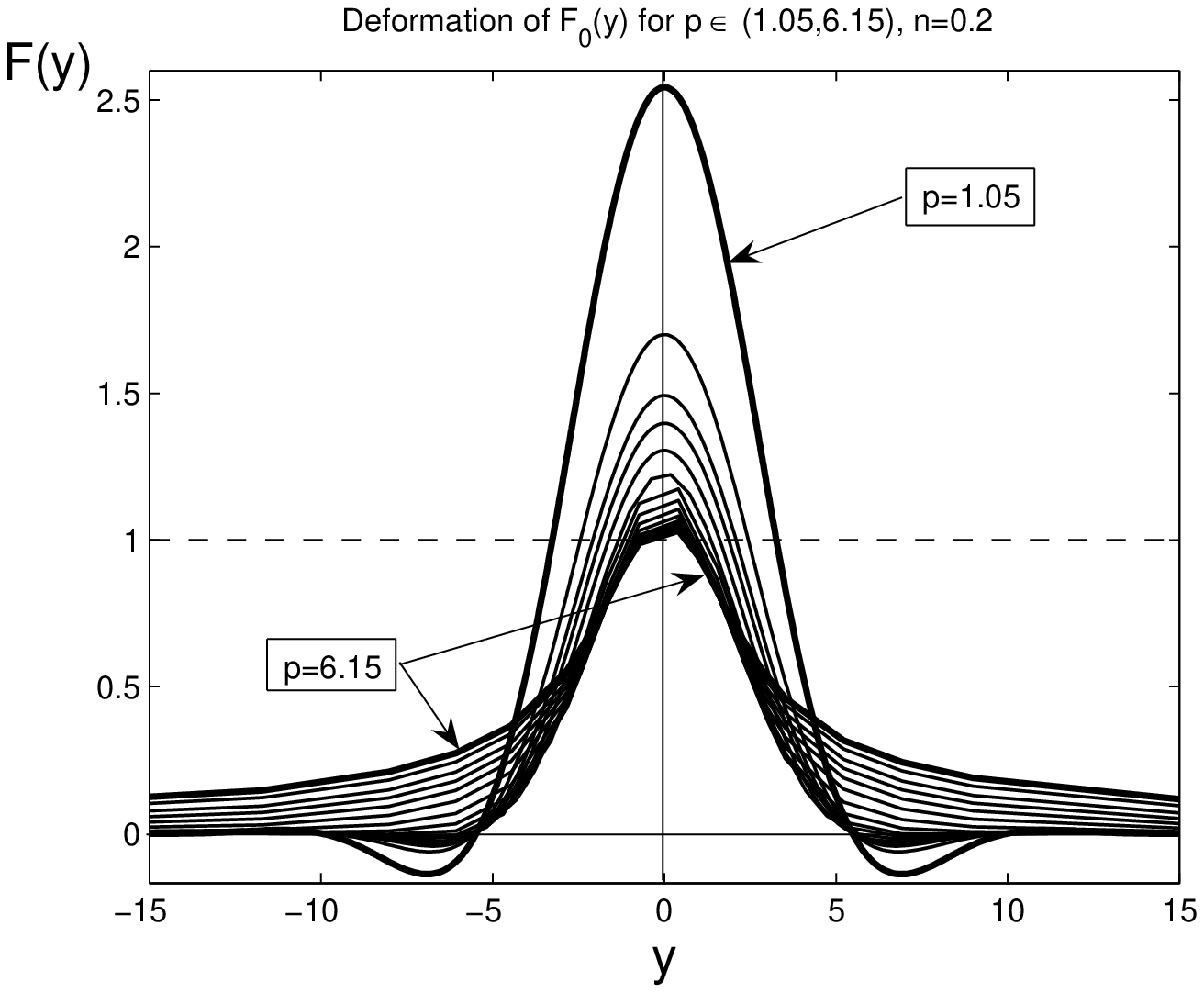}                     
}
 \vskip -.2cm
\caption{\rm\small  The first $p_0$-branch of solutions $F_0(y)$
of equation (\ref{1N}) for $n=0.2$, $p \in (1.05,6.15)$ (a);
corresponding deformation of $F_0$ (b).}
 \label{Fp0p}
\end{figure}

 The basic $p_0$-branch of $F_0(y)$ of the simplest shape (as well as
the $p_1$-one for $l=1$) exists for all $p>1$. In Figure
\ref{Fp0p}, we show the first $p$-branch of $F_0$ in (a) and the
deformation of the profiles $F_0(y)$ in  (b), for $n=0.2$ and $p
\in (1.05,6.15)$.
 This branch is  extended to the global blow-up case $p<n+1$,
 to be  discussed next in Section \ref{SectHS}.
 We
expect that, as usual in blow-up analysis, this {\em first} $p_0$
branch is composed from structurally stable solutions and hence
represents the generic  blow-up behaviour for the parabolic PDE
(\ref{1.5}).

\begin{figure}
 \centering
 {
\includegraphics[scale=0.65]{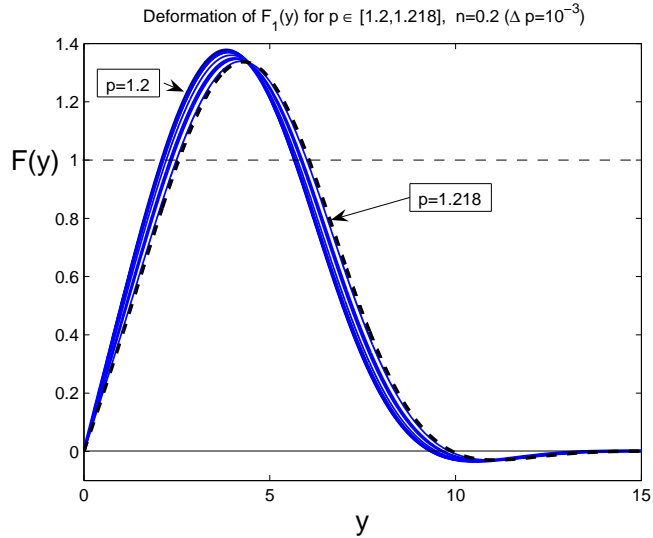}
}
 \vskip -.3cm
\caption{\rm\small  Deformation of the dipole $F_1(y)$ of equation
(\ref{1N}) for $p \in [1.2,1.218]$, $n=0.2$.
}
 \label{Fp1NN}
\end{figure}

Deformation with $p$ of $F_1(y)$ (the  part for $y >0$ is shown
only) for $p$ slightly above the variational exponent $p=1+n=1.2$
for $n=0.2$
 is presented in Figure \ref{Fp1NN}.
 Further extension of this branch beyond $p=1.218$  leads to
 strong numerical
instabilities that possibly reflects the actual nonexistence of
such solutions far away from $p=n+1$.

\ssk

\begin{figure}
 \centering \subfigure[$p$-branch]{
\includegraphics[scale=0.52]{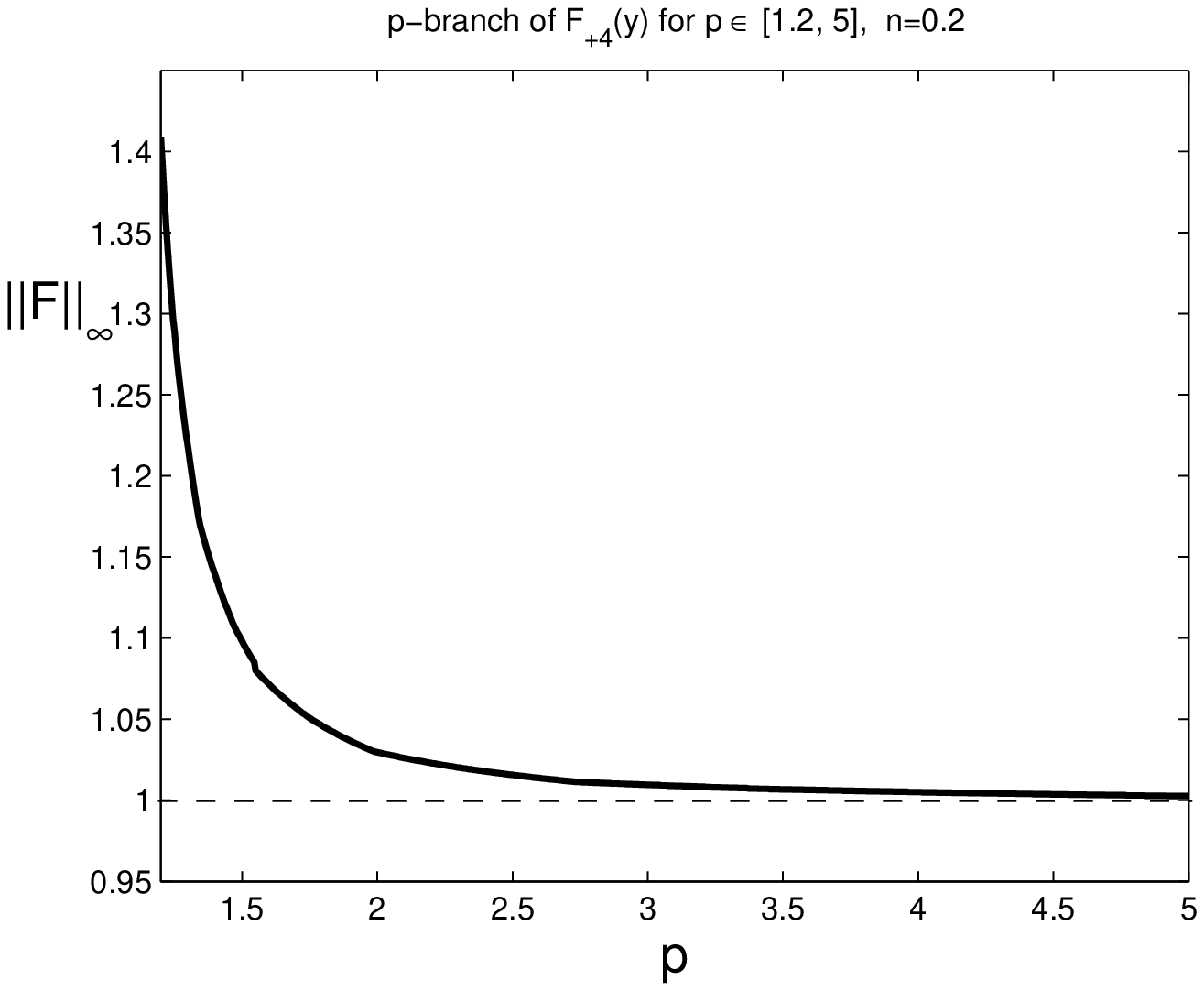}
} \centering \subfigure[$F_{+4}$ profiles]{
\includegraphics[scale=0.52]{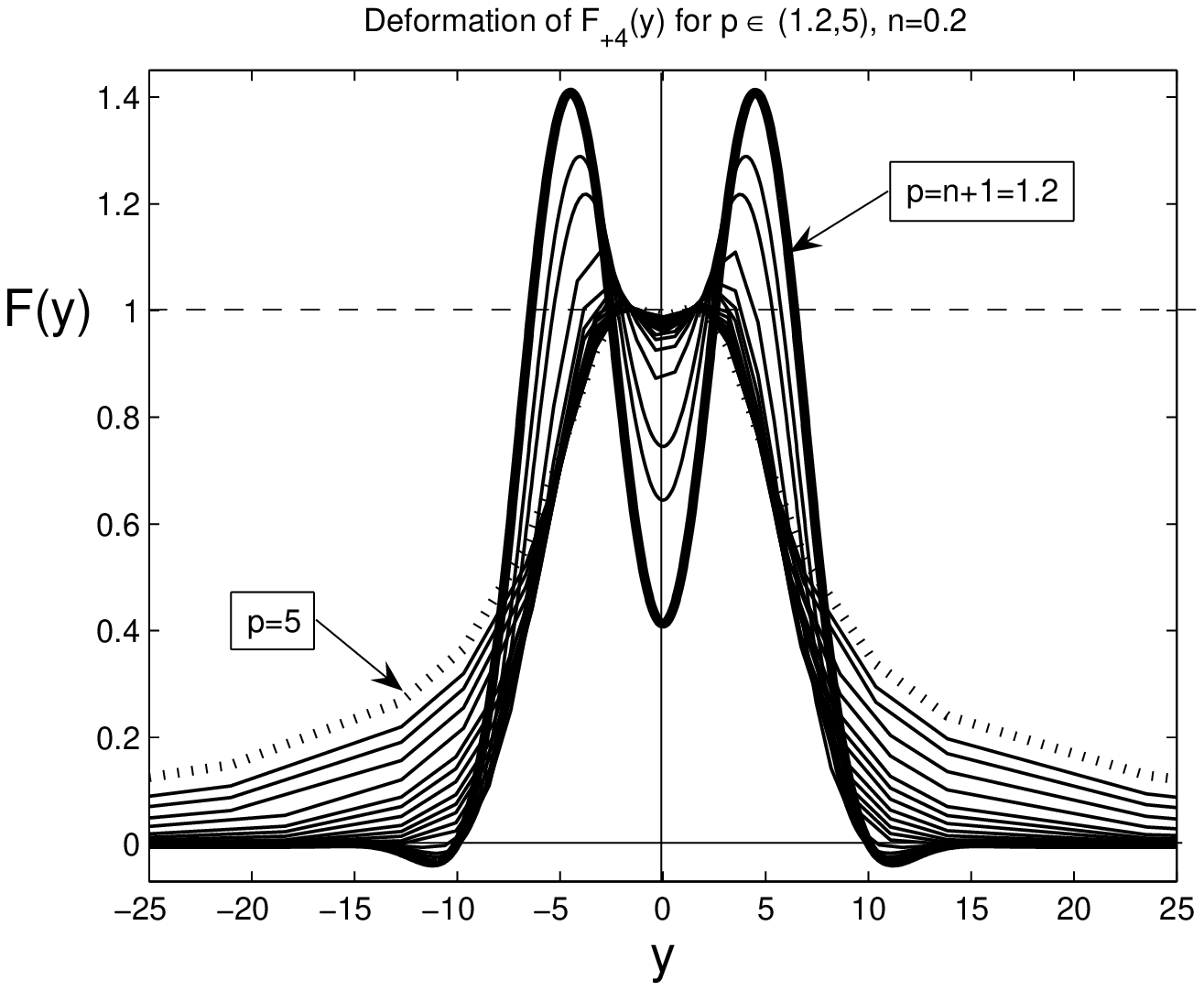}                  
}
 \vskip -.2cm
\caption{\rm\small  The  $p$-branch of solutions $F_{+4}(y)$ of
equation (\ref{1N}) for $n=0.2$ (a); corresponding deformation of
$F_{+4}(y)$ (b).}
 \label{Fp1SS}
\end{figure}

Concerning other, more complicated profiles for $p=n+1$, such as
$F_{+4}(y)$ and others containing structures shown in Figures
\ref{G6}, \ref{FF4p}, and \ref{G8}, numerical results suggest that
some of them  cannot be extended for $p>n+1$ (then any version of
(\ref{ind1}) is not valid).
 Nevertheless, for $F_{+4}$ this is not the case; see Figure
 \ref{Fp1SS}, where the extension is shown to exist for all
 $p>n+1$ and that
  \be
  \label{ff1}
  \|F_{+4}\|_\infty \to 1^+ \quad \mbox{as} \quad p \to +\infty.
   \ee
   It seems that all the global $p$-branches satisfy (\ref{ff1});
   cf. an analogous result in \cite{GHUni} for global similarity solutions.
We expect that a similar $p$-branch of $F_{+4}$ is originated at a
saddle-node bifurcation for some  $p_* \in (1,n+1)$, at which it
appears together with the $p$-branch of the profiles $F_{+2,2+2}$;
see further comments below.


The $p$-branches can connect various profiles, with rather obscure
understanding of possible geometry of such branches and their
saddle-node bifurcation (turning) points.
 For $p=n+1$, the questions on connections with respect to regularization
 parameters as in (\ref{4.1})
  are addressed in \cite{GMPSobII} posing problems of
  homotopy classification of patterns in variational problems and
  approximate ``Sturm's Index" of solutions.

\section{On global blow-up similarity profiles for $p \in (1,n+1)$}
 \label{SectHS}

\subsection{Local oscillatory behaviour close to interfaces
remains the same}

Indeed, the ODE (\ref{f11E}) now reads for $f \approx 0$ as
 $$
 -(|f''|^n f'')'' - \b y f'+...=0 \quad (\b<0),
 $$
 and reflecting near interface $y \mapsto y-y_0$, on integration for small $y >0$, we have
  $$
(|f''|^n f'')' =  \b y_0 f+... \, .
 $$
 This is precisely (\ref{le2}) with $\l =\b y_0 <0$ to be
 reduced to $-1$ by scaling.
 Hence, for $p \in(1,n+1)$, the similarity profiles are
 equally oscillatory near interfaces as for $p=n+1$.
Therefore, according to (\ref{as55}), this local 2D asymptotic
family looks sufficient to be matched with two symmetry boundary
conditions (\ref{BCs})
 (or (\ref{as9}))
at the origin, though the proof of existence remains open.

  \subsection{On similarity profiles and $p$-branches}

\begin{figure}
 \centering
\includegraphics[scale=0.75]{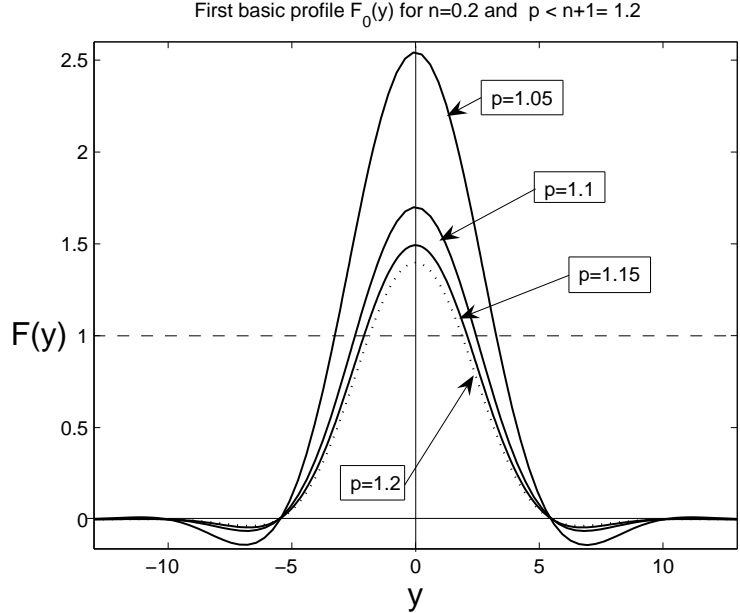}
 \vskip -.4cm
\caption{\rm\small Deformation of the first basic profile $F_0(y)$ of the ODE
(\ref{1N}) for $n=0.2$ and $p=1.2$, 1.15, 1.1, 1.05.}
   \vskip -.3cm
 \label{FHS1p}
\end{figure}

  For $1<p<n+1$, the rescaled ODE (\ref{1N}) is more difficult to solve
  numerically than for $p \ge n+1$.
  Figure \ref{FHS1p} shows deformation of similarity profiles
  $F_0(y)$ for $n=0.2$ and $p \in [1.05,1.2]$.
  We observe an easy visible growth of solutions as $p \to 1^-$.
 Structurally, the first basic profile $F_0(y)$ remains of a similar
 geometric shape as in the variational case $p=n+1=1.2$.

\begin{figure}
 \centering \subfigure[$p_0$-branch]{
\includegraphics[scale=0.52]{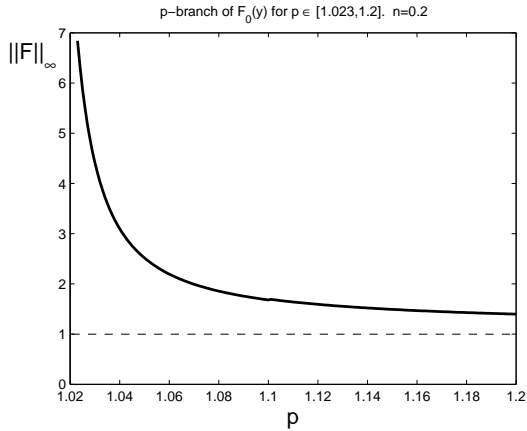}
} \centering \subfigure[$F_0$ profiles]{
\includegraphics[scale=0.52]{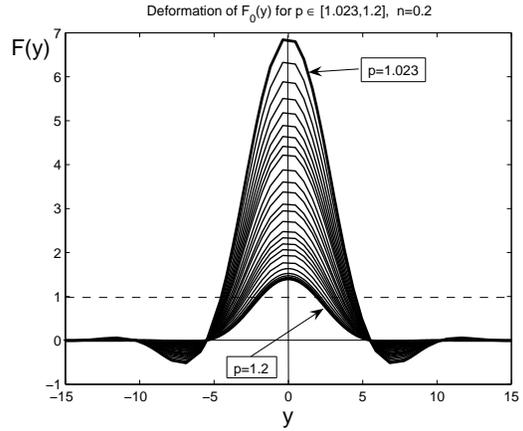}  
}
 \vskip -.2cm
\caption{\rm\small  The first $p$-branch of solutions $F_0(y)$ of
equation (\ref{1N}) for $n=0.2$ (a); corresponding deformation of
$F_0$ (b).}
 \label{FHSN}
\end{figure}

A part of the corresponding $p$-branch of profiles $F_0$ was shown
earlier in Figure \ref{Fp0p}(a). More detailed and sharp results
are presented in Figure \ref{FHSN} for $p \in [1.023,1.2]$,
$n=0.2$, where we used branching from the variational profile for
$p=1.2$ (with the step size $\D p=-10^{-3}$). From (a), we
definitely observe that this $p_0$-branch is going to blow-up as
$p \to 1^-$, as suggested before. Note   that the $p_0$-branch is
expected to consist of asymptotically (structurally) stable
blow-up profiles $F_0(y)$, but we cannot prove this even in the
linearized approximation. The linearized operator is a difficult
non-self-adjoint one with unknown spectrum and proper functional
setting.

\begin{figure}
 \centering \subfigure[$p_1$-branch]{
\includegraphics[scale=0.52]{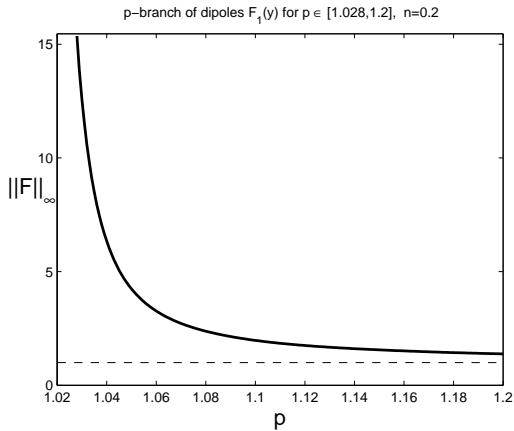}
} \centering \subfigure[$F_1$ profiles]{
\includegraphics[scale=0.52]{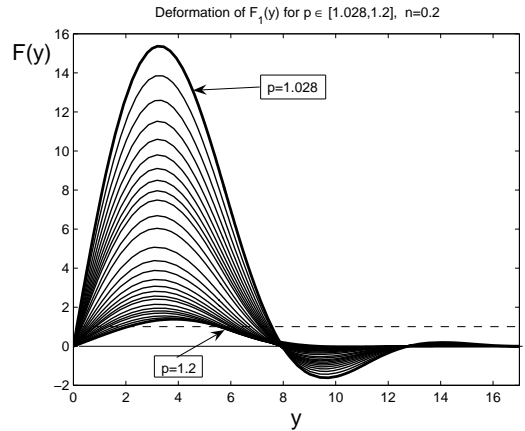}  
}
 \vskip -.2cm
\caption{\rm\small  The second $p$-branch of dipole
 profiles  $F_1(y)$ of
equation (\ref{1N}) for $n=0.2$ (a); the corresponding deformation
of $F_1(y)$ (b).}
 \label{FHSN1}
\end{figure}

The next $p$-branch of dipole-like profiles $F_1(y)$ is shown in
Figure \ref{FHSN1}(a), together with the deformation (b) of  the
functions $F_1(y)$. It is seen that this $p$-branch is global and
blows up as $p \to 1^-$.

\begin{figure}
 \centering \subfigure[$p_2$-branch]{
\includegraphics[scale=0.52]{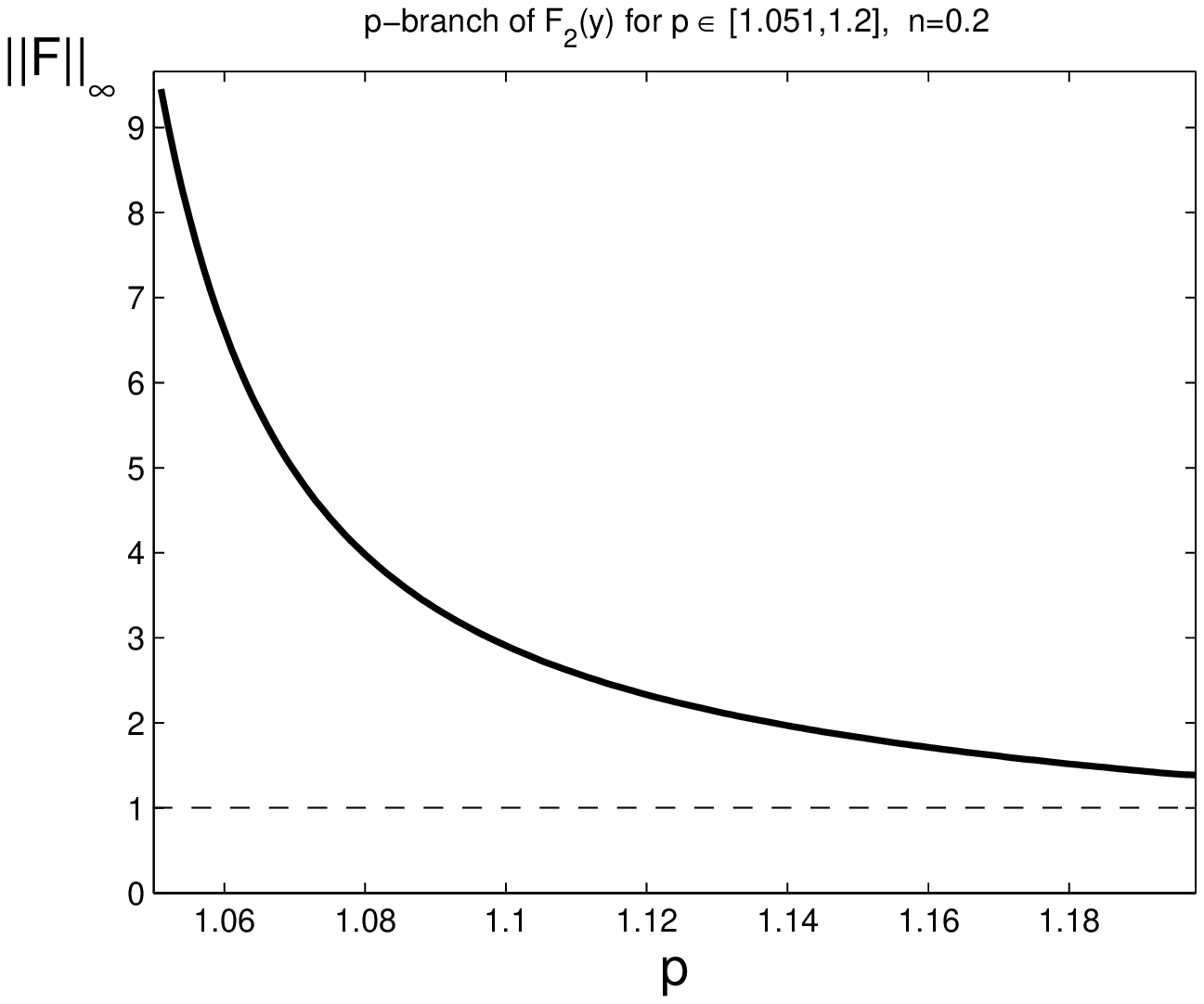}
} \centering \subfigure[$F_2$ profiles]{
\includegraphics[scale=0.52]{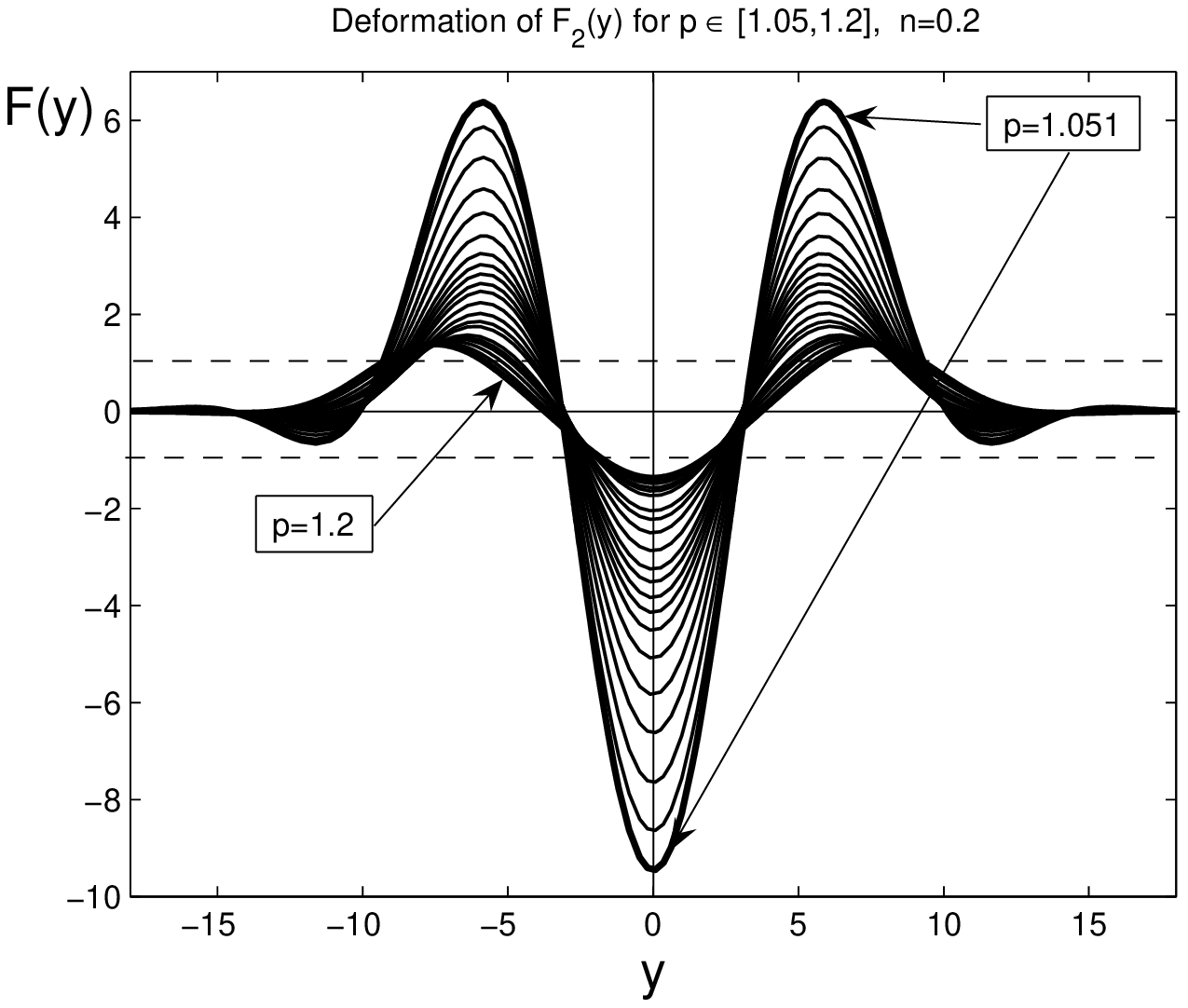}  
}
 \vskip -.2cm
\caption{\rm\small  The third $p$-branch of
 profiles  $F_2(y)$ of
equation (\ref{1N}) for $n=0.2$ (a); the corresponding deformation
of $F_2 (y)$ (b).}
 \label{FF2}
\end{figure}

 We claim that the $p$-branches of the basic similarity profiles
$\{f_l(y)\}$ ({\em q.v.} (\ref{sc1})) are extended up to $p=1^-$,
with a blow-up behaviour as in
 (\ref{g1}). We expect that these
$p$-branches can be connected as $n \to 0$ with those predicted by
the linear problem with patterns (\ref{u11}). We refer to
\cite[p.~1090]{GW2} for an example of such an analysis. For
instance, in Figure \ref{FF2}, we present the $p$-branch and the
corresponding deformation of the third basic profile $F_2(y)$,
which for $p=n+1$ is given in Figure \ref{G4}(c) (with the
opposite sign).

\begin{figure}
 \centering
\includegraphics[scale=0.65]{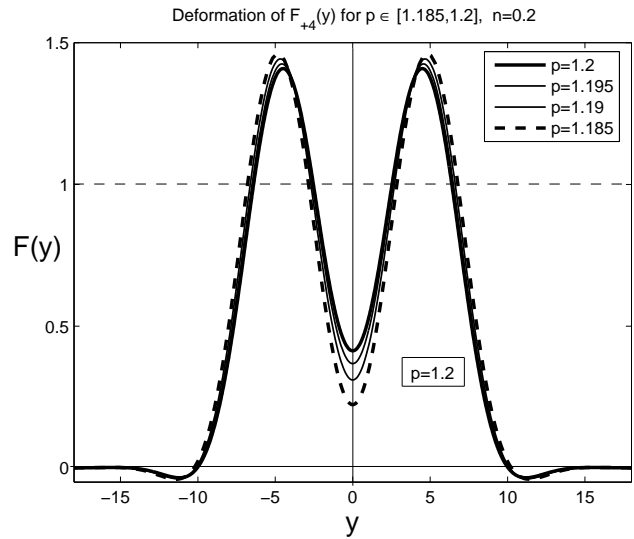}
 \vskip -.4cm
\caption{\rm\small The deformation of  $F_{+4}(y)$ of the ODE
(\ref{1N}) for $n=0.2$ and $p=1.2$, 1.195, 1.19, 1.185.}
   \vskip -.3cm
 \label{F44p}
\end{figure}

Also,  a principal fact of existence of the $p$-branch of the
non-basic profiles is explained in Figure \ref{F44p}, where a
local $p$-branch of $F_{+4}(y)$ (see Figure \ref{FF4p} for
$p=n+1$) is shown to exist for $p<n+1=1.2$ for $n=0.2$.

It is key that this branch cannot be extended for all $1<p<n+1$.
We expect that, as $p<n+1$ decreases, the $p$-branch of $F_{+4}$
meets the $p$-branch of the  ``geometrically similar" profile
$F_{+2,2,+2}$ shown in Figure \ref{G6}(a) (both have two dominant
maxima and a single minimum in between) in a turning {\em
saddle-node} bifurcation point (another branch scenario is also
possible,  \cite[\S~7]{GPME}). Such scenarios were detected in
variational problems; see \cite{GMPSobII}. In the present
non-variational case, both analytic and even a reliable numerical
 description of such bifurcations become much more difficult and still obscure.

\ssk

\ssk

{\bf Acknowledgements.}  The author would like to thank
S.I.~Pohozaev for a discussion on
genus theory.




\end{document}